\documentclass[a4paper,leqno,10pt]{article}
\usepackage{amssymb,amsfonts,amsmath,amsthm}
\usepackage{epsfig,epsf,curves}
\usepackage{bbm,epic,eepic}
\usepackage[mathscr]{eucal}
\theoremstyle{remark}

\begin{document}

\sloppy

%------------------------- Definitions --------------------------

\def\Mat#1#2#3#4{\left(\!\!\!\begin{array}{cc}{#1}&{#2}\\{#3}&{#4}\\ \end{array}\!\!\!\right)}

\def\tAA{{\Bbb A}}
\def\CC{{\Bbb C}}
\def\HH{{\Bbb H}}
\def\NN{{\Bbb N}}
\def\QQ{{\Bbb Q}}
\def\RR{{\Bbb R}}
\def\ZZ{{\Bbb Z}}
\def\FF{{\Bbb F}}
\def\SS{{\Bbb S}}
\def\GG{{\Bbb G}}
\def\PP{{\Bbb P}}
\def\LL{{\Bbb L}}

%----------------------------------------------------------------

\title{On mod $p^c$ transfer and applications}

\author{J. Mahnkopf, U. of Vienna}

\maketitle

{\small 
We study a mod $p^c$ analog of the notion of transfer for automorphic forms. Instead of existence of eigenforms, such transfers yield 
congruences between eigenforms but, like transfers, we show that they can be established by a comparison of trace formulas. This rests on the properties of 
mod $p^c$ reduced multiplicities which count congruences between eigenforms. As an application we construct finite slope $p$-adic {\it continuous} families 
of Siegel eigenforms using a comparison of trace formulas.
}

%describe an elementary comparison of trace formulas which yields the existence of finite slope $p$-adic {\it continuous} families of Siegel eigenforms.
 
\bigskip

{\bf (0.1) } In this article we were motivated by the idea of the {comparison} of trace formulas as a universal principle to relate the existence of different types of automorphic 
representations. If such an 
idea holds true then the comparison of trace formulas should be applicable to a wider class of statements relating different types of automorphic 
representations. Probably the most important statement of this kind is the functoriality principle. This principle yields a transfer from automorphic 
representations on a reductive group $G$ to the set of automorphic representations on another reductive group $G'$. Thus, it relates the existence of automorphic 
representations on different groups and special cases of functoriality have been proven via a comparison of trace formulas.

Another example of such a statement which also is of a general natural is the theory of $p$-adic families of automorphic forms. Given an 
eigenform $f_0$ it predicts the existence of a $p$-adic analytic family of eigenforms passing through $f_0$, i.e. it predicts the existence of 
infinitely many eigenforms which are related to each other and not only to $f_0$ and which are determined by $f$ only modulo a power of $p$ (but which are all on the 
same group $G$). Thus, like the functoriality principle, the theory of $p$-adic families relates the existence of automorphic forms (in varying weights) but it is of 
a different nature. An approach based on a comparison of trace formulas therefore is not obvious but would confirm the idea of the 
comparison of trace formulas as a universal principle.

%: ?????????????? functoriality and $p$-families can be deduced from the same source.

\bigskip

{\bf (0.2) } In this article we describe a comparison of trace formulas which solves the problem of existence of  
{\it continuous} families of finite slope. Since a continuous family is defined by a system of congruences between its eigenforms we need a 
comparison of trace formulas which instead of existence of eigenforms yields congruences between eigenforms. We are led to such a kind of comparison by considering  
a mod $p^c$ analog of the notion of transfer for automorphic forms. Unlike transfers, a mod $p^c$ transfer only yields a mod $p^c$ approximation to the predicted 
eigenform, thus, it yields a {\it congruence} which is satisfied by an eigenform. 
%; as a consequence, 
%continuous families can be constructed by applying mod $p^c$ transfers to the initial eigenform $f_0$ (cf. 0.4 below). 
Like transfers, mod $p^c$ transfers can be established by a comparison of trace formulas, i.e. the existence of a mod $p^c$ transfer follows from certain congruences between 
traces of Hecke operators. 
%if one replaces multiplicities by mod $p^c$ reduced multiplicities.
This is based on the properties of mod $p^c$ reduced multiplicities. Unlike multiplicities, reduced multiplicities count congruences 
between eigenforms, hence, their relation to mod $p^c$ transfers. On the other hand, like multiplicities, they  
can be computed as traces of certain Hecke operators, hence, their relation to the trace formula; cf. (0.3) for more details.

\medskip

In part B we verify the necessary congruences for the group ${\rm GSp}_{2n}$ by comparing the geometric sides of two simple topological trace formulas. 
This essentially comes down to a problem about eigenvalues of certain 
symplectic matrices (cf. section 7, in particular (7.1) Lemma). 

\medskip

We only construct continuous families of eigenforms but our method is very different 
from the ones used in the construction of analytic families. In particular, we do not make use of overconvergent cohomology, $p$-adic Fredholm theory 
or rigid analytic geometry and the proof is of an elementary nature. We hope that the present method also might apply to yield analyticity of families. 

\medskip

Since, here, we compare two trace formulas on the same group, we avoid the deep problems which have to be solved if the trace formula 
is applied to the functoriality principle.

%Formally it might even apply to yield congruences between eigenforms on different groups $G$ and $G'$ but
%this needs to relate harmonic analysis on $G$ and $G'$ which is the main difficulty in proving functoriality.

\medskip

%The proof that we obtain is elementary. The geometric side of the simple and elementary topological trace formula that we use only consists of archimedean components of 
%orbital integrals, i.e. of values of finite dimensional characters and the {\it Weyl character formula} reduces their evaluation to an elementary problem about eigenvalues of 
%certain symplectic matrices (cf. (7.1) Lemma). 

%its main ingredients are
%Lagrange interpolation in the proof of (2.3) Lemma (leading to (2.4) Theorem); the Weyl character formula and Fermat's Little Theorem in the proof of 
%(7.4) Proposition (leading to (7.5) Theorem) and a simple topological trace formula that we use is of an elementary nature (cf. (6.5) Remark). 

%It is typical for the present approach that the result of the mod $p^c$ comparison depends on the dimension of the involved spaces of automorphic forms. In 
%the application to $p$-adic families we therefore need to know in advance that the dimension of the slope subspaces is bounded. This is exactly the first part of the 
%Gouvea-Mazur Conjecture and has been proven, again, by elementary means in [Ma 1].

%This should correspond to congruence relation between orbital integrals. 
%and may be seen as a kind of $mod $p^c$ functoriality.

\bigskip

{\bf (0.3) Mod $p^c$-transfer. } We explain the mod $p^c$ transfer on which our construction of $p$-adic families is based in more detail. Let ${\cal H}$, ${\cal H}'$ be free 
commutative ${\Bbb Z}$-algebras 
in countably many generators (e.g. Hecke algebras attached to reductive groups $G,G'$) and denote by $\hat{\cal H}$ the set of characters $\Theta:\,{\cal H}\rightarrow\bar{\Bbb Q}_p$. 
For any ${\cal H}$-module ${\bf H}$ we denote by ${\cal E}({\bf H})\subseteq\hat{\cal H}$ the set of 
{eigencharacters} occuring in ${\bf H}$, i.e. ${\cal E}({\bf H)}$ consists of all characters $\Theta$ such that the 
corresponding generalized (simultaneous) eigenspace ${\bf H}(\Theta)$ does not vanish. Let $\Phi:\,{\cal H}'\rightarrow{\cal H}$ be an algebra morphism and 
denote by
$$
\Phi^\vee:\,\hat{\cal H}\rightarrow\hat{\cal H}'
$$ 
the dual map; cf. (2.1). Let ${\bf H}$ resp. ${\bf H}'$ be a ${\cal H}$ resp. a ${\cal H}'$-module. The fundamental problem is to 
examine if $\Phi^\vee$ defines a map on {\it eigencharacters} $\Phi^\vee:\,{\cal E}({\bf H})\rightarrow{\cal E}({\bf H}')$. 
(If ${\bf H},{\bf H}'$ are Hecke modules of automorphic forms this means to examine whether $\Phi^\vee$ defines a transfer for automorphic 
representations.)
In this article, with the application to the theory of $p$-adic families in mind, we want to examine whether $\Phi^\vee$ defines a map on 
eigencharacters if we reduce modulo a given power of $p$. More precisely, we want to examine if there is a map
$$
\Psi^{[c]}:\,{\cal E}({\bf H})\rightarrow{\cal E}({\bf H}')\leqno(1)
$$
satisfying $\Psi^{[c]}(\Theta)\equiv\Phi^\vee(\Theta)\pmod{p^c}$ for all $\Theta\in{\cal E}({\bf H})$. We show that the 
existence of such a mod $p^c$ transfer  
$\Psi^{[c]}$ corresponding to $\Phi$ can be established by a comparison of trace formulas if we replace the notion of multiplicity by that of a 
mod $p^c$ reduced multiplicity. The mod $p^c$ reduced multiplicity of $\Theta\in\hat{\cal H}$ is defined as
$$
m_{\bf H}(\Theta,c)=\sum_{\mu\equiv\Theta\pmod{p^c}} {\rm dim}\, {\bf H}(\mu), 
$$
where $\mu$ runs over all characters of ${\cal H}$ which are congruent to $\Theta$ modulo $p^c$. Thus, $m_{\bf H}(\Theta,c)$ counts the number of 
eigencharacters of ${\bf H}$ which are congruent to $\Theta$ mod $p^c$. In particular, if
$$
m_{\bf H}(\Theta,c)=m_{{\bf H}'}(\Phi^\vee(\Theta),c)\leqno(2)
$$
for all $\Theta\in\hat{\cal H}$ then for any eigencharacter $\mu\in{\cal E}({\bf H})$ there is an eigencharacter 
$\mu'\in{\cal E}({\bf H}')$ such that $\mu'\equiv\Phi^\vee(\mu)\pmod{p^{c}}$, i.e. a mod $p^c$ transfer $\Psi^{[c]}$ correspondong to $\Phi$ exists. 
On the other hand, it is crucial that the 
reduced multiplicities in (2) can be computed as traces, i.e. there is an element $e'\in{\cal H}'$ such that
$$
{\rm tr}(\Phi(e')|{\bf H})\equiv m_{\bf H}(\Theta,c) \quad\mbox{and}\quad{\rm tr}\,(e'|{\bf H}')\equiv m_{{\bf H}'}(\Phi^\vee(\Theta),c) 
$$
modulo a "high" power of $p$; cf. (2.3) Lemma and (2.4) Remark. Using this we obtain:

%We will show there is an element 
%${\bf e}\in{\cal H}'$ such that
%$$
%{\bf m}_H(\Theta,c)={\rm tr}\,{\bf e}|{\bf H}\quad\mbox{and}\quad {\bf m}_H(\Theta,c)={\rm tr}\,(\Phi^\vee|{\bf H}')
%$$

\bigskip

{\bf Theorem} (cf. 2.4 Theorem).  {\it Let ${\rm dim}\,{\bf H},\,{\rm dim}\,{\bf H}'\le \frac{M}{2}$ and assume that there is an $s$ such that
$$
{\rm tr}\,(\Phi(T')|{\bf H})\equiv{\rm tr}\,(T'|{\bf H}')\pmod{p^s}\leqno(3)
$$
for all $T'\in{\cal H}'$. Then for any $\Theta\in\hat{\cal H}$ there is $c=c(\Theta)>\frac{s}{M}-(M+2) \log_p M$ such that equation (2) holds.
%$$
%m_{\bf H}(\Theta,c)=m_{{\bf H}'}(\Phi^\vee(\Theta),c)
%$$ 
In particular, a mod $p^{\frac{s}{M}-(M+2)\log_p M}$ transfer corresponding to $\Phi$ exists.
}

\bigskip

{\bf (0.4) Application to $p$-adic families. }  The theory of $p$-adic continuous families is a special case of mod $p^c$ transfers. To explain this, let ${\bf H}_\lambda$ be a 
family of ${\cal H}$-modules indexed by their ``weight'' $\lambda\in {\Bbb Z}^n$. If there are ${\sf a},{\sf b}$ such that a congruence 
$\lambda\equiv\lambda_0\pmod{(p-1)p^m}$ implies that there is a mod $p^{{\sf a}(m+1)+{\sf b}}$-transfer
$$
\Psi_\lambda:\,{\cal E}({\bf H}_{\lambda_0})\rightarrow{\cal E}({\bf H}_\lambda)
$$
corresponding to the identity map $\Phi={\rm id}$, i.e. $\Psi_\lambda(\Theta)\equiv\Theta\pmod{p^{{\sf a}(m+1)+{\sf b}}}$ then the collection of transfers 
$(\Psi_\lambda(\Theta_0))_\lambda$ is a $p$-adic contionuous family passing 
through a given initial eigencharacter $\Theta_0\in{\cal E}({\bf H}_{\lambda_0})$ (cf. (1.9) Proposition). Thus, the existence of $p$-adic continuous families 
follows from a system of congruences of the type in (3) (cf. (3.7) Proposition).

\medskip

%Let ${\bf H}_\lambda$ be a family of ${\cal H}$-modules indexed by their ``weight'' $\lambda\in {\Bbb Z}^n$. If any 
%$\Theta_0\in{\cal E}({\bf H}_{\lambda_0})$ fits in a $p$-adic 
%continuous family $\Theta_\lambda\in{\cal E}({\bf H}_\lambda)$, i.e. there are numbers 
%${\sf a},{\sf b}$ such that $\lambda\equiv \lambda'\pmod{(p-1)p^m}$ implies $\Theta_\lambda\equiv \Theta_{\lambda'}\pmod{p^{{\sf a}(m+1)+{\sf b}}}$, then 
%there is a family of maps
%$$
%\Psi_\lambda:\,{\cal E}({\bf H}_{\lambda_0})\rightarrow{\cal E}({\bf H}_\lambda),\quad(\lambda\in{\Bbb Z}^n)
%$$
%such that $\Psi_\lambda(\Theta)\equiv\Theta\pmod{p^{{\sf a}(m+1)}}$ if $\lambda\equiv\lambda_0\pmod{(p-1)p^m}$. Thus, $\Psi_\lambda$ is a mod 
%$p^{{\sf a}(m+1)}$ transfer corresponding to the identity map $\Phi$. On the other hand, the transfers $\Psi_\lambda$
%already imply the existence of continuous families (cf. (1.9) Proposition). Hence, the existence of $p$-adic continuous families follows from 
%a system of congruences of the type in (3) (cf. (3.7) Proposition). 

%\medskip

In part B we show that the family of slope subspaces of the cohomology of Siegel upper half plane 
satisfies these congruences by comparing trace formulas (cf. (7.5) Theorem). This the main technical work. As a consequence we obtain

\bigskip

{\bf Corollary } (cf. (7.7) Corollary). {\it Any Siegel modular eigenform $f_0$ of slope $\alpha$ fits in a $p$-adic continuos family of eigenforms of slope $\alpha$.

%Moreover, the dimension of the spaces $H^i(\Gamma\backslash {\Bbb H}_{2n},{\cal L}_\lambda)^{\le\beta}$ is locally constant as a function of  $\lambda$}
}

\medskip

We also obtain local constancy of the dimension of the slope spaces.

\bigskip

{\bf (0.5) } Starting with the work of Hida (cf. [H 1], [H 2]) $p$-adic families of automorphic eigenforms have been constructed by several authors; we mention
Hida, Ash-Stevens, Buzzard, Coleman, Emerton, Tilouine, Urban, Harder ... Moreover, there is the work of Koike  [K] who applies an explicit Selberg trace 
formula to prove mod $p^m$ congruences between the traces of Hecke operators on the space of elliptic modular forms for varying weights $k$. Since he only 
considers the Hecke operators $T(p^m)$ at $p$, he can not make statements concerning (existence 
of) eigenforms, hence, there is no kind of transfer or relation between eigenforms in different weights and, consequently, he does not have to set up a 
comparison of trace formulas as we described in part A (and which involves comparing traces of all Hecke operators). 
%and he need not apply any kind of comparison of trace formulas which seems to be unavoidable if one has to considers traces of all Hecke operators (compare the content of our part A). 
%Consequently his trace identity is restricted to ${\rm GL}_2$ whereas we consider the more general case of ${\rm Sp}_n$-Hecke operators acting on slope spaces. 
We also mention work of Urban who $p$-analytically interpolates the traces of Hecke operators for varying weight. Here, Franke's trace formula is applied in 
the construction but as far as we understand in an technical way to reduce from the whole spectrum to the cuspidal spectrum; essentially, his construction 
is based on the work of [A-S], in particular, on their notion of overconvergent cohomology. Thus, like Koike he does not set up a comparison of trace 
formulas (note that Franke's trace formula only has a spectral side and no geometric side, hence, it cannot function in a comparison of trace formulas) 
and his work seems to be very different from ours.

%We also mention work of Urban. Like Koike, this does not seem to set up a comparison of trace formulas to deduce the existence of $p$-adic families. As far as we understand  
%he $p$-analytically interpolate the traces of Hecke operators for varying weight. But e does not deduce this from a trace formula, instead he applies the work of [A-S], in particular their 
%notion of overconvergent 
%cohomology, which in [A-S]leads to an analytic interpolation of the eigenvalues of Hecke operators. In particular, there is no evaluation of
%spectral sides as in our part B. (Franke's trace formula which appears only seems to be used to reduce spectral problems to the cuspidal spectrum; it has no geometric side and, 
%hence, can not relate spectral properties to geometric properties which is the basic characteristic of the Arthur-Selberg or topological trace formulas.) 

%and he does not deduce from trace identities to existence of eigenforms (ware auch unsinn denn [A-S] beweisen mit ihren Methoden, die Urban ubernimmt 

\medskip

Part of this article has been described in [Ma 2,3].

\newpage

\newpage

\centerline{\bf \Large A. Reduced Multiplicities}

\bigskip

\section{Mod $p^c$ reduced Multiplicities}

We define mod $p^c$ reduced multiplicities and we describe their connection to congruences between 
eigenforms.

\bigskip

{\bf (1.1) } We fix a prime $p\in{\Bbb N}$. We denote by $v_p$ the $p$-adic valuation on $\bar{\Bbb Q}_p$ normalized by $v_p(p)=1$. We write 
${\cal O}$ for the ring of integers in $\bar{\Bbb Q}_p$ and we say that $x\equiv y\pmod{p^t}$, $t\in{\Bbb R}$, if $v_p(x-y)\ge t$ ($x,y\in\bar{\Bbb Q}_p$). 
Finally, $\lceil x\rceil$ denotes the smallest integer larger than or equal to $x$ and $\log_p$ is the complex logarithm with base $p$.

\bigskip

{\bf (1.2) } We let ${\cal H}=\bar{\Bbb Q}_p[T_\ell,\,\ell\in I]$ be the polynomial algebra over $\bar{\Bbb Q}_p$ generated by a countable number of 
elements $T_\ell$, $\ell\in I$. We set ${\cal H}_{\cal O}={\cal O}[T_\ell,\,\ell\in I]$, hence, ${\cal H}={\cal H}_{\cal O}\otimes\bar{\Bbb Q}_p$ and ${\cal H}_{\cal O}$
is an order in ${\cal H}$. We denote by $\hat{\cal H}_{{\cal O}}$ the set of all $\bar{\Bbb Q}_p$-algebra characters $\Theta:\,{\cal H}\rightarrow \bar{\Bbb Q}_p$ 
%d.h. $\Theta|_{\bar{\Bbb Q}}={\rm id}$; Beachte: summe von zwei charakteren ist kein charakter mehr, da $\mu+\nu$ nicht mehr multiplikativ ist. D.h. 
%$\hat{\cal H}_{{\cal O}}$ ist kein modul und keine algebra (die menge aller Abbildungen von ${\cal H}_{\cal O}$ nach ${\cal O}$ ist eine ${\cal O}$-algebra)
which are defined over ${\cal O}$, i.e. which satisfy $\Theta({\cal H}_{\cal O})\subseteq{\cal O}$. 
%d.h. $\Theta(t_\ell)\in{\cal O}$ for all $\ell\in I
In the following we will also understand by $\Theta\in\hat{\cal H}_{{\cal O}}$ the induced character
$$
\Theta:\,{\cal H}_{\cal O}\rightarrow {\cal O}
$$ 
given by restriction of $\Theta$ to ${\cal H}_{\cal O}$.  
%Bemerkung: ${\Bbb Q}_p[T]$ ist uber ${\Bbb Z}_p$ definierte ${\Bbb Q}_p$-algebra aber nicht jeder char von ${\Bbb Q}-p[T]$ ist uber ${\Bbb Z}_p$ definiert. Z.B. 
%$Theta$ definiert durch $\Theta|_{{\Bbb Q}_p}={\rm id}$ und $\Theta(T)=1/2$. D.h die bedingung uber ${\cal O}$ definiert ist nicht immer (automatisch) erfullt also nicht trivial.
Any $\Theta\in\hat{\cal H}_{{\cal O}}$ is determined by the collection of values $\Theta_\ell:=\Theta(T_\ell)$, $\ell\in I$, and we obtain an embedding
$$
\hat{\cal H}_{{\cal O}}\hookrightarrow {\cal O}^I.
$$
%In the following we will identify $\Theta\in \hat{\cal H}_{{\cal O}}$ with its belonging collection $\Theta=(\Theta_{\ell})_{\ell\in I}$. 
For any character $\Theta$ of ${\cal H}$ we set
$$
v_p(\Theta)={\rm inf}_{T\in{\cal H}_{\cal O}} \, v_p(\Theta(T))\in{\Bbb Q}\cup\{-\infty\}.
$$
We note that $\Theta\in\hat{\cal H}_{\cal O}$ precisely if $v_p(\Theta)>-\infty$ and in this case we obtain 
$$
v_p(\Theta)={\rm inf}_{\ell\in I} \, v_p(\Theta_{\ell})\in{\Bbb Q}_{\ge 0}.
$$
%beachte, dass $\lambda(T)$ jetzt werte in ${\cal O}$ hat fur $T\in{\cal H}_{\cal O}$, d.h. $v_p(\Theta(T_\ell))\ge 0$  !!!!
%d.h. Fälle wie $\Theta:\,{\cal O}[T]\rightarrow\bar{\Bbb Q}$, $t\mapsto p^{-1}\not\in{\cal O}$ sind nicht moglich. 
%Allgemeiner: $\Theta$ hat genau dann ein beschranktes Bild in $\bar{\Bbb Q}$, wenn $\Theta$ uber ${\cal O}$ definiert ist.
%d.h. $v_p(\Theta)>-\infty\Leftrightarrow v_p(\Theta)\ge 0$
We say that $\Theta$ is congruent to a character $\mu$ of ${\cal H}$ modulo $p^t$, $t\in{\Bbb R}$, written as $\Theta\equiv \mu\pmod{p^t}$, if
$$
v_p(\Theta-\mu)\ge t.
$$
%Beachte, dass $\hat{\cal H}_{{\cal O}}$ kEIN ${\cal O}$-module ist, denn: produkt von zwei char kein char da additivitat verletzt; summe von zwei char ist kein char, das multiplikativitat verletzt; skalares vielfahes eines chaqr ist kein char, da multipllikativitat verletzte. Es gilt also gar nichts !!
%Insbesondere ist $p^e\Theta$ nicht definiert !!!!! und wir konnen eine kongruenz $\Theta\equiv\mu\pmod{p^c}$, $c\in{\Bbb N}$, nicht als $\Theta-\mu\in p^c{\cal H}_{{\cal O}}$ ausdrücken
%
%We note that the definition of ``$v_p$'' and, hence, "$\equiv$" depend on the choice of the order ${\cal H}_{\cal O}$; a different chcoice of an order 
%${\cal H}_{\cal O}$ would yield another valuation ``$v_p$'' on characters

\bigskip

{\bf (1.3) } Let ${\bf H}$ be a ${\cal H}$-module which is finite dimensional as $\bar{\Bbb Q}_p$-vector space. We assume that 
${\bf H}$ is defined over ${\cal O}$ meaning that ${\bf H}$ contains a ${\cal O}$-submodule ${\bf H}_{{\cal O}}$ 
which is free and stable under the action of ${\cal H}_{{\cal O}}$ such that ${\bf H}={\bf H}_{\cal O}\otimes\bar{\Bbb Q}_p$. 
%d.h. eine ${\cal O}$-basis von $H_{\cal O}$ ist eine $\bar{\Bbb Q}$-basis von $H$
For any $\Theta\in\hat{\cal H}_{{\cal O}}$ we denote by
$$
{\bf H}(\Theta)=\{v\in {\bf H}:\,\mbox{ for all $T\in{\cal H}$ there is $n_T\in{\Bbb N}$ such that}\;(T-\Theta(T))^{n_T}(v)=0\}
$$
the generalized simultaneous eigenspace attached to the character $\Theta$ (or, equivalently, to the (system of) eigenvalue(s) $\Theta=(\Theta_{\ell})_{{\ell}\in I}$). 
We denote by ${\cal E}({\bf H})={\cal E}_{\cal H}({\bf H})$ the set of all eigencharacters occuring in ${\bf H}$, i.e. ${\cal E}({\bf H})$ consists of all characters $\Theta\in\hat{\cal H}_{\cal O}$ 
such that ${\bf H}(\Theta)\not=0$. Thus, elements in ${\cal E}({\bf H})$ correspond to simultaneous ${\cal H}$-eigenforms in ${\bf H}$. We note that any character 
$\Theta$ of ${\cal H}$ which occurs in ${\bf H}$ is defined over ${\cal O}$ because ${\bf H}$ is defined over ${\cal O}$, hence, $\Theta\in\hat{\cal H}_{\cal O}$.  
%Da ${\cal H}$ auf dem endl.dim $\bar{\Bbb Q}$-Vrm $H$ operiert ist bjeder EW von einem $T\in{\cal H}$ eine alg. Zahl, also $\in\bar{\Bbb Q}$.
%Da $T_\ell\in{\cal H}_{\cal O}$ einen freien ${\cal O}$-modul stabilisiert, hat die Darstellungsmatrix von $T$ eintrage aus ${\cal O}$, also ist das char polynom
%von $T$ ein normiertes polynom mit ganzen Koeffizienten, also sind die Eigenwerte von $T$, die ja NS des char polynims sind, ganz sind uber 
%${\cal O}$, also im ganzabschluss ${\cal O}$ enthalten.
%Da jedes $\Theta(T)$ der EW von $\Theta$ auf einem gewissen vektor $v$ ist, folgt, dass $\Theta(T)$ ganz ist.
%
%beachte: $H_{\cal O}$ ist ein ${\cal O}$-lattice, also insbes. ein freier ${\cal O}$-modul !!
%Anmerkung: Das gilt nicht fur Operatoren $T$ die in ${\cal H}$ aber nicht in ${\cal H}_{\cal O}$ liegen; die haben nicht ganze EW !!!
%
%BEMRKUNG (WICHTIG): Da wir mit der Darstellungsmatrix von $T$ auf ${\bf H}_{\cal O}$ argumentiert haben brauchen wir hier, dass ${\bf H}-{\cal O}$ ein freier ${\cal O}$-modul ist !!!!!!!!
Since ${\cal H}$ is commutative we thus obtain a decomposition
$$
{\bf H}=\bigoplus_{\Theta\in{\cal E}({\bf H})} {\bf H}(\Theta).
$$
%s. [Storch-Wieb, Lehrbuch der Mathmatik, Bd 2], p. 326, Satz 11.B.17 und den Kommentar nach dem Satz; beachte unsere ${\cal H}$-moduln $H^i(\gamma\backslash X,{\cal L}_\Theta)$ sind endlich dimensional !
%Bessere Referenz: [Jacobson Lectures in Abstract Algebra Bd 2], chapter IV, 9, Theorem 7, p. 134. Ais diesem Satz 7 folgt dass $H=\bigoplus_\mu H(\mu)$ und das $H$ iber ${\cal O}$ definiert ist folgt dass alle $\mu$ die vorkommen auch uber ${\cal O}$ definiert sind, d.h. in $\hat{\cal H}_{\cal O}$ liegen.

\bigskip

Let ${\bf H}$ be a ${\cal H}$-module which is finite dimensional as $\bar{\Bbb Q}_p$-vector space and which is defined over ${\cal O}$ with respect to the 
free ${\cal O}$-submodule ${\bf H}_{\cal O}$.  
%meaning that $H$ contains an ${\cal O}$-lattice $H_{{\cal O}}$ which is stable under ${\cal H}_{{\cal O}}$.

\bigskip

{\bf (1.4) Definition. } Let $\Theta\in\hat{\cal H}_{{\cal O}}$ and let $c\in{\Bbb Q}$. We define the $\pmod{p^c}$-reduced multiplicity of $\Theta$ in ${\bf H}$ as
$$
{m}_{\bf H}(\Theta,c,p)=\sum_{\mu\in\hat{\cal H}_{{\cal O}}\atop \mu\equiv\Theta\pmod{p^c}} {\rm dim}\,{\bf H}(\mu).  
$$
If the prime $p$ is understood we will omit "$p$" and write more simply ${m}_{\bf H}(\Theta,c)$ instead.

\bigskip

{\bf (1.5) Remark. } 
%1.) The definition of the $\pmod{p^c}$-reduced multiplicities does not make use of the integral structure on $H$ (the definition of congruence between eigencharacters $\Theta$ makes use of the integral structure on $H$ and in fact depends on the choice of the integral structure).
%; they can be defined in the same way even if $H$ is not defined over $\cal O$. 
%(cf. (2.3) below). 
%if $H$ is defined over $\cal O$ then we have $\lambda\equiv\mu\pmod{c}$ with $c\ge 0$ for all $\mu,\lambda$. THis need no longer hold if $H$
%is not defined over $\cal O$. it an even happen that $\lambda\eqiov\mu\pmod{p^c}$ only for $c=-\infty$
%
%
%\medskip
%
%2.) 
The reduced multiplicities are related to the multiplicity ${m}_{\bf H}(\Theta):={\rm dim}\,{\bf H}(\Theta)$ by
$$
\lim_{c\rightarrow \infty\atop c\in{\Bbb N}} {m}_{\bf H}(\Theta,c)={m}_{\bf H}(\Theta).
$$

%\medskip

%{\it Proof. } By definition we know that ${\bf m}_H(\Theta,c)\ge {\bf m}_H(\Theta)$ for all $c$. Morover, the sequence ${\bf m}_H(\Theta,c)$ is decreasing in $c$, hence, it becomes constant: ${\bf m}_H(\Theta,c)={\bf m}$ for some ${\bf m}\in{\Bbb N}$ and all $c\ge c_0$ (in particular, the sequence converges).  
%We set ${\bf E}(c)=\{\mu\in{\bf E}(H):\,\mu\equiv\Theta\pmod{p^c}\}$, thus, in the definition of $ {\bf m}_H(\Theta,c)$ we only have to sum over ${\bf E}(c)$ and ${\rm dim}\,H(\mu)>0$ is strictly positive for all $\mu\in {\bf E}(c)$. Since the sequence ${\bf m}_H(\Theta,c)$ is constant for all $c\ge c_0$ we obtain ${\bf E}(c_0)={\bf E}(c)$ for all $c\ge c_0$. Hence, ${\bf E}(c_0)=\{\Theta\}$ 
%%ist $\mu\equiv\Theta\pmod{p^c}$ fur alle $c\ge c_0$, dann folgt $\mu=\Theta$
%and ${\bf m}_H(\Theta,c)={\bf m}_H(\Theta)$ for all $c\ge c_0$. 

\bigskip

{\bf (1.6) Mod ${p^c}$ reduced multiplicities and mod $p^c$ reduction. } We set ${\bf H}_{{\cal O}}(\mu)={\bf H}(\mu)\cap {\bf H}_{{\cal O}}$, hence, 
$$
\bigoplus_{\mu\in{\cal E}({\bf H})}{\bf H}_{\cal O}(\mu)\subseteq{\bf H}_{\cal O}.
$$
We let $c\in{\Bbb Q}_{\ge 0}$ and define the ideal ${\mathfrak a}=\{x\in{\cal O}:\,v_p(x)\ge c\}\le{\cal O}$. We denote by 
$$
\bar{\cal O}=\frac{\cal O}{\mathfrak a},\quad \bar{\cal H}_{\cal O}=\frac{{\cal H}_{\cal O}}{{\mathfrak a}{\cal H}_{\cal O}}\quad\mbox{and}\quad\bar{\bf H}_{{\cal O}}
=\frac{{\bf H}_{{\cal O}}}{{\mathfrak a} {\bf H}_{{\cal O}}}
$$ 
the ${\rm mod}\,{\mathfrak a}$-reductions and we denote by 
$\bar{\bf H}_{{\cal O}}(\mu)$ the image of ${\bf H}_{{\cal O}}(\mu)$ in $\bar{\bf H}_{{\cal O}}$. Hence, $\bar{\bf H}_{{\cal O}}$ and 
$\bar{\bf H}_{{\cal O}}(\mu)$ are $\bar{\cal H}_{\cal O}$-modules.   
%bea.: ${\cal H}$ bildet $H_{\cal O}$ auf sich ab (da $H$ uber ${\cal O}$ definiert) und ebenso $H_{\cal O}(\mu)$
We assume that ${\bf H}_{\cal O}(\mu)\le{\bf H}_{\cal O}$ is a free ${\cal O}$-submodule (in our later applications ${\bf H}$ will be defined over a finite extension 
$E/{\Bbb Q}_p$, hence, we may replace ${\cal O}$ by ${\cal O}_E$ which is a P.I.D and the assumption holds). Since ${\bf H}_{\cal O}(\mu)\le{\bf H}_{\cal O}$ is 
saturated,  
%denn ist $v\in{\bf H}_{\cal O}$ mit $\alpha v\in{\bf H}_{\cal O}(\mu)$ dann folgt $v\in{\bf H}(\mu)$ und $v\in{\bf H}_{\cal O}$ also $v\in{\bf H}_{\cal O}(\mu)$.
$\bar{\bf H}_{{\cal O}}(\mu)$ is a free $\bar{\cal O}$-module 
%denn ${\bf H}_{\cal O}(\mu)$ hat ein direktes Komplement in ${\bf H}_{\cal O}$ (da eben saturated), d.h. ${\bf H}_{\cal O}={\bf H}_{\cal O}(\mu)\oplus V$. Ist
%Sei $b_i$ eine basis von ${\bf H}_{\cal O}(\mu)$. Die lasst sich zu basis von ${\bf H}_{\cal O}$ erganzen. Dann folgt aber aus  
%$\sum \alpha_i b_i\in{\mathfrak a}{\bf H}_{\cal O}$ dass $\alpha_i\in{\mathfrak a}$, d.h. die $b_i$ sind l.u. auch in 
%$\bar{\bf H}_{{\cal O}}=\frac{{\bf H}_{{\cal O}}}{{\mathfrak a} {\bf H}_{{\cal O}}}$.
%
which is annihilated by $(T-\bar{\mu}(T))^{n_T}$ for all $T\in\bar{\cal H}_{\cal O}$; here,  
$\bar{\mu}=\mu\,{\rm mod}\,p^c:\, \bar{\cal H}_{\cal O}\rightarrow\bar{\cal O}$, $T+{\mathfrak a}{\cal H}_{\cal O}\mapsto \mu(T)+{\mathfrak a}$ 
is the ${\rm mod}\,{\mathfrak a}$-reduced character. 
%(note that $\bar{\bf H}_{{\cal O}}(\mu)$ and, hence, $\bar{\bf H}_{{\cal O}}$ in general is not semi simple).
%note that $\mu(s)\in{\cal O}$ for all $\mu\in\hat{\cal H}_{\cal O}$
If the decomposition of ${\bf H}$ as a sum of generalized eigenspaces is defined over $\cal O$, i.e. if
$\bigoplus_{\mu\in{\cal E}({\bf H})} {\bf H}_{{\cal O}}(\mu)={\bf H}_{{\cal O}}$ then 
$$
\bigoplus_{\mu\in{\cal E}({\bf H})}\bar{\bf H}_{\cal O}(\mu)=\bar{\bf H}_{\cal O}\leqno(1)
$$
%
%denn: sei $\sum_i \alpha_i (v_i+{\mathfrak a}{\bf H}_{\cal O})\in{ßmathfrak a}{\bf H}_{\cal O}$, d.h. $\sum_i \alpha_i (v_i+{\mathfrak a}{\bf H}_{\cal O})=\sum_i \beta_i (b_i+{\mathfrak a}{\bf H}_{\cal O})$ fur gewisse $\beta_i\in{\mathfrak a}$ wobei $(b_i)$ eine ${\cal O}$-Basis von ${\bf H}_{\cal O}$ ist, dann folgt da die 
%Summe der ${\bf H}_{\cal O}(\mu)$ direkt ist, dass $\alpha_i\in{ßmathfrak a}$, d.h. $\alpha_i=0$ in $\bar{\cal O}$
%
as $\bar{\cal O}$-modules and, hence, we obtain for the multiplicity of ${\Theta}\,{\rm mod}\,p^c$ in $\bar{\bf H}_{\cal O}$
$$
m_{\bar{\bf H}_{{\cal O}}}(\Theta\,{\rm mod}\,p^c):={\rm dim}_{\bar{\cal O}}\,\sum_{\mu\in\hat{\cal H}_{\cal O}\atop\,\bar{\mu}=\bar{\Theta}}\bar{\bf H}_{\cal O}(\mu)=m_{\bf H}(\Theta,c,p). 
$$
%bea., dass $ \mu\equiv\Theta\pmod{p^c}  \Leftrightarrow   v_p(\mu-\Theta(T_\ell))\ge c \forall \ell  \Leftrightarrow   \mu-\Theta(T_\ell)\in{\mathfrak a}\Leftrightarrow \mu(T_\ell)\equiv \Theta(T_\ell)\pmod{\mathfrak a}   \Leftrightarrow     \mu(T_\ell)+ \mu(T_\ell)=\Theta(T_\ell)+ \mu(T_\ell)\equiv \Theta(T_\ell)    \Leftrightarrow    \bar{\mu}_c(T_\ell)=\bar{\Theta}_c(T_\ell)$. Beachte noch dass $dim\,\bar{H}_{{\cal O}}(\mu)=dim\,\bar{H}(\mu)$.
In general, we only obtain an inclusion of a (not necessarily direct) sum 
$$
\sum_{\mu\in{\cal E}({\bf H})}\bar{\bf H}_{\cal O}(\mu)\subseteq\bar{\bf H}_{\cal O},
$$
%i.e. $\bar{H}_{{\cal O}}$ is nearly a semi simple ${\cal H}$-module (if the $\bar{H}_{{\cal O}}(\mu)$ are s.s)
%
%$\bar{H}_{{\cal O}}=\sum_{\mu} \bar{H}_{{\cal O}}(\mu)$ da $H_{{\cal O}}=\sum_{\mu} H_{{\cal O}}(\mu)$
%Direktheit der Summe gilt denn: ist $\sum_mu \bar{v_\mu}=0$, $\bar{v_\mu}=v_\mu+{\mathfrak a}$ mit $v_\mu\in H_{{\cal O}}(\mu)$, dann folgt $v\in {\mathfrak a} H_{{\cal O}}=\bigoplus_\mu {\mathfrak a} H_{{\cal O}}(\mu)$ woraus folgt dass $v_\mu\in {\mathfrak a} H_{{\cal O}}$, d.h. \bar{v_\mu}=0$.
which yields an inequality
$$
m_{\bar{\bf H}_{{\cal O}}}(\Theta\,{\rm mod}\,p^c)\le m_{\bf H}(\Theta,c,p).
$$
In this sense, the reduced multiplicity ${m}_{\bf H}(\Theta,c,p)$ is a substitute for the multiplicity of the mod $p^c$-reduction of 
${\Theta}$ in $\bar{H}_{{\cal O}}$ if the primary decomposition of ${\bf H}$ is not defined over ${\cal O}$.

\bigskip

{\bf (1.7) Higher Congruences }. %We explain the application of reduced multiplicities to proving congruences between automorphic forms. 
Let ${\cal H}=\bar{\Bbb Q}_p[T_\ell,\,\ell\in I]$ and ${\cal H}'=\bar{\Bbb Q}_p[T_\ell,\ell\in I']$ and let ${\bf H}$ resp. ${\bf H}'$ be a ${\cal H}$ 
resp. ${\cal H}'$-module which is defined over ${\cal O}$ with respect to the lattice ${\bf H}_{\cal O}$ resp. ${\bf H}'_{\cal O}$ as in (1.3). Let 
$$
\Phi^\vee:\,\hat{\cal H}_{\cal O}\rightarrow\hat{\cal H}'_{\cal O}
$$
be a map and let $m\in{\Bbb N}$. 
%We define a map
%$$
%\Phi^\vee:\,{\cal E}({\bf H})\rightarrow{\cal E}({\bf H}')
%$$
%by sending $\mu$ to $\Phi^\vee(\mu)=\mu\circ \Phi$, i.e. the diagram
%$$
%\begin{array}{ccc}
%{\cal H}&\stackrel{\Phi}{\leftarrow}&{\cal H}'\\
%\mu\;\searrow&&\swarrow\;\Phi^\vee(\mu)\\
%&\bar{\Bbb Q}&\\
%\end{array}
%$$
%commutes.
If for any $\Theta\in\hat{\cal H}_{\cal O}$ there is a rational number $c=c(\Theta)\ge m$ such that
$$
m_{\bf H}(\Theta,c)=m_{{\bf H}'}(\Phi^\vee(\Theta),c)
$$ 
then for any $\Theta\in{\cal E}({\bf H})$ there is a eigencharacter $\Theta'\in{\cal E}({\bf H}')$ such that
$$
\Theta'\equiv\Phi^\vee(\Theta)\pmod{p^c}.
$$
In different words there is a map on {\it eigencharacters}
$$
\Psi^{[c]}:\,{\cal E}_{\cal H}({\bf H})\rightarrow{\cal E}_{{\cal H}'}({\bf H}')
$$
such that $\Psi^{[c]}(\Theta)\equiv \Phi^\vee(\Theta)\pmod{p^c}$. Thus, by comparing reduced multiplicities we can establish congruences between 
eigencharacters 
in ${\cal E}({\bf H}')$ and lifts of eigencharacters in ${\cal E}({\bf H})$ or, equivalently, a mod $p^c$ transfer from 
${\cal E}_{\cal H}({\bf H})$ to ${\cal E}_{{\cal H}'}({\bf H}')$. In particular, if ${\cal H}={\cal H}'$ and $\Phi^\vee={\rm id}$ then the 
set of identities 
$$
m_{\bf H}(\Theta,c)=m_{{\bf H}'}(\Theta,c),\quad\Theta\in\hat{\cal H}_{\cal O},
$$
implies that for any $\Theta\in{\cal E}({\bf H})$ a congruence
$$
\Theta\equiv\Theta'\pmod{p^m}
$$
holds for some $\Theta'\in{\cal E}({\bf H}')$.

\bigskip

{\bf (1.8) Systems of higher congruences. } Slightly refining the above discussion we can give a set of simple 
identitites between mod $p^c$ reduced multiplicities which implies the existence of $p$-adic {\it continuous} 
families of eigencharacters (note that such families are defined by a system of congruences between their members). To be more precise, we let 
${\bf G}/{\Bbb Q}$ be a connected reductive algebraic group with maximal split torus ${\bf T}$ and we denote by $X({\bf T})$ the group of 
${\Bbb Q}$-characters of ${\bf T}$. $X({\bf T})$ is a finitely generated, free abelian group which we write using additive notation. 
%note that ${\bf T}$ is connnected, hence, $X({\bf T})$ is a finitely generated, free abelian group
For any $\lambda\in X({\bf T})$ we denote by $v_p(\lambda)$ the largest integer $m$ such that $\lambda\in p^m X({\bf T})$.
Thus, if we identify $X({\bf T})\cong {\Bbb Z}^k$ via the choice of a basis $(\gamma_i)_i$ of $X({\bf T})$ then $v_p(\sum_i z_i \gamma_i)=\inf_{i} v_p(z_i)$.
%Dies ist unabhangig von der gewahlten basis $(b_i)_i$, denn $p^m X({\bf T})$ entspricht immer $ßsum_i z_i\gamma_i$ mit $p^m|z_i$ fur alle $i$
%compare section (6.2)

\bigskip

{\bf (1.9) Proposition. }{\it Let ${\cal R}\subset X({\bf T})$ be a subset and let $({\bf H}_\lambda)$, $\lambda\in {\cal R}$, be a family of finite dimensional 
${\cal H}$-modules which are defined over ${\cal O}$. Assume there are ${\sf a},{\sf b}\in{\Bbb Q}$ with the following property: if 
$\lambda\equiv\lambda'\pmod{(p-1)p^m X({\bf T})}$, then for any $\Theta\in\hat{\cal H}_{\cal O}$ there is $c=c(\Theta)\ge {\sf a}(m+1)+{\sf b}$ with
$$
m_{{\bf H}_\lambda}(\Theta,c)=m_{{\bf H}_{\lambda'}}(\Theta,c),
$$
i.e. the transfer $\Psi^{[{\sf a}(m+1)+{\sf b}]}:\,{\cal E}({\bf H}_{\lambda})\rightarrow {\cal E}({\bf H}_{\lambda'})$ corresponding to $\Phi^\vee={\rm id}$ exists. Then, any 
$\Theta\in{\cal E}({\bf H}_{\lambda_0})$ fits in a $p$-adic continuous family of eigencharacters, i.e. there is a family 
$(\Theta_\lambda)$, $\lambda\in{\cal R}$, such that

\begin{itemize}

\item $\Theta_\lambda\in{\cal E}({\bf H}_\lambda)$

\item $\Theta_{\lambda_0}=\Theta$

\item $\lambda\equiv\lambda'\pmod{(p-1)p^m X({\bf T})}$ implies $\Theta_\lambda\equiv\Theta_{\lambda'}\pmod{p^{{\sf a}(m+1)+{\sf b}}}$.

\end{itemize}

}

\medskip

{\it Proof. } For any weight $\mu\in X({\bf T})$ we set ${\cal R}_\mu=\{\lambda\in{\cal R}:\,\lambda\equiv\mu\pmod{(p-1) X({\bf T})}\}$. We first 
construct a $p$-adic family $(\Theta_\lambda)_\lambda$ satisyfing the above conditions with $\lambda$ only running over ${\cal R}_{\lambda_0}$. To 
this end, we enumerate the weights $\lambda$ in ${\cal R}_{\lambda_0}$ in a sequence $\lambda_0,\lambda_1,\lambda_2,\lambda_3,\ldots$. 
%note that $X({\bf T})$ is countable
We inductively construct elements $\Theta_{\lambda_i}\in{\cal E}({\bf H}_{\lambda_i})$, $i=0,1,2,3,\ldots$ such that $\Theta_{\lambda_0}=\Theta$ and 
${\lambda_i}\equiv {\lambda_j}\pmod{(p-1)p^m X({\bf T})}$ implies $\Theta_{\lambda_i}\equiv \Theta_{\lambda_j}\pmod{p^{{\sf a}(m+1)+{\sf b}}}$. Clearly, we 
set $\Theta_{\lambda_0}=\Theta$. Assume that $\Theta_{\lambda_0},\ldots,\Theta_{\lambda_n}$ have been defined such that 
$\lambda_i\equiv \lambda_j\pmod{(p-1)p^m X({\bf T})}$ implies that $\Theta_{\lambda_i}\equiv \Theta_{\lambda_j}\pmod{p^{{\sf a}(m+1)+{\sf b}}}$ for all 
$i,j=0,\ldots,n$. To define $\Theta_{\lambda_{n+1}}$ we select $a\in\{0,1,2,\ldots,n\}$ such that
$$
v_p(\lambda_{n+1}-\lambda_a)\ge v_p(\lambda_{n+1}-\lambda_i)\quad\mbox{for all}\; i=0,\ldots,n
$$
We set $w_1=v_p(\lambda_{n+1}-\lambda_a)$, hence, $\lambda_{n+1}-\lambda_a\in (p-1)p^{w_1} X({\bf T})$ (note that 
$\lambda_a-\lambda_{n+1}\in (p-1)X({\bf T})$ because $\lambda_{n+1},\lambda_a\in{\cal R}_{\lambda_0}$). By (2.4) Remark there is 
$\Theta\in{\cal E}({\bf H}_{\lambda_{n+1}})$ such that $\Theta\equiv\Theta_{\lambda_a}\pmod{p^{{\sf a}(w_1+1)+{\sf b}}}$. We then set 
$\Theta_{\lambda_{n+1}}$ equal to this $\Theta$. 

Let $i\in\{0,\ldots,n\}$ be arbitrary and set $w_3=v_p(\lambda_{n+1}-\lambda_i)$, hence, $\lambda_{n+1}\equiv\lambda_i\pmod{(p-1)p^{w_3} X({\bf T})}$. 
We have to show that $\Theta_{\lambda_{n+1}}\equiv \Theta_{\lambda_i}\pmod{p^{{\sf a}(w_3+1)+{\sf b}}}$. To this end we set 
$w_2=v_p(\lambda_a-\lambda_i)$.
$$
\left.\begin{array}{cccc}
&&&\bullet \lambda_{n+1}\\
w_1\{&&/&\\
&\lambda_a\bullet&&\\
w_2\{&|&&\\
&\lambda_i\bullet&&\\
\end{array}\right\}w_3.
$$
We know that $\Theta_{\lambda_{n+1}}\equiv \Theta_{\lambda_a}\pmod{p^{{\sf a}(w_1+1)+{\sf b}}}$ by definition of 
$\Theta_{\lambda_{n+1}}$ and that $\Theta_{\lambda_a}\equiv \Theta_{\lambda_i}\pmod{p^{{\sf a}(w_2+1)+{\sf b}}}$ by our induction 
hypotheses, hence,
$$
\Theta_{\lambda_{n+1}}\equiv \Theta_{\lambda_i}\pmod{p^{{\sf a}({\rm min}\{w_1,w_2\}+1)+{\sf b}}}.\leqno(2)
$$ 
We distinguish cases.

{\it Case A} $w_2> w_1$. In this case ${\rm min}\{w_1,w_2\}=w_1$ and $w_3=w_1$ by the $p$-adic triangle inequality. Hence, equation (2) implies 
that $\Theta_{\lambda_{n+1}}\equiv \Theta_{\lambda_i}\pmod{p^{{\sf a}(w_3+1)+{\sf b}}}$. 

{\it Case B} $w_2< w_1$. In this case ${\rm min}\{w_1,w_2\}=w_2$ and $w_3=w_2$. Hence, equation (2) implies that 
$\Theta_{\lambda_{n+1}}\equiv \Theta_{\lambda_i}\pmod{p^{{\sf a}(w_3+1)+{\sf b}}}$. 

{\it Case C} $w_2=w_1$. In this case ${\rm min}\{w_1,w_2\}=w_1$. On the other hand, by the choice of $a$ we know that $w_1\ge w_3$; thus 
equation (2) yields $\Theta_{\lambda_{n+1}}\equiv \Theta_{\lambda_i}\pmod{p^{{\sf a}(w_3+1)+{\sf b}}}$.

Thus, $(\Theta_\lambda)_\lambda$, $\lambda\in{\cal R}_{\lambda_0}$, is a $p$-adic continuous family. To obtain a $p$-adic family 
$(\Theta_\lambda)$ with $\lambda$ running through all of ${\cal R}$ we denote by $\{\mu_0=\lambda_0,\mu_1,\ldots,\mu_r\}$ a system of 
representatives for $X({\bf T})/(p-1) X({\bf T})$. For any $i=1,\ldots,r$ we construct in the same way as above a $p$-adic family 
$(\Theta_\lambda)_\lambda$ with $\lambda$ running through ${\cal R}_{\mu_i}$. 
Since $\mu_i\not\equiv\mu_j\pmod{(p-1)X({\bf T})}$ if $i\not=j$ the union $(\Theta_\lambda)$, $\lambda\in\bigcup_i {\cal R}_{\mu_i}$ then is a 
$p$-adic family satisfying the requirements of the Proposition. This completes the proof.
%$bigcup_i {\cal R}_{\mu_i}={\cal R}$

\bigskip

%{\bf (2.7) Remark. } (2.5) Remark and (2.6) Proposition show that in order to prove the existence of congruences between eigencharacters or 
%the existence of $p$-adic families of eigencharacters we have to compare reduced multiplicites. In the next section we will explain how 
%this can be done by using a trace formula.

\section{Mod $p^c$ Transfer}

We show that reduced multiplicities can be computed as traces of certain operators and we use this to establish a ``mod $p^c$ transfer'' for 
eigencharacters.

\bigskip

{\bf (2.1) The dual map. } We let ${\cal H}=\bar{\Bbb Q}_p[T_\ell,\,\ell\in I]$ and ${\cal H}'=\bar{\Bbb Q}_p[T_\ell,\,\ell\in I']$ be countably generated 
polynomial algebras algebras and we set ${\cal H}_{\cal O}={\cal O}[T_\ell,\,\ell\in I]$ and ${\cal H}'_{\cal O}={\cal O}[T_\ell,\,\ell\in I']$. We assume 
that there is a morphism of ${\cal O}$-algebras
$$
\Phi:\,{\cal H}'_{\cal O}\rightarrow {\cal H}_{\cal O}.
$$
The morphism $\Phi$ induces a dual map
$$
\Phi^\vee:\,\hat{\cal H}_{\cal O}\rightarrow\hat{\cal H}'_{\cal O}\leqno(1)
$$
by sending $\mu$ to $\Phi^\vee(\mu)=\mu\circ \Phi$, i.e. the diagram
$$
\begin{array}{ccc}
{\cal H}&\stackrel{\Phi}{\leftarrow}&{\cal H}'\\
\mu\;\searrow&&\swarrow\;\Phi^\vee(\mu)\\
&\bar{\Bbb Q}_p&\\
\end{array}
$$
commutes. Thus, for any $v_\mu\in {\bf H}(\mu)$, $\mu\in\hat{\cal H}_{\cal O}$, and any $T'\in{\cal H}'$ we obtain
$$
\Phi(T') v_\mu=\Phi^\vee(\mu)(T')v_\mu.
$$
%denn $\Phi^\vee(\mu)(T')v_\mu=\mu(\Phi(T')) v_\mu=\Phi(T') v_\mu$

\medskip

{\bf (2.2) Remark. } Let ${\bf H}$ resp. ${\bf H}'$ be a ${\cal H}$ resp. ${\cal H}'$-module which is defined over ${\cal O}$ with respect to the lattice 
${\bf H}_{\cal O}$ resp. ${\bf H}'_{\cal O}$. In general, $\Phi^\vee$ does not induce a mapping on eigencharacters
$$
\Phi^\vee:\,{\cal E}({\bf H})\rightarrow{\cal E}({\bf H}').
$$
We want to examine whether this is the case modulo powers of $p$, i.e. whether for any $\Theta\in{\cal E}({\bf H})$ there is a $\Theta'\in{\cal E}({\bf H}')$ such that
$$
\Phi^\vee(\Theta)\equiv \Theta'\pmod{p^c}.
$$
This is equivalent to the existence of a map 
$$
\Psi^{[c]}:\,{\cal E}({\bf H})\rightarrow{\cal E}({\bf H}')
$$
such that $\Psi^{[c]}(\Theta)\equiv \Phi^\vee(\Theta)\pmod{p^c}$ for all $\Theta\in{\cal E}({\bf H})$. By (1.7) we can establish the 
existence of the map $\Psi^{[c]}$ by comparing the mod $p^c$ reduced multiplicties of $\Theta$ and $\Phi^\vee(\Theta)$ in ${\bf H}$ 
and ${\bf H}'$, i.e. we need to compute reduced multiplicities. This will be based on the following 
Lemma which expresses the reduced multiplicites $m_{\bf H}(\Theta,c)$ and $m_{{\bf H}'}(\Phi^\vee (\Theta),c)$ as traces of a certain 
element in ${\cal H}'$.

%{\bf Proposition. } {\it We fix a finite subset ${\cal S}\subseteq\hat{\cal H}_{{\cal O}}$ and integers $M,m\in{\Bbb N}$ such that $m\ge (2|{\cal S}|+1)\log_p M$. Then there is a rational number $c=c({\cal S},M,m)$ satisfying $m-(2|{\cal S}|+1)\log_p M\le c\le m$ with the following property. For any $\Theta\in{\cal H}_{{\cal O}}$ there is a Hecke operator $e(\Theta)\in{\cal H}_{\bar{\Bbb Q}}$ satisyfing 
%
%(i) $e(\Theta)\in\frac{1}{p^{c|{\cal S}|}}{\cal H}_{{\cal O}}$
%i.e. the $p$-adic value of the denominator of $e(\Theta)$ is bounded by $c|{\cal S}|$
%
%
%(ii) for any $H$ be any finite dimensional ${\cal H}_{\bar{\Bbb Q}}$-module defined over ${\cal O}$ such that ${\rm dim}\,H\le M$ and ${\bf E}(H)\subset {\cal S}$ we have
%$$
%{\rm tr}\,e(\Theta)|_H \equiv {\bf m}_H(\Theta,c) \pmod{p^{\log_p M}} 
%$$
%
%
%}
%
%We set
%$$
%v_p(T)={\rm inf}\,\{v_p(x),\,x\in{\cal O}:\,xT\in{\cal H}_{{\cal O}}\}
%$$
%Thus, $\xi T\in{\cal H}_{{\cal O}}$ for all $\xi\in{\cal O}$ satisfying $v_p(\xi)\> v_p(T)$. 
%
%
%\bigskip

\bigskip

{\bf (2.3) Reduced multiplicities as traces. } From now on we assume that $\Phi:\,{\cal H}'_{\cal O}\rightarrow{\cal H}_{\cal O}$ is 
{\it surjective}. 

%In the following Lemma we compute reduced multiplicities as traces of certain elements in the algebras ${\cal H}$, ${\cal H}'$.

\bigskip

{\bf Lemma. } {\it Assume that ${\rm dim}\,{\bf H},\,{\rm dim}\,{\bf H}'\le \frac{1}{2} M$ for some $M\in 2{\Bbb N}$ and let
$m\in {\Bbb N}$. For any $\Theta\in\hat{\cal H}_{\cal O}$ there is an element
$e(\Theta)\in{\cal H}'$ and a rational number $c=c(\Theta)\ge m-(M+\frac{3}{2})\log_p M$ such that the following holds

\begin{itemize}

\item $e(\Theta)\in\frac{1}{\xi}{\cal H}_{\cal O}'$, where $\xi\in{\cal O}$ with $v_p(\xi)\le Mm$

\item ${\rm tr}\,(e(\Theta)|{\bf H}')\equiv m_{{\bf H}'}(\Phi^\vee(\Theta),c)\pmod{p^{\log_p M}}$

\item ${\rm tr}\,(\Phi(e(\Theta))|{\bf H})\equiv m_{\bf H}(\Theta,c)\pmod{p^{\log_p M}}.$

\end{itemize}

(We note that $\log_p M\in{\Bbb R}_{\ge 0}$ and the congruence has to be understood as in (1.1).)
}

%beachte: $M\ge 2$, da $M\in 2{\Bbb N}$. insbesondere ist $\log_p M\ge 0$. ware $M=1$ so folgte ${\rm dim}\,{\bf H},\,{\rm dim}\,{\bf H}'=0$, d.h. 
%ein fall den wir nicht betrachten mussen.

\medskip

{\it Proof. } Let $\Theta\in\hat{\cal H}_{\cal O}$. We proceed in steps.

\medskip

a.) We first define $c$. We abbreviate $l=\log_p M$. 
% $l=\log_p (M/2)$ wurde reichen, da ${\rm dim}\, {\bf H},{\bf H}'ßle M/2$gilt 
We set
$$
\Omega=\Omega(\Theta)=\{v_p(\mu-\Theta),\,\mu\in{\cal E}({\bf H})\}\cup \{v_p(\mu'-\Phi^\vee(\Theta)),\,\mu'\in{\cal E}({\bf H}')\}\subseteq{\Bbb Q}
$$
and we define the interval
$$
{\bf I}%=[m-(|{\cal E}({\bf H})|+|{\cal E}({\bf H}')|+3/2)l,m]
=\{r\in{\Bbb Q}:\,m-(M+3/2)l\le r\le m\}\subset{\Bbb Q}_{\ge 0}.
$$
%beachte: da $l\ge 0$ ist das ein echtes, d.h. nicht entartetes Intervall (linker endpunkt kleiner rechter endpunkt)
Since ${\bf I}$ has length $(M+3/2)l$ and since $|\Omega|\le|{\cal E}({\bf H})|+|{\cal E}({\bf H}')|\le {\rm dim}\,{\bf H}+{\rm dim}\,{\bf H}'\le M$ 
there is a $c\in{\bf I}$ such that
$$
[c,c+l]\cap \Omega=\emptyset.
$$
%(it would suffice if ${\bf I}$ had length $(|{\cal E}({\bf H})|+|{\cal E}({\bf H}')|+1+\epsilon)l$ with $\epsilon>0$). 
Thus, if $\gamma\in{\cal E}({\bf H})$ with $v_p(\gamma-\Theta)\ge c$ then we know that $v_p(\gamma-\Theta)\in\Omega\cap[c,\infty)$, hence,
$$
v_p(\gamma-\Theta)> c+l.\leqno(2)
$$
Similarly, if $\gamma'\in {\cal E}({\bf H}')$ with $v_p(\gamma'-\Phi^\vee(\Theta))\ge c$ then
$$
v_p(\gamma'-\Phi^\vee(\Theta))> c+l.\leqno(2')
$$
We note that the number $c$ obviously satisfies 
$$
m-(M+3/2)l\le c\le m.
$$ 

\medskip

b.) Next we define $e(\Theta)$. We note that for any $\mu\in{\cal E}({\bf H})$ with $\mu\not\equiv \Theta\pmod{p^c}$, i.e.  $v_p(\mu-\Theta)<c$, there is a element $T_\mu'\in{\cal H}'_{\cal O}$ such that
$$
\mu(\Phi(T_\mu'))\not\equiv \Theta(\Phi(T_\mu'))\pmod{p^c}.
$$
(note that we assume $\Phi$ to be surjective); hence, 
$$
\Phi^\vee(\mu)(T_\mu')\not\equiv \Phi^\vee(\Theta)(T_\mu')\pmod{p^c}.\leqno(3)
$$
Similarly, for any $\mu'\in {\cal E}({\bf H}')$ with $\mu'\not\equiv \Phi^\vee(\Theta)\pmod{p^c}$, i.e.  $v_p(\mu'-\Phi^\vee(\Theta))<c$,
there is a element $T_{\mu'}\in{\cal H}'_{\cal O}$ such that
$$
\mu'(T_{\mu'})\not\equiv \Phi^\vee(\Theta)(T_{\mu'})\pmod{p^c}.\leqno(3')
$$
We set
$$
\xi=\prod_{\mu\in{\cal E}({\bf H})\atop \mu\not\equiv \Theta\pmod{p^c}} \Phi^\vee(\Theta)(T_\mu')-\Phi^\vee(\mu)(T_\mu')\,
\prod_{\mu'\in{\cal E}({\bf H}')\atop \mu'\not\equiv \Phi^\vee(\Theta)\pmod{p^c}} \Phi^\vee(\Theta)(T_{\mu'})-\mu'(T_{\mu'}) \;(\in{\cal O})
$$
%$\xi\in{\cal O}$ da $\Theta,\mu in\bar{\cal H}_{\cal O}$ (bea.: ${\cal S}\subseteq\bar{\cal H}_{\cal O}$)
and we define
$$
e(\Theta)=\frac{1}{\xi}\,\prod_{\mu\in{\cal E}({\bf H})\atop \mu\not\equiv \Theta\pmod{p^c}} T_\mu'-\Phi^\vee(\mu)(T_\mu')\,
\prod_{\mu'\in{\cal E}({\bf H}')\atop \mu'\not\equiv \Phi^\vee(\Theta)\pmod{p^c}}T_{\mu'}-\mu'(T_{\mu'})\;\in\frac{1}{\xi}{\cal H}_{\cal O}'.
$$
%bea.: $\mu_{T_\mu}\in{\cal O}$ a $\mu\in{\cal E}({\bf H})\subseteq\bar{\cal H}_{\cal O}$
Equations $(3)$ and $(3')$ imply that 
$$
v_p(\xi)\le |{\cal E}({\bf H})|c+|{\cal E}({\bf H}')|c\le Mc\le Mm.
$$
which in particular implies the first claim of the Lemma.

\medskip

c.) We compute the trace of $e(\Theta)$ on ${\bf H}'$. We write
$$
{\bf H}'=\bigoplus_{\gamma'\in{\cal E}({\bf H}')} {\bf H}'(\gamma').
$$
Since ${\cal H}'$ is commutative, there is for any $\gamma'\in{\cal E}({\bf H}')$ a basis ${\cal B}(\gamma')$ of ${\bf H}'(\gamma')$ such that all 
$T'\in{\cal H}'$ are represented on ${\bf H}'(\gamma')$ by an upper triangular matrix:
$$
{\cal D}_{{\cal B}(\gamma')}(T')=\left(\begin{array}{ccc}\gamma'(T') &&*\\&\ddots&\\&&\gamma'(T') \end{array} \right).
$$
Then $e(\Theta)$ is represented on ${\bf H}'(\gamma')$ by the matrix
$$
{\cal D}_{{\cal B}(\gamma')}(e(\Theta))=\left(\begin{array}{ccc}x &&*\\&\ddots&\\&&x' \end{array} \right),
$$
where
$$
x=x(\gamma')=\frac{1}{\xi}\,\prod_{\mu\in{\cal E}({\bf H})\atop \mu\not\equiv \Theta\pmod{p^c}} \gamma'(T_\mu')-\Phi^\vee(\mu)(T_\mu')\,
\prod_{\mu'\in{\cal E}({\bf H}')\atop \mu'\not\equiv \Phi^\vee(\Theta)\pmod{p^c}}\gamma'(T_{\mu'})-\mu'(T_{\mu'})   \in\frac{1}{\xi}{\cal O}.
$$
%d.h. $x$ muss nicht in ${\cal O}$ liegen, insbes. ist  ${\rm tr}\,(e(\Theta)|H(\gamma))$ wohl keine ganze zahl, d.h. nicht in ${\cal O}$ gelegen
Using this we determine the trace of $e(\Theta)$ on ${\bf H}'(\gamma')$, $\gamma'\in{\cal E}({\bf H}')$, as follows. 
%(note that ${\bf E}(H)\subseteq {\cal S}$). 

\medskip

If $\gamma'\not\equiv\Phi^\vee(\Theta)\pmod{p^c}$ then $x=0$, hence, ${\rm tr}\,(e(\Theta)|{{\bf H}'(\gamma')})=0$.

\medskip

If $\gamma'\equiv \Phi^\vee(\Theta)\pmod{p^c}$ then we write for any $\mu\in{\cal E}({\bf H})$, $\mu\not\equiv\Theta\pmod{p^c}$
$$
\gamma'(T_\mu')=\Phi^\vee(\Theta)(T_\mu')+\delta_\mu',
$$ 
where $\delta_\mu'\in{\cal O}$.
%beachte, dass $\ganmma,\Theta$ ${\cal O}$-wertig sind, d.h. $\delta_\mu$ liegt ebenfalls in ${\cal O}$.
Since 
$v_p(\gamma'-\Phi^\vee(\Theta))\ge c$, equation (2') implies that $v_p(\gamma'-\Phi^\vee(\Theta))>c+l$, hence,
$$
v_p(\delta_\mu')> c+l.\leqno(4)
$$ 
Similarly, we write 
$$
\gamma'(T_{\mu'})=\Phi^\vee(\Theta)(T_{\mu'})+\delta_{\mu'}
$$
and find
$$
v_p(\delta_{\mu'})>c+l.\leqno(4')
$$
Recalling the definition of $\xi$ we obtain
\begin{eqnarray*}
x&=&\frac{\prod_{\mu\in{\cal E}({\bf H})\atop \mu\not\equiv \Theta\pmod{p^c}} \Phi^\vee(\Theta)(T_\mu')-\Phi^\vee(\mu)(T_\mu')+\delta_\mu'}{\prod_{\mu\in{\cal E}({\bf H})\atop \mu\not\equiv \Theta\pmod{p^c}} \Phi^\vee(\Theta)(T_\mu')-\Phi^\vee(\mu)(T_\mu')}\,
\frac{\prod_{\mu'\in{\cal E}({\bf H}')\atop \mu'\not\equiv \Phi^\vee(\Theta)\pmod{p^c}}\Phi^\vee(\Theta)(T_{\mu'})-\mu'(T_{\mu'})+\delta_{\mu'}}{\prod_{\mu'\in{\cal E}({\bf H}')\atop \mu'\not\equiv \Phi^\vee(\Theta)\pmod{p^c}}\Phi^\vee(\Theta)(T_{\mu'})-\mu'(T_{\mu'})}\\
%&=&\left(1+\frac{\delta}{\prod_{\mu\in{\cal E}({\bf H})\atop \mu\not\equiv \Theta\pmod{p^c}} \Phi^\vee(\Theta)(T_\mu')-\Phi^\vee(\mu)(T_\mu')}\right)\left(1+\frac{\delta'}{\prod_{\mu'\in{\cal E}({\bf H}')}\mu'(\Theta)-\mu'(T_{\mu'})}\right),\\
\end{eqnarray*}
%where $\delta$ is a sum of terms of the form $\prod_{\mu\in{\cal E}({\bf H})\atop \mu\not\equiv \Theta\pmod{p^c}} \Theta(T_\mu)-\mu(T_\mu)$ with at least one factor $\Theta(T_\mu)-\mu(T_\mu)$ replaced by $\delta_\mu$. 

Since $v_p(\Phi^\vee(\Theta)(T_\mu')-\Phi^\vee(\mu)(T_\mu'))<c$ for all $\mu\not\equiv\Theta\pmod{p^c}$ (cf. equation(3)) and since, 
$v_p(\delta_\mu')>c+l$ (cf. equation (4)) we obtain that the first factor on the right hand side of the above equation for $x$ is congruent to $1$ modulo $p^l$. 
Similarly, using equations (3') and (4') we find that the second factor on the right hand side of the above equation for $x$ is congruent to $1$ modulo $p^l$. 
Thus, we obtain $x\equiv 1\pmod{p^l}$, hence,
$$
{\rm tr}(e(\Theta)|{\bf H}'(\gamma'))\equiv{\rm dim}\,{\bf H}'(\gamma')\pmod{p^l}.
$$ 
Summing over all $\gamma'\in{\cal E}({\bf H}')$ finally yields
$$
{\rm tr}\,(e(\Theta))|{\bf H}')\equiv m_{\bf H}'(\Phi^\vee(\Theta),c)\pmod{p^l}.
$$

\medskip

d.) Quite analogous we compute the trace of $\Phi(e(\Theta))$ on ${\bf H}$. We note that
$$
\Phi(e(\Theta))=\frac{1}{\xi}\,\prod_{\mu\in{\cal E}({\bf H})\atop \mu\not\equiv \Theta\pmod{p^c}} \Phi(T_\mu')-\Phi^\vee(\mu)(T_\mu')\,
\prod_{\mu'\in{\cal E}({\bf H}')\atop \mu'\not\equiv \Phi^\vee(\Theta)\pmod{p^c}}\Phi(T_{\mu'})-\mu'(T_{\mu'})\;\in\frac{1}{\xi}{\cal H}_{\cal O}.
$$
Let $\gamma\in{\cal E}({\bf H})$. We choose a basis ${\cal B}(\gamma)$ of ${\bf H}(\gamma)$ such that any $T\in{\cal H}$ is upper triangular on ${\bf H}(\gamma)$: 
$$
{\cal D}_{{\cal B}(\gamma)}(T)=\left(\begin{array}{ccc}\gamma(T) &&*\\&\ddots&\\&&\gamma(T) \end{array} \right).
$$
Then $\Phi(e(\Theta))$ is represented on ${\bf H}(\gamma)$ by the matrix
$$
{\cal D}_{{\cal B}(\gamma)}(\Phi(e(\Theta)))=\left(\begin{array}{ccc}x &&*\\&\ddots&\\&&x \end{array} \right),
$$
where
$$
x=x(\gamma)=\frac{1}{\xi}\,\prod_{\mu\in{\cal E}({\bf H})\atop \mu\not\equiv \Theta\pmod{p^c}} \Phi^\vee(\gamma)(T_\mu')-\Phi^\vee(\mu)(T_\mu')\,
\prod_{\mu'\in{\cal E}({\bf H}')\atop \mu'\not\equiv \Phi^\vee(\Theta)\pmod{p^c}}\Phi^\vee(\gamma)(T_{\mu'})-\mu'(T_{\mu'})   \in\frac{1}{\xi}{\cal O}.
$$

\medskip

If $\gamma\not\equiv\Theta\pmod{p^c}$ then $x=0$, hence, ${\rm tr}\,(\Phi(e(\Theta))|{{\bf H}(\gamma)})=0$.

\medskip

If $\gamma\equiv \Theta\pmod{p^c}$ then equation (2) implies that $v_p(\gamma-\Theta)>c+l$. We write $\gamma(\Phi(T_\mu'))=\Theta(\Phi(T_\mu'))+\delta_\mu'$, or , 
equivalently,  
$$
\Phi^\vee(\gamma)(T_\mu')=\Phi^\vee(\Theta)(T_\mu')+\delta_\mu',
$$
where
$$
v_p(\delta_\mu')> c+l.\leqno(5)
$$ 
Similarly, we write
$$
\Phi^\vee(\gamma)(T_{\mu'})=\Phi^\vee(\Theta)(T_{\mu'})+\delta_{\mu'},
$$
where
$$
v_p(\delta_{\mu'})> c+l.\leqno(5')
$$ 
Recalling the definition of $\xi$ we obtain
\begin{eqnarray*}
x&=&\frac{\prod_{\mu\in{\cal E}({\bf H})\atop \mu\not\equiv \Theta\pmod{p^c}} \Phi^\vee(\Theta)(T_\mu')-\Phi^\vee(\mu)(T_\mu')+\delta_\mu'}{\prod_{\mu\in{\cal E}({\bf H})\atop \mu\not\equiv \Theta\pmod{p^c}} \Phi^\vee(\Theta)(T_\mu')-\Phi^\vee(\mu)(T_\mu')}\,
\frac{\prod_{\mu'\in{\cal E}({\bf H}')\atop \mu'\not\equiv \Phi^\vee(\Theta)\pmod{p^c}}\Phi^\vee(\Theta)(T_{\mu'})-\mu'(T_{\mu'})+\delta_{\mu'}}{\prod_{\mu'\in{\cal E}({\bf H}')\atop \mu'\not\equiv \Phi^\vee(\Theta)\pmod{p^c}}\Phi^\vee(\Theta)(T_{\mu'})-\mu'(T_{\mu'}).}\\
%&=&\left(1+\frac{\delta}{\prod_{\mu\in{\cal E}({\bf H})\atop \mu\not\equiv \Theta\pmod{p^c}} \Phi^\vee(\Theta)(T_\mu')-\Phi^\vee(\mu)(T_\mu')}\right)\left(1+\frac{\delta'}{\prod_{\mu'\in{\cal E}({\bf H}')}\mu'(\Theta)-\mu'(T_{\mu'})}\right),\\
\end{eqnarray*}
%where $\delta$ is a sum of terms of the form $\prod_{\mu\in{\cal E}({\bf H})\atop \mu\not\equiv \Theta\pmod{p^c}} \Theta(T_\mu)-\mu(T_\mu)$ with at least one factor $\Theta(T_\mu)-\mu(T_\mu)$ replaced by $\delta_\mu$. 
As above, using equations (3), (3'), (5), (5') we obtain $x\equiv 1\pmod{p^l}$, hence,
$$
{\rm tr}(\Phi(e(\Theta))|{\bf H}(\gamma))\equiv{\rm dim}\,{\bf H}(\gamma)\pmod{p^l}.
$$  
Summing over all $\gamma\in{\cal E}({\bf H})$ finally yields
$$
{\rm tr}\,(\Phi(e(\Theta))|{\bf H})\equiv m_{\bf H}(\Theta,c)\pmod{p^l}
$$
Hence, the Lemma is proven.

\bigskip

%{\bf Remark. } If $m, m'\in{\Bbb N}_0$ are natural numbers then $m\equiv m'\pmod{p^{\log_p M}}$ as defined in (1.1) is equivalent to $m\equiv m'\pmod{\lceil M\rceil_p}$, where the congruence on the right hand side is defined in the usual way, i.e. the difference $m-m'$ is divisible by $\lceil M\rceil_p$. 

%\bigskip

{\it Remark. } 1.) Since $m_{\bf H}(\Theta,c)$ and $m_{{\bf H}'}(\Phi^\vee(\Theta),c)$ are smaller than or equal to $M$ the congruences in (2.3) Lemma 
uniquely determine $m_{\bf H}(\Theta,c)$ and $m_{{\bf H}'}(\Phi^\vee(\Theta),c)$.

%denn: Seien $0\le a,b\le M$ mit $a\equiv b\pmod{p^l}$. O.E. sei $b\ge a$ also $0\le b-a \le M$. Da $v_p(b-a)\ge l$ folgt $v_p(b-a)\ge \lceil l\rceil$, d.h.
%$a\equiv b\pmod{p^{\lceil l\rceil}}$. Also existiert ein $h\in{\Bbb N}_0$ so dass $b-a=hp^{\lceil l \rceil}$ $\Rightarrow hp^{\lceil l \rceil}\le M$. Da $l=\log_p M$ ist aber
%$p^{\lceil l\rceil}\ge p^l=M$, d.h. $hM\le hp^{\lceil l \rceil}$. Es folgt $hM\le M$ woraus $h\le 0$ folgt da $M>0$. Da $h\in{\Bbb N}_0$ heisst das $h=0$ 
%also $a=b$.

2.) Since $c\ge m-(M+\frac{3}{2})\log_p M$, by increasing $m$ we can search for eigencharacters which are arbitrarily close to $\Phi^ \vee(\Theta)$ or $\Theta$; 
in doing so the denominators of $e(\Theta)$ will not grow unreasonably fast because $v_p(\xi)\le Mm$.
 
3.) The proof is not constructive, e.g. we do not obtain the value of $c$; in particular, we do not know whether we can choose $c=m-(M+\frac{3}{2})\log_p M$.

\bigskip

{\bf (2.4) } Using (2.3) Lemma we can give the following criterion for the existence of  the ``mod $p^c$ transfer'' $\Psi^{[c]}$ which makes
it possible to establish ``mod $p^c$ tansfer'' {via} a comparison of trace formulas.

\bigskip

{\bf Theorem. } {\it Assume that ${\rm dim}\,{\bf H},\,{\rm dim}\,{\bf H}' \le \frac{1}{2}M$ for some $M\in 2{\Bbb N}$. Assume that there is a rational number $s$ such that
$$
{\rm tr}\,(\Phi(T')|{\bf H})\equiv {\rm tr}\,(T'|{\bf H}')\pmod{p^s}
$$
for all $T'\in{\cal H}_{\cal O}'$. Then, for any $\Theta\in\hat{\cal H}_{\cal O}$ there is a rational number 
$c=c(\Theta)\ge \frac{s}{M}-(M+2)\log_p M$ such that
$$
m_{\bf H}(\Theta,c)= m_{{\bf H}'}(\Phi^\vee(\Theta),c).
$$ 
Hence, for any $\Theta\in{\cal E}({\bf H})$ there is an element $\Theta'\in{\cal E}({\bf H}')$ such that 
$$
\Theta'\equiv \Phi^\vee(\Theta)\pmod{p^{\frac{s}{M}-(M+2)\log_p M}},
$$
or, equivalently, there is a map on eigencharacters
$$
\Psi^{[c]}:\,{\cal E}({\bf H})\rightarrow{\cal E}({\bf H}')
$$
such that
$$
\Psi^{[c]}(\Theta)\equiv \Phi^\vee(\Theta)\pmod{p^{\frac{s}{M}-(M+2)\log_p M}}
$$
for all $\Theta\in{\cal E}({\bf H})$. 

}

\medskip

{\it Proof. } We set $m=\frac{s-\log_p M}{M}$. Let $\Theta\in\hat{\cal H}_{\cal O}$. According to (2.3) Lemma there is an element 
$e(\Theta)\in\frac{1}{\xi}{\cal H}'_{\cal O}$, where $v_p(\xi)\le Mm=s-\log_p M$, and a rational number
$$
c=c(\Theta)\ge m-(M+\frac{3}{2})\log_pM=\frac{s}{M}-(M+2)\log_p M
$$ 
(note that $\frac{1}{M}\le \frac{1}{2}$) such that ${\rm tr}\,(\Phi(e(\Theta))|{\bf H})\equiv m_{\bf H}(\Theta,c)$ and ${\rm tr}\,(e(\Theta)|{\bf H}')\equiv m_{{\bf H}'}(\Phi^\vee(\Theta),c)\pmod{p^{\log_p M}}$. Since $v_p(\xi)\le Mm=s-\log_p M$, the assumption of the Theorem implies 
$$
{\rm tr}\,(\Phi(e(\Theta))|{\bf H})\equiv{\rm tr}\,(e(\Theta)|{\bf H}')\pmod{p^{\log_p M}},
$$ 
hence, 
$$
m_{\bf H}(\Theta,c)\equiv m_{{\bf H}'}(\Phi^\vee(\Theta),c)\pmod{p^{\log_p M}}.
$$
Since $m_{\bf H}(\Theta,c)$ and $m_{{\bf H}'}(\Phi^\vee(\Theta),c)$ are natural numbers which are smaller than $\frac{1}{2}M\le p^{\log_p M}$ this implies
$$
m_{\bf H}(\Theta,c)= m_{{\bf H}'}(\Phi^\vee(\Theta),c).
$$ 
In particular, for any $\Theta\in{\cal E}({\bf H})$ there is an element 
$\Theta'\in{\cal E}({\bf H}')$ such that $\Theta'\equiv \Phi^\vee(\Theta)\pmod{p^c}$. Hence, the proof is complete. 

\medskip

%{\it Remark. } The Theorem may be seen as a mod $p^c$ analog of Lemma 16.1.1 in [J-L]

\bigskip

{\bf (2.5) Corollary. }{\it Assume that ${\rm dim}\,{\bf H},\,{\rm dim}\,{\bf H}' < \frac{1}{2}M$ for some $M\in 2{\Bbb N}$. Assume that there is a rational number $s$ such that
$$
{\rm tr}\,(T|{\bf H})\equiv {\rm tr}\,(T|{\bf H}')\pmod{p^s}
$$
for all $T\in{\cal H}_{\cal O}$. Then, for any $\Theta\in\hat{\cal H}_{\cal O}$ there is a rational number 
$c\ge \frac{s}{M}-(M+2)\log_p M$ such that
$$
m_{\bf H}(\Theta,c)= m_{{\bf H}'}(\Theta,c).
$$ 
Hence, for any $\Theta\in{\cal E}({\bf H})$ there is an eigencharacter $\Theta'\in{\cal E}({\bf H}')$
such that 
$$
\Theta'\equiv \Theta\pmod{p^{\frac{s}{M}-(M+2)\log_p M}}.
$$
}

\medskip

{\it Proof. } This is the special case ${\cal H}={\cal H}'$ and $\Phi={\rm id}$.

\bigskip

\section{Systems of higher congruences in the finite slope case}

Combining (1.9) Proposition with (2.5) Corollary we would obtain a trace criterion for the existence of $p$-adic 
continuous families of eigencharacters passing through a given eigencharacter. However, in this section we will construct certain elements which behave 
like an approximate idempotent attached to the slope $\alpha$ subspace of a ${\cal H}$-module ${\bf H}$; using this we then will obtain a 
criterion for the existence of finite slope $p$-adic continuous families of eigencharacters.

\bigskip

{\bf (3.1) Slope subspaces. } As in (1.2) we set ${\cal H}=\bar{\Bbb Q}_p[T_\ell,\,\ell\in I]$ and ${\cal H}_{\cal O}={\cal O}[T_\ell,\,\ell\in I]$ (countably generated). We fix an element $T\in{\cal H}_{{\cal O}}$. For any finite dimensional ${\cal H}$-module ${\bf H}$ 
which is defined over ${\cal O}$ with respect to the lattice ${\bf H}_{\cal O}$ we denote by ${\bf H}^\alpha$, $\alpha\in{\Bbb Q}_{\ge 0}$, the slope $\alpha$ subspace of ${\bf H}$ with respect to the operator $T$, i.e. 
$$
{\bf H}^\alpha=\bigoplus_{\mu\in{\cal O}\atop v_p(\mu)=\alpha} {\bf H}(\mu)
$$
where ${\bf H}(\mu)$ is the generalized eigenspace attached to $T$ and $\mu$ and we set
$$
{\bf H}^{\le\alpha}=\bigoplus_{\beta\le\alpha} {\bf H}^\beta.
$$ 
We denote by $\Phi_{\bf H}\subset {\cal O}$ the set of eigenvalues of $T$ acting on ${\bf H}$ and by $\Phi_{\bf H}^{\le\alpha}\subseteq \Phi_{\bf H}$ the subset of all eigenvalues $\gamma$ of $T$ 
satisfying $v_p(\gamma)\le \alpha$. Hence, ${\bf H}^{\le\alpha}=\bigoplus_{\gamma\in \Phi_{\bf H}^{\le\alpha}}{\bf H}(\gamma)$.

%VVVVVVVVVVVVVVVVV
%
%We set 
%$E_{\bf H}^{\le\alpha}={\Bbb Q}(\Phi_{\bf H}^{\le\alpha})$, i.e. $E_{\bf H}^{\le\alpha}$ is the field obtained by adjoining to ${\Bbb Q}$ all eigenvalues 
%of $p$-adic value less than or equal to $\alpha$. 
%
%Hence, $[E_{{\bf H},{\bf H}'}^{\le\alpha}:{\Bbb Q}]\le [E_{\bf H}^{\le\alpha}:{\Bbb Q}][E_{{\bf H}'}^{\le\alpha}:{\Bbb Q}]\le M(\alpha)^2$ and $e^{\le\alpha}\le M(\alpha)^2$. 

%AAAAAAAAAAAAAAAAAAA

%We denote by $\Xi_H^{\le \alpha}$ the characteristic polynomial of $T$ acting on $H^{\le \alpha}$ and $E_H^{\le\alpha}$ is the field obtained by adjoining all 
%roots of $\Xi_H^{\le\alpha}$ to ${\Bbb Q}$. By $\Phi_H$ we denote the set of all roots of the characteristic polynomial $\Xi_H$ of $T$ acting on $H$ and 
%$\Phi_H^{\le\alpha}$ is the set of all roots of $\Xi_H^{\le\alpha}$. Thus, $\Phi_H^{\le\alpha}$ consists of all complex numbers $\gamma$ such that 
%$v_p(\gamma)\le\alpha$ and $H(\gamma)\not=0$.

\bigskip

{\bf (3.2) Construction of approximate idempotents. } Let ${\bf H},{\bf H}'$ be arbitrary ${\cal H}$-modules, finite dimensional and defined over ${\cal O}$ and select a slope $\alpha\in{\Bbb Q}_{\ge 0}$. 
%
%VVVVVVVVVVV
%We denote by
%$E_{{\bf H},{\bf H}'}^{\le\alpha}=E_{\bf H}^{\le\alpha}E_{\bf H'}^{\le\alpha}$ the compositum of the fields $E_{\bf H}^{\le\alpha}$ and 
%$E_{{\bf H}'}^{\le\alpha}$ and $e^{\le\alpha}=e_{{\bf H},{\bf H}'}^{\le\alpha}$ is the ramification degree over $p$ of the prime ideal ${\mathfrak p}$ 
%in $E_{{\bf H},{\bf H}'}^{\le\alpha}$ which corresponds to $v_p$. 
%AAAAAAAAAAAAAAAA
%
For any $\gamma\in\Phi_{\bf H}$ we choose a basis 
${\cal B}={\cal B}(\gamma)$ of ${\bf H}(\gamma)$ such that the representing matrix 
${\cal D}_{\cal B}(T|_{{\bf H}(\gamma)})$ of $T$ on ${\bf H}(\gamma)$ is upper triangular:
$$
{\cal D}_{\cal B}(T|_{{\bf H}(\gamma)})=
\left(
\begin{array}{ccc}
\gamma&&*\\
&\ddots&\\
&&\gamma\\
\end{array}
\right).\leqno(1)
$$
%where $"*"\in {\Bbb Z}$ (e.g. the Jordan canonical form). 
Similarly, for any $\gamma\in\Phi_{{\bf H}'}$ we choose a basis ${\cal B}'={\cal B}'(\gamma)$ of 
${\bf H}'(\gamma)$ such that the matrix representing $T$ on ${\bf H}'(\gamma)$ is upper triangular: 
$$
{\cal D}_{{\cal B}'}(T|_{{\bf H}'(\gamma)})=
\left(
\begin{array}{ccc}
\gamma&&*\\
&\ddots&\\
&&\gamma\\
\end{array}
\right).
$$
%where $"*"\in{\Bbb Z}$. 
We define the element 
$$
{\mathbf e}^{\le\alpha}={\mathbf e}_{{\bf H},{\bf H}'}^{\le\alpha}=1-\prod_{\mu\in \Phi_{\bf H}^{\le\alpha}\cup\Phi_{{\bf H}'}^{\le\alpha}}\frac{T-\mu}{-\mu}\in\bar{\Bbb Q}_p[T].
$$
Clearly, ${\mathbf e}^{\le\alpha}=p^{\le\alpha}(T)$, where the polynomial $p^{\le\alpha}(X)$ is given by
$$
p^{\le\alpha}=p_{{\bf H},{\bf H}'}^{\le\alpha}=1-\prod_{\mu\in \Phi_{\bf H}^{\le\alpha}\cup\Phi_{{\bf H}'}^{\le\alpha}}\frac{X-\mu}{-\mu}\in \bar{\Bbb Q}_p[X].
$$
We want to collect some properties of ${\mathbf e}^{\le\alpha}$ and $p^{\le\alpha}$. To this end we let $M(\alpha)\in 2{\Bbb N}$, $\alpha\in{\Bbb Q}_{\ge 0}$, be a collection of natural numbers such that 
$$
{\rm dim}\,{\bf H}^{\le\alpha}\le \frac{1}{2}M(\alpha)\quad\mbox{and}\quad{\rm dim}\,{\bf H}'^{\le \alpha}\le \frac{1}{2}M(\alpha)\leqno(2)
$$
for all $\alpha\in{\Bbb Q}_{\ge 0}$. For an arbitrary polynomial $p=\sum_{i\ge 0}a_iX^i\in\bar{\Bbb Q}_p[X]$ we define its slope as 
$$
{\bf S}(p)={\rm sup}\,\{s\in{\Bbb Q}\cup\{-\infty\}: v_p(a_i)\ge si\,\mbox{for all $i\ge 0$}\},
$$
hence, ${\bf S}(p)>-\infty$ implies $v_p(a_0)\ge 0$. Easy calculation shows that
$$
{\bf S}(pq)\ge {\rm min}\{{\bf S}(p),{\bf S}(q)\}\leqno(3a)
$$
and 
$$
{\bf S}(p+q)\ge {\rm min}\{{\bf S}(p),{\bf S}(q)\}.\leqno(3b)
$$

\medskip

{\bf (3.3) Lemma. }{\it 1.) For any $\gamma\in\Phi_{\bf H}- \Phi_{\bf H}^{\le \alpha}$ we have
$$
{\cal D}_{\cal B}({\mathbf e}_{{\bf H},{\bf H}'}^{\le\alpha}|_{{\bf H}(\gamma)})=
\left(
\begin{array}{ccc}
\zeta&&*\\
&\ddots&\\
&&\zeta\\
\end{array}
\right)
$$
where $\zeta\in{\cal O}$ satisfies $v_p({\zeta})\ge \frac{2}{M(\alpha+1)}$. The analogous statement holds for $\gamma\in \Phi_{{\bf H}'}- \Phi_{{\bf H}'}^{\le \alpha}$.

2.) For any $\gamma\in\Phi_{\bf H}^{\le \alpha}$ we have
$$
{\cal D}_{\cal B}({\mathbf e}_{{\bf H},{\bf H}'}^{\le\alpha}|_{{\bf H}(\gamma)})=
\left(
\begin{array}{ccc}
1&&*\\
&\ddots&\\
&&1\\
\end{array}
\right).
$$
Again, the analogous statement holds for $\gamma\in \Phi_{{\bf H}'}^{\le \alpha}$.

3.) $$
{\rm deg}\,p_{{\bf H},{\bf H}'}^{\le\alpha}=|\Phi_{\bf H}^{\le\alpha}|+|\Phi_{{\bf H}'}^{\le\alpha}|\le {\rm dim}\, {\bf H}^{\le\alpha}+{\rm dim}\,{{\bf H}'}^{\le\alpha}\le M(\alpha).
$$

4.) $$
p_{{\bf H},{\bf H}'}^{\le\alpha}(0)=0
$$

5.) $${\bf S}(p_{{\bf H},{\bf H}'}^{\le\alpha})\ge -\alpha.
$$
In particular, $(p^{\le\alpha})^L=\sum_{h\ge L} b_hX^h$, where $v_p(b_h)\ge -h\alpha$ for all $h\ge L$.

}

\medskip

{\it Proof. } 1.) Let $\gamma\in\Phi_{\bf H}-\Phi_{\bf H}^{\le \alpha}$. Equation (1) implies that with respect to ${\cal B}={\cal B}(\gamma)$
$$
{\cal D}_{\cal B}({\mathbf e}^{\le\alpha}|_{{\bf H}(\gamma)})=
\left(
\begin{array}{ccc}
\zeta&&*\\
&\ddots&\\
&&\zeta\\
\end{array}
\right)
$$
where 
$$
\zeta=1-\prod_{\mu\in \Phi_{\bf H}^{\le\alpha}\cup\Phi_{{\bf H}'}^{\le\alpha}} \frac{\gamma-\mu}{-\mu}.
$$
Let $\mu\in\Phi_{\bf H}^{\le\alpha}\cup \Phi_{{\bf H}'}^{\le\alpha}$ be arbitrary. We distinguish two cases. First, if $\gamma\not\in\Phi_{\bf H}^{\le\alpha+1}$, 
i.e. $v_p(\gamma)>\alpha+1$, then we obtain $v_p(\gamma)>v_p(\mu)+1\ge v_p(\mu)+2/M(\alpha+1)$. 
%da $M(\alpha+1)\ge 2$ 
Second, if $\gamma\in\Phi_{\bf H}^{\le\alpha+1}$, then $\gamma$ is an eigenvalue of $T$ acting on ${\bf H}^{\le\alpha+1}$, hence, it is a root of the characteristic 
polynomial of $T$ acting on ${\bf H}^{\le\alpha+1}$
%denn $\gamma$ ist ein EW von $T$ auf ${\bf H}^{\le\alpha+1}$ (da ${\bf H}^{\le\alpha}=\bigoplus_{\gamma\in \Phi_{\bf H}^{\le\alpha}}{\bf H}(\gamma)$)
which has degree ${\dim}\,{\bf H}^{\le\alpha+1}\le \frac{1}{2} M(\alpha+1)$. We deduce that $\gamma$ is contained in an extension of ${\Bbb Q}_p$ of degree less than or 
equal to $\frac{1}{2} M(\alpha+1)$ which implies that $v_p(\gamma)\in\frac{2}{M(\alpha+1)}{\Bbb N}_0$. 
%da der Verzweigungsindex von $p$ in $K/{\Bbb Q}$ kleiner oder gleich $[K:{\Bbb Q}]$ ist
On the other hand, since $\mu\in\Phi_{\bf H}^{\le\alpha+1}\cup \Phi_{{\bf H}'}^{\le\alpha+1}$ we obtain quite similarly that $v_p(\mu)\in\frac{2}{M(\alpha+1)}{\Bbb N}_0$. 
Since $\gamma\not\in\Phi_{\bf H}^{\le\alpha}$ we know that $v_p(\gamma)>v_p(\mu)$, hence, $v_p(\gamma)\ge v_p(\mu)+2/M(\alpha+1)$. 
Thus, in both cases we find $v_p(\frac{\gamma}{\mu})\ge 2/M(\alpha+1)$. Since
$$
\zeta=1-\prod_{\mu\in \Phi_{\bf H}^{\le\alpha}\cup\Phi_{{\bf H}'}^{\le\alpha}} 1-\frac{\gamma}{\mu}
$$
is a sum of products of the form $\pm\prod_\mu \frac{\gamma}{\mu}$, where $\mu$ runs over a non-empty subset of $\Phi_{\bf H}^{\le\alpha}\cup\Phi_{{\bf H}'}^{\le\alpha}$ 
(the summand "$1$" cancels), we deduce that $v_p(\zeta)\ge 2/M(\alpha+1)$.

2.) Immediate by the definition of ${\bf e}^{\le\alpha}$.

3.) and 4.) Clear

5.) By definition it is immediate that ${\bf S}(\frac{T-\mu}{-\mu})={\bf S}(\frac{T}{-\mu}+1)= -v_p(\mu)$. All $\mu$ appearing in the definition of 
${\bf e}^{\le\alpha}$ satisfy $v_p(\mu)\le\alpha$; hence, using equation (3a) and (3b) we deduce ${\bf S}(p^{\le \alpha})\ge {\rm min}\{0,-\alpha\}=-\alpha$. 
The second statement follows because $X$ divides $p^{\le\alpha}$ and ${\bf S}((p^{\le\alpha})^L)\ge{\bf S}(p^{\le\alpha})$. This finishes the proof of the Lemma.

%We set $\Phi_{p,k,k'}^{\le\alpha}=\Phi_{p,k}^{\le\alpha}\cup\Phi_{p,k'}^{\le\alpha}$. Straightforward computation shows that
%$$
%c_i=\sum_{\{\mu_1,\ldots,\mu_{{\rm deg}\,p-i}\}\subset\Phi_{p,k,k'}^{\le\alpha}} \frac{\mu_1\cdot\ldots\cdot\mu_{{\rm deg}\,p-i}}{\prod_{\mu\in\Phi_{p,k,k'}^{\le\alpha}}\mu}.
%$$
%Since $v_p(\mu)\le\alpha$ for all $\mu\in\Phi_{p,k,k'}^{\le\alpha}$ this implies the first claim. To prove the second claim we note that
%$$
%b_h=\sum_{(i_1,\ldots,i_L)\atop i_1+\cdots+i_L=h} c_{i_1}\cdot\ldots\cdot c_{i_L}.
%$$
%Since $v_p(c_{i_k})\ge -i_k\alpha$ we find 
%$$
%v_p(c_{i_1}\cdot\ldots\cdot c_{i_L})=\sum_{k=1}^L v_p(c_{i_k})\ge \sum_{k=1}^L -i_k\alpha=-\alpha\sum_{k=1}^L i_k=-h\alpha.
%$$ 
%Thus, we obtain $v_p(b_h)\ge -h\alpha$. This completes the proof of the Lemma.

\bigskip

%In particular, we find using the above Lemma 

%{\bf Corollary. } {\it 1.) For any $\gamma\in\Phi_{p,k}- \Phi_{p,k}^{\le \alpha}$ we have
%$$
%{\cal D}_{\cal B}({\mathbf e}_{k,k'}^{\le\alpha}|_{H(\gamma)})=
%\left(
%\begin{array}{ccc}
%\zeta&&*\\
%&\ddots&\\
%&&\zeta\\
%\end{array}
%\right)
%$$
%where $\zeta\in {\cal O}$ and $v_p({\zeta})\ge \frac{1}{2M(\alpha)}$. An analogous statement holds for $\gamma\in \Phi_{p,k'}- \Phi_{p,k'}^{\le \alpha}$.

%2.) $$
%{\rm deg}\,p_{k,k'}^{\le\alpha}\le 2M(\alpha).
%$$
%}

{\bf (3.4) Proposition. }{\it For any pair of finite dimensional ${\cal H}$-modules ${\bf H}$, ${\bf H}'$ which are defined over ${\cal O}$ and satisfy equation (2) 
and for any $T\in{\cal H}_{\cal O}$ we have
$$
{\rm tr}\,(T ({{\bf e}_{{\bf H},{\bf H}'}^{\le\alpha}})^L|{\bf H})\equiv {\rm tr}\,(T|{{\bf H}^{\le\alpha}})\pmod{p^{\frac{2L}{M(\alpha+1)}}}.
$$
The same congruence holds for ${\bf H}'$ in place of ${\bf H}$.
}

\medskip

{\it Proof. } Let $\gamma\in\Phi_{\bf H}$ and let ${\cal B}$ be a basis of ${\bf H}(\gamma)$ such that equation (1) holds. (3.3) Lemma implies that 
$T({\bf e}^{\le\alpha})^L$ is represented on ${\bf H}(\gamma)$ by the matrix
$$
\left(
\begin{array}{ccc}
\gamma\zeta^L&&*\\
&\ddots&\\
&&\gamma\zeta^L\\
\end{array}
\right),
$$
where $\zeta\equiv 0\pmod{\frac{2}{M(\alpha+1)}}$ if $v_p(\gamma)>\alpha$ and $\zeta=1$ if $v_p(\gamma)\le \alpha$. Since $v_p(\gamma)\ge 0$ ($T\in{\cal H}_{\cal O}$) 
this implies the claim. The same argument works if we replace ${\bf H}$ by ${\bf H}'$, hence the proof is complete.

\bigskip

{\bf (3.5) Remark. } In particular, we obtain
$$
\lim_{L\rightarrow\infty} \,{\rm tr}\,(T ({{\bf e}^{\le\alpha}})^L|{\bf H})= {\rm tr}\,(T|{{\bf H}^{\le\alpha}})
$$
and the same holds if we replace ${\bf H}$ by ${\bf H}'$. Thus, ${\bf e}^{\le\alpha}$ behaves like an approximate idempotent attached to the slope 
$\le\alpha$-subspaces of ${\bf H}$ and ${\bf H}'$. On the other hand, ${\bf e}^{\le\alpha}$ is not universal, i.e. it not only depends on $\alpha$ but also on the pair 
of modules ${\bf H}$ and ${\bf H}'$.

\bigskip

{\bf (3.6) A criterion for the existence of $p$-adic continuous families of finite slope. }  We give the synthesis of our results obtained so far. We use the notations 
from (1.8), i.e. ${\bf G}/{\Bbb Q}$ is a reductive algebraic group with maximal split torus ${\bf T}$. Let ${\cal R}\subseteq X({\bf T})$ be a subset and let $({\bf H}_\lambda)$, $\lambda\in{\cal R}$, be a family of finite dimensional ${\cal H}$-modules. From now on we will always 
assume that the following assumptions hold for the family $({\bf H}_\lambda)$:

\begin{enumerate}

\item Any ${\bf H}_\lambda$ is defined over ${\cal O}$ 

\item There are numbers $M(\alpha)\in 2{\Bbb N}$, $\alpha\in{\Bbb Q}_{\ge 0}$, such that
$$
{\rm dim}\,{\bf H}^{\le\alpha}_\lambda\le \frac{1}{2}M(\alpha) \quad\mbox{for all $\alpha\in{\Bbb Q}_{\ge 0}$ and all $\lambda\in {\cal R}$}.
$$

\end{enumerate}

\medskip

We select a slope $\alpha\in{\Bbb Q}_{\ge 0}$ and we put ${\bf e}_{\lambda,\lambda'}={\bf e}_{{\bf H}_\lambda,{\bf H}_{\lambda'}}^{\le \alpha}$. 

\bigskip

{\bf (3.7) Proposition. } {\it Let $({\bf H}_\lambda)_{\lambda\in{\cal R}}$ be a family of ${\cal H}$-modules. Assume that there 
is a collection of rational numbers $a'=a'(\alpha),a=a(\alpha)\in{\Bbb Q}_{>0}$ and $b=b(\alpha)\in{\Bbb Q}_{\le 0}$ ($\alpha\in{\Bbb Q}_{\ge 0}$) with 
the following property: i) $a'/M(\alpha+1)$, $a$ and $b$ are decreasing in $\alpha$ ii) for any $\alpha\in{\Bbb Q}_{\ge 0}$ and any pair 
$\lambda,\lambda'\in X({\bf T})$ with $\lambda\equiv \lambda'\pmod{(p-1)p^m X({\bf T})}$ there is a natural number $L\ge a'(m+1)$ such that
$$
{\rm tr}\,({\bf e}_{\lambda,\lambda'}^L T|{\bf H}_\lambda)\equiv {\rm tr}\,({\bf e}_{\lambda,\lambda'}^L T|{\bf H}_{\lambda'})\pmod{p^{a(m+1)+b}}\leqno(\dag)
$$
for all $T\in{\cal H}_{\cal O}$. Then 

\medskip

1.) ${\rm dim}\,{\bf H}_\lambda^\alpha$ is locally constant as a function of $\lambda$, i.e. there is $D=D(\alpha)\in{\Bbb N}$ only 
depending on $\alpha$ such that $\lambda\equiv\lambda'\pmod{(p-1)p^D}$ implies ${\rm dim}\,{\bf H}_\lambda^\alpha={\rm dim}\,{\bf H}_{\lambda'}^\alpha$.

\medskip

2.) Any $\Theta\in{\cal E}({\bf H}_{\lambda_0}^{\alpha})$ fits in a $p$-adic continuous family of eigencharacters
of slope $\alpha$, i.e. there is a family $(\Theta_\lambda)_{\lambda\in{\cal R}}$ such that

\begin{enumerate}

\item $\Theta_\lambda\in{\cal E}({\bf H}_\lambda^{\alpha})$

\item $\Theta_{\lambda_0}=\Theta$

\item $\lambda\equiv\lambda'\pmod{(p-1)p^m X({\bf T})}$ implies 
$$
\Theta_\lambda\equiv\Theta_{\lambda'}\pmod{p^{{\sf a}(m+1)+{\sf b}}},
$$
where ${\sf a}=\frac{1}{M(\alpha)}{\rm min}\,(a,\frac{2 a'}{M(\alpha+1)})$ and ${\sf b}=\frac{b}{M(\alpha)}-(M(\alpha)+2)\log_p (M(\alpha))$.

\end{enumerate}

}

\medskip

{\it Proof. }  We select $\alpha\in{\Bbb Q}_{\ge 0}$ and we proceed in steps.

\medskip

a.) We let $\lambda\equiv\lambda'\pmod{(p-1)p^m X({\bf T})}$. We choose $L\ge a'(m+1)$ 
such that equation ($\dag$) holds. Together with (3.4) Proposition we obtain for all $T\in{\cal H}_{\cal O}$
$$
{\rm tr}\,(T|{\bf H}_\lambda^{\le\alpha})\equiv{\rm tr}\,({\bf e}_{\lambda,\lambda'}^LT|{\bf H}_\lambda)
\equiv {\rm tr}\,({\bf e}_{\lambda,\lambda'}^LT|{\bf H}_{\lambda'})
\equiv {\rm tr}\,(T|{\bf H}_{\lambda'}^{\le\alpha})\pmod{p^s},\leqno(4)
$$
where $s=\bar{a}(m+1)+b$ with $\bar{a}={\rm min}\,(\frac{2 a'}{M(\alpha+1)},a)$ (note that $b\le 0$). We note that our 
assumptions imply that $s$ is decreasing in $\alpha$. 

\medskip

b.) We next show that ${\rm dim}\,{\bf H}_\lambda^{\le\alpha}$ is locally constant in $\lambda$. To this end we set $D=D(\alpha)
=\frac{\log_p M(\alpha)-b}{\bar{a}}-1$ and we let $\lambda\equiv \lambda'\pmod{(p-1)p^D X({\bf T})}$. Since
$$
{\rm dim}\,{\bf H}_\lambda^{\le\alpha}={\rm tr}\,({\bf 1}|{{\bf H}_\lambda^{\le\alpha}})
\equiv{\rm tr}\,({\bf 1}|{{\bf H}_{\lambda'}^{\le\alpha}})={\rm dim}\,{\bf H}_{\lambda'}^{\le\alpha}  \pmod{p^{\bar{a}(D+1)+b}}.
$$
and since $p^{\bar{a}(D+1)+b}= M(\alpha)> {\rm dim}\, {\bf H}_\lambda^{\le\alpha}$, ${\rm dim}\, {\bf H}_{\lambda'}^{\le\alpha}$ we obtain
${\rm dim}\,{\bf H}_\lambda^{\le\alpha}={\rm dim}\,{\bf H}_{\lambda'}^{\le\alpha}$. Since $D(\alpha)$ is increasing in $\alpha$ the congruence 
$\lambda\equiv \lambda'\pmod{(p-1)p^D X({\bf T})}$ even implies 
${\rm dim}\,{\bf H}_\lambda^{\le\beta}={\rm dim}\,{\bf H}_{\lambda'}^{\le\beta}$ for all $\beta\le \alpha$.

\medskip

c.) Let $\lambda_i\in X({\bf T})$ and let $0\le\alpha_1<\cdots<\alpha_s\le\alpha$ be the non-trivial slopes appearing in 
${\bf H}_{\lambda_i}^{\le\alpha}$. Part b.) implies that for any $\lambda\in \lambda_i+p^{D(\alpha)} X({\bf T})$ the non trivial slopes appearing in ${\bf H}_\lambda^{\le\alpha}$
again are $\alpha_1<\cdots<\alpha_s$ (note that $D(\alpha)$ is increasing in $\alpha$). 
%hence, part b) implies that ${\rm dim}\,{\bf H}_\lambda^{\le\beta}={\rm dim}\,{\bf H}_{\lambda_i}^{\le\beta}$ for all $\beta\le\alpha$. 
Since $X({\bf T})\cong {\Bbb Z}^n$ is covered by 
finitely many cosets $\lambda_i+p^D X({\bf T})$, $i=1,\ldots,s$ we can select non negative rational numbers $0\le \alpha_1<\cdots<\alpha_s\le\alpha$ such
 that  for any $\lambda\in{\cal R}$ the inequality ${\bf H}_\lambda^\alpha\not=0$ implies that $\alpha$ is one of the $\alpha_i$.

\medskip

d.) We denote by $\beta$ the largest of the numbers $\alpha_1<\cdots<\alpha_s$ which is strictly smaller than $\alpha$. We obtain ${\rm dim}\,{\bf H}_\lambda^{\alpha}={\rm dim}\,{\bf H}_\lambda^{\le\alpha}-
{\rm dim}\,{\bf H}_\lambda^{\le\beta}$ for any $\lambda\in{\cal R}$. Part b.) then implies that 
${\rm dim}\,{\bf H}_\lambda^{\alpha}={\rm dim}\,{\bf H}_{\lambda'}^{\alpha}$ if $\lambda\equiv\lambda'\pmod{(p-1)p^{D(\alpha)} X({\bf T})}$.

\medskip

e.) Finally, since $s$ is decreasing in $\alpha$ equation (4) still holds if we replace 
"$\le\alpha$" by "$\le\beta$". Hence, by substracting we obtain
$$
{\rm tr}\,(T|{\bf H}_\lambda^{\alpha})\equiv {\rm tr}\,(T|{\bf H}_{\lambda'}^{\alpha})\pmod{p^s}.
$$
Thus, (2.5) Corollary implies that for all $\Theta\in\hat{\cal H}_{\cal O}$
$$
m_{{\bf H}_\lambda^{\alpha}}(\Theta,c)=m_{{\bf H}_{\lambda'}^{\alpha}}(\Theta,c)
$$
for some $c=c(\Theta)\ge {\sf a}(m+1)+{\sf b}$ with ${\sf a}$ and ${\sf b}$ as in the Proposition. (1.9) Proposition now
implies the claim. This completes the proof of the Proposition.

\bigskip

{\it Remark. } Since ${\bf H}_\lambda^{\le\beta}\subseteq{\bf H}_\lambda^{\le\alpha}$ if $\beta\le\alpha$ it is natural that on ${\bf H}_\lambda^{\le\alpha}$ 
weaker congruences for eigenvalues of operators hold, i.e. $a$ and $b$ are decreasing in $\alpha$. In particular, smaller power of ${\bf e}_{\lambda,\lambda'}$ should be sufficient, i.e. $a'/M(\alpha+1)$ also should be decreasing in $\alpha$.

\bigskip

\centerline{\bf \Large B. Cohomology of the Siegel upper half plane}

\bigskip

In this second part we consider an example: we will show that the family of cohomology groups of the Siegel upper half plane with 
coefficients in the irreducible representation of varying highest weight $\lambda$ satisfies equation 
$(\dag)$ in (3.7). As a consequence we obtain the existence of $p$-adic 
continuous families of Siegel eigenforms of finite slope $\alpha$. We start by fixing some notation.

\bigskip

\section{Notations}

{\bf (4.1) The symplectic group. } From now on we set ${\bf G}={\bf GSp}_{2n}$. Hence, for any ${\Bbb Z}$-algebra $K$ we have 
$$
{\bf G}(K)=\{g\in {\rm GL}_{2n}(K):\,g^t\Mat{}{I_n}{-I_n}{}g=\nu(g)\Mat{}{I_n}{-I_n}{}\;\mbox{for some $\nu(g)\in K^*$}\}.
$$
The multiplier defines a character $\nu:\,{\bf G}\rightarrow{\Bbb G}_m$ and the derived group of ${\bf G}$ is the symplectic group ${\bf G}^0={\bf Sp}_{2n}$ which is the kernel of $\nu$, i.e.  
$$
{\bf G}^0(K)=\{g\in{\bf GSp}_{2n}(K):\,\nu(g)=1\}.
$$
Thus, ${\bf G}^0(K)$ consists of all matrices $g=\Mat{A}{B}{C}{D}$ satisfying
$$
A^tB=B^tA,\quad C^tD=D^tC,\quad A^tD-B^tC=1\leqno(1)
$$
%cf. [Hein, Struktur und darstellungstheorie der klass gruppen], chapter I.5.3, Satz 3, p. 85
%s. auch [Freitag, Siegel modular forms], chapter I, Bemrekung 1.2, p. 25; hier sind aber $B$ und $C$ vertauscht gegenuber [Hein]
%
and is simply connected. 
%cf. [Malle, testerman, Lie groups of finte type], table 9.2, p. 72
We set ${\bf Z}={\bf Z}_{\bf G}$, hence, ${\bf Z}(K)=\{\lambda\cdot I_{2n},\;\lambda\in K^*\}$ and we denote by ${\bf T}$ the maximal ${\Bbb Q}$-split torus in ${\bf G}$ whose $K$-points are given by 
$$
{\bf T}(K)=\{{\rm diag}(\alpha_1,\ldots,\alpha_n,\frac{\nu}{\alpha_1},\ldots,\frac{\nu}{\alpha_n}),\,\alpha_i\in K^*,\,\nu\in K^*\}.\leqno(2)
$$
%$$
%{\bf T}(K)=\{\left(
%\begin{array}{cccccc}
%\alpha_1&$$$$\\
%&\ddots&&&&\\
%&&\alpha_n&&&\\
%&&&\nu/\alpha_1&&\\
%&&&&\ddots&\\
%&&&&&\nu/\alpha_n
%\end{array}
%\right),\,\alpha_i\in K^*,\,\nu\in K^*
%\}.
%$$
The intersection ${\bf T}^0={\bf T}\cap {\bf G}$ is a maximal split torus in ${\bf G}^0$ and 
$$
{\bf T}^0(K)=\{{\rm diag}(\alpha_1,\ldots,\alpha_n,\alpha_1^{-1},\ldots,\alpha_n^{-1}),\,\alpha_i\in K^*\}.
$$ 
We denote by ${\mathfrak g}$ resp. ${\mathfrak g}^0$ resp. ${\mathfrak h}$ resp. ${\mathfrak h}^0$ the complexified Lie Algebra of ${\bf G}/{\Bbb Q}$ resp. ${\bf G}^0/{\Bbb Q}$ resp. ${\bf T}$ resp. ${\bf T }^0$. Thus,
$$
{\mathfrak g}^0=\{\Mat{A}{B}{C}{D}\in M_{2n}({\Bbb C}):\,A=-D^t, B=B^t, C=C^t\}
$$
%cf. [Hein, Struktur und darstellungstheorie der klass gruppen], chapter II.2.2, p. 110
and 
$$
{\mathfrak h}^0=\{{\rm diag}(a_1,\ldots,a_n,-a_1,\ldots,-a_n),\,a_i\in{\Bbb C}\}.
$$
We denote by $\Phi=\Phi({\bf T}^0,{\mathfrak g}^0)$ the root system of ${\bf G}^0$ with respect to ${\bf T}^0$; explicitly, the roots 
$\alpha\in\Phi$, a generator $X_\alpha$ of the corresponding root space ${\mathfrak g}_\alpha\le{\mathfrak g}^0$ and the $1$-parameter subgroup 
$\exp tX_\alpha$ are given as follows:
$$
\begin{array}{lclcc}
\hline\\
\qquad \alpha&X_\alpha\qquad &\exp tX_\alpha&\\
\hline\\
\epsilon_i-\epsilon_j,\;1\le j<i\le n&E_{i,j}-E_{j+n,i+n}&1+t(E_{i,j}-E_{j+n,i+n})&{\rm positive}&\\    %{\rm ht}_\Delta(\alpha)=i-j
\epsilon_i-\epsilon_j,\;1\le i<j\le n&E_{i,j}-E_{j+n,i+n}&1+t(E_{i,j}-E_{j+n,i+n})&&\Phi_1\\    %i-j
\epsilon_i+\epsilon_j,\;1\le i<j\le n&E_{i,n+j}-E_{j,n+i}&1+t(E_{i,n+j}+E_{j,n+i})&{\rm positive}&\\    %i-j+1
-\epsilon_i-\epsilon_j,\;1\le i<j\le n&E_{n+i,j}-E_{n+j,i}&1+t(E_{n+i,j}+E_{n+j,i})&&\Phi_2\\  %-i+j-1
2\epsilon_i,\;1\le i\le n&E_{i,n+i}&1+tE_{i,n+i}&{\rm positive}&\\    %2i-1
-2\epsilon_i,\;1\le i\le n&E_{n+i,i}&1+tE_{n+i,i}&&\Phi_2\\   %-2i+1
\hline\\
\end{array}
$$
%die Aussage über die Hohen folgt aus: $e_i-e_j=\sum e_i-e_{i-1}+\cdots+e_{j+1}-e_j$, $e_i+e_j=\sum e_i-e_{i-1}+\cdots+e_{j+1}-e_j+2e_j$, $2e_i=2e_1+2(e_2-e_1)+\cdots+2(e_i-e_{i-1})$
Here, $\epsilon_i:\,{\bf T}^0\rightarrow{\Bbb C}$ is defined by mapping ${\rm diag}(\alpha_1,\ldots,\alpha_n,\alpha_1^{-1},\ldots,\alpha_n^{-1})$ to $\alpha_i$ ($1\le i\le n$) and $\exp:\,{\mathfrak g}^0\rightarrow {\bf G}({\Bbb C})$ is the exponential.  
%das lasst sich unter die theorie der Chev gruppen unterordnen. Sei $\rho$ die kanonischeDarstellung von ${\mathfrak sp}_{2n}$ auf ${\Bbb C}^{2n}$. dann gilt
%${\cal D}(\rho(g))=g$, d.h. $\exp (t \rho(x)=\exp tx$ fur $x \in {\mathfrak g}_\alpha). D.h. die exp funktion ist gerade die matrix exp funktion, die wir beutzt haben
%
%
%(to verify the column "$\exp tX_\alpha$" use $E_{ij} E_{ab}=\delta_{ja} E_{ib}$ plus the power series defining $\exp$). 
We choose as a basis of the root system 
$$
\Delta=\{2\epsilon_1,\epsilon_{i+1}-\epsilon_{i},\,i=1,\ldots,n-1\}.
$$
%$\epsilon_{i}-\epsilon_{j}=\epsilon_{i}-\epsilon_{i+1}+\cdots+\epsilon_{j-1}-\epsilon_{j}$
%2\epsilon_i=2\epsilon_1+2(\epsilon_{1}-\epsilon_{i})
%\epsilon_{i}+\epsilon_{i+1}=2\epsilon_i-\epsilon_{i}+\epsilon_{i+1}$

\medskip

We extend the roots $\alpha\in \Phi$ to ${\bf T}$ by setting them equal to $1$ on the center ${\bf Z}\le{\bf T}$. 
%wurzeln einer reduktiven gruppe sind immer $=1$ auf dem zentrum

\bigskip

{\bf (4.2) } We denote by $X({\bf T})$ the (additively written) group of (${\Bbb Q}$-)characters of ${\bf T}$. Since ${\bf T}$ is connected, 
$X({\bf T})$ is a finitely generated, free abelian group. Let $\gamma_1,\ldots,\gamma_{n+1}$ be a ${\Bbb Z}$-basis of $X({\bf T})$. For any 
$\lambda=\sum_{i=1}^{n+1} \lambda_i\gamma_i\in X({\bf T})$ ($\lambda_i\in{\Bbb Z}$) we set
$$
v_p(\lambda)={\rm min}_i\,v_p(\lambda_i)={\rm max}\,\{m\in{\Bbb N}_0: \lambda\in p^m X({\bf T})\}.
$$
Hence, $\lambda\equiv \lambda'\pmod{p^m X({\bf T})}$ is equivalent to $v_p(\lambda-\lambda')\ge m$ (compare section (1.8)).

%We choose a basis 
%$\{\omega_\alpha\}_{\alpha\in\Delta}$ of $X({\bf T}^0)$ and a character $\kappa\in X({\bf T})$ such that 
%$$
%X({\bf T})=\bigoplus_\alpha{\Bbb Z}\omega_\alpha\oplus{\Bbb Z}\kappa.
%$$

\bigskip

{\bf (4.3) Irreducible Representations. } Let $\lambda\in X({\bf T})$ be a dominant character, i.e. the restriction $\lambda^\circ=\lambda|_{{\bf T}^0}$ is a dominant 
character of ${\bf T}^0$ and $\lambda|_{\bf Z}$ is an (algebraic) character. We denote by $(\pi_\lambda,L_\lambda)$ the irreducible representation of 
${\bf G}(\bar{\Bbb Q}_p)$ of highest weight $\lambda$ on a $\bar{\Bbb Q}_p$-vector space $L_\lambda$. The representation $L_\lambda$ is defined over ${\Bbb Z}$, 
i.e. there is a ${\Bbb Z}$-submodule $L_\lambda({\Bbb Z})$ in $L_\lambda$ such that $L_\lambda=L_\lambda({\Bbb Z})\otimes\bar{\Bbb Q}_p$ 
and which is stable under ${\bf G}({\Bbb Z})$. Thus, for any ${\Bbb Z}$-algebra $R$ we obtain a ${\bf G}(R)$-module 
$L_\lambda(R)=L_\lambda({\Bbb Z})\otimes R$.

\bigskip

%---------------------------------------------
%
%??? brauchen wir spater auch die reelle Lie alg von ${\bf G}$ ??? die kopm,lexe gruppe ${\bf G}({\Bbb C})$ tuacht als koeffzsystem auf, ${\bf G}({\Bbb R})$ als symmetrischer Raum, der $\Gamma$ enthalt.
%
%---------------------------------------------

\section{The Hecke algebra ${\cal H}$ (attached to ${\bf GSp}_{2n}$)}

\bigskip

We define the Hecke algebra ${\cal H}$ which we shall be using. Essentially ${\cal H}$ omits all Hecke operators at primes 
dividing the level and its local component at the prime $p$ is generated by one single Hecke operator $T_p$. 

\bigskip

{\bf (5.1) The local level subgroup ${\cal I}$. } Let $g=(g_{ij})_{ij}\in{\bf G}(K)$ be arbitrary. We partition $g$ as 
$$
g=\Mat{\sf A}{\sf B}{\sf C}{\sf D}=\left(
\begin{array}{ccc|ccc}
\delta_{1}&&A_-&b_1&&B\\
&\ddots&&&\ddots&\\
A_+&&\delta_{n}&B'&&b_n\\
\hline\\
c_1&&C&\delta_{n+1}&&A_+'\\
&\ddots&&&\ddots&\\
C'&&c_n&A_-'&&\delta_{2n}\\
\end{array}
\right)\leqno(1)
$$
where $A_+=(a_{ij}^+)$, $A_-=(a_{ij}^-)$, $B=(b_{ij})$, $B'=(b'_{ij})$ and so on (a prime indicates that the submatrix and its corresponding 
primed submatrix are connected via one of the relations in equation (1) in section 4; the positive root spaces "correspond to entries in $A_+,B,B'$ and 
$b_1,\ldots,b_n$"). We denote by ${\bf B}\le {\bf G}$ the minimal parabolic subgroup attached to the basis $-\Delta$ of $\Phi$ and we define 
$$
{\cal I}\le {\bf G}({\Bbb Z}_p)
$$ 
as the subgroup consisting of all matrices $g$ such that $g\pmod{p{\Bbb Z}_p}$ is contained in ${\bf B}({\Bbb Z}_p/(p))$. Thus, $g\in{\bf G}({\Bbb Z}_p)$ as in equation (1) is contained in ${\cal I}$ precisely if
$$
a^+_{ij}\in p{\Bbb Z}_p,\,b_{ij}\in p{\Bbb Z}_p,\,b_1,\ldots,b_n\in p{\Bbb Z}_p.\leqno(2)
$$

%VVVV
%
%Thus, ${\cal I}$ is generated by $\exp t_\alpha x_\alpha$ where $\alpha\in \Phi$ and $t_\alpha\in p{\Bbb Z}_p$ if $\alpha$ is positive and $t_\alpha\in{\Bbb Z}_p$ if $\alpha$ is negative. %Explicitly ${\cal I}$ consists of all elements $g=(g_{ij})$ in ${\bf G}({\Bbb Z}_p)^0$ satisfying $g_{ij}\in p{\Bbb Z}_p$ if $1\le j<i\le n$ and $g_{ij+n}\in p{\Bbb Z}_p$ if $1\le i\le j\le n$.
%
%AAAA

{\bf (5.2) An auxiliary Lemma. } For any prime $\ell$ we set
$$
{\bf T}({\Bbb Q}_\ell)^+=\{h\in {\bf T}({\Bbb Q}_\ell):\,v_\ell(\alpha(h))\ge 0\;\mbox{for all}\;\alpha\in\Delta\}.
$$

\bigskip

{\bf Lemma. } {\it 1.) Let $s\in {\bf T}({\Bbb Q}_p)^+$ and set
$$
{\cal V}_s=\{1+t_\alpha X_\alpha,\,\alpha\in\Phi^-,\,t_\alpha\in{\Bbb Z}_p/p^{-v_p(\alpha(s))}{\Bbb Z}_p\}\quad (\subseteq{\bf B}({\Bbb Z}_p)).
$$
%(note that for negative roots $\alpha\in\Phi^-$ we have $-v_p(\alpha(h))\ge 0$). 
Then, ${\cal I}s{\cal I}=\bigcup_{v\in{\cal V}_h}{\cal I}sv$. 

\medskip

2.) Let $s_1,s_2\in {\bf T}({\Bbb Q}_p)^+$. Then 
$$
{\cal V}_{s_1s_2}={\rm Ad}(s_2^{-1})({\cal V}_{s_1}){\cal V}_{s_2}.
$$
}

\bigskip

{\it Proof. } 1.) Since ${\cal I}s{\cal I}=\bigcup_{v} {\cal I}sv$ with $v$ running over $s^{-1}{\cal I}s\cap{\cal I}\backslash {\cal I}$, we have to show that ${\cal V}_s$ is a system of representatives for $s^{-1}{\cal I}s\cap{\cal I}\backslash {\cal I}$. 
To this end, we set $\alpha_{ij}=\epsilon_i-\epsilon_j$ and $\beta_{ij}=-\epsilon_i-\epsilon_j$, $1\le i<j\le n$ and $\Phi_1=\{\alpha_{ij},\,1\le i<j\le n\}$ and $\Phi_2=\{\beta_{ij},\,1\le i\le j\le n\}$;  
hence, $\Phi^-=\Phi_1\cup \Phi_2$. Since $s^{-1} xs={\rm Ad}(s^{-1})(x)$ and since ${\rm Ad}(s)$ acts on ${\mathfrak g}_\alpha$ via $\alpha$ we obtain that $s^{-1}{\cal I}s\cap{\cal I}$ consists of all matrices $g=(g_{ij})\in{\cal I}$ satisfying
$$
g_{ij},g_{n+j,n+i}\in p^{-v_p(\alpha_{ij}(s))},\, 1\le i<j\le n\quad\mbox{and}\quad g_{i+n,j},g_{j+n,i}\in p^{-v_p(\beta_{ij}(s))},\,1\le i\le j\le n.\leqno(3)
$$
%(Bemerkung: Es genugt dass $s^{-1}{\cal I}s\cap{\cal I}$ in der menge aller matrizen in ${\cal I}$ enthalten ist die obige gleichung erfullen; wir bruachen nicht gleichheit der beiden mengen.)
Now let $g=(g_{ij})\in{\cal I}$ be arbitrary. We will perform two sets of elementary row operations to transform $g$ into the unit matrix. 

a.) Multiplying $g$ from the left by matrices of the form $1+t(E_{ij}-E_{j+n,i+n})=(1+tE_{ij})(1-tE_{n+j,n+i})$ with $1\le i<j\le n$ and 
$t\in {\Bbb Z}_p/p^{-v_p(\alpha_{ij}(s))}{\Bbb Z}_p$, i.e. by performing (two) elementary row operations we can achieve that all entries $g_{ij}$ with 
$1\le i<j \le n$ are contained in $p^{-v_p(\alpha_{ij}(s))}{\Bbb Z}_p$. (First eliminate the entries in the most right column in $A_-$, then eliminate the 
entries in column $n-1$ and so on; 
note that the $\delta_i=g_{ii}$ are units in ${\Bbb Z}_p$.) Analogously, multiplying from the left by matrices 
$1+t(E_{n+i,j}+E_{n+j,i})=(1+tE_{n+i,j})(1+tE_{n+j,i})$ with $i\le j$ and $t\in{\Bbb Z}_p/p^{-v_p(\beta_{ij}(s))}{\Bbb Z}_p$ we can achieve that 
$g_{n+i,j}\in p^{-v_p(\beta_{ij}(s))}{\Bbb Z}_p$ for all $1\le i\le j\le n$. We note that the matrices by which we multiplied are contained in 
${\cal V}_s$ because they are contained in $1$-parameter subgroups corresponding to negative roots (namely $\alpha_{ij},\beta_{ij}$; cf. the table in section (4.1)). 

\medskip

b.) Now, multiplying further from the left 

- by matrices $1+t(E_{ij}-E_{j+n,i+n})$ with $i<j$ and $t\in p^{-v_p(\alpha_{ij}(s))}$ we can achieve $A_-=0$

- by matrices $1+t(E_{i+n,j}+E_{j+n,i})$ with $i\le j$, $t\in p^{-v_p(\beta_{ij}(s))}$ we can achieve $C=0$ and $c_{11}=0,\ldots,c_{nn}=0$

- by matrices $1+t(E_{ij}-E_{j+n,i+n})$ with $i>j$, $p\in p{\Bbb Z}_p$, we can achieve $A_+=0$

- by matrices $1+t(E_{i,j+n}+E_{j,i+n})$ with $j\ge i$, $p\in p{\Bbb Z}_p$, we can achieve $B=0$ and $b_{11}=0,\ldots,b_{nn}=0$

- by matrices $(\delta_1^{-1},\ldots,\delta_n^{-1},\delta_1,\ldots,\delta_n)$, $\delta_i\in{\Bbb Z}_p^*$, we can achieve $\delta_1=\cdots=\delta_n=1$.

We note that the matrices by which we multiplied are contained in ${\cal I}\cap s^{-1}{\cal I}s$ by equations (2) and (3). Since ${\sf A}=1$, the first relation in equation (1) in  section 4 implies ${\sf B}={\sf B}^t$, hence, ${\sf B}=0$. The last relation in equation (1) in section 4 then implies
${\sf D}=1$ and the second relation in (1) in  section 4 finally implies ${\sf C}={\sf C}^t$, hence, ${\sf C}=0$. Thus, we have found that
$$
\prod_j k_j \prod_i v_i \, g=1
$$
with certain $k_j\in{\cal I}\cap s^{-1}{\cal I}s$ and $v_i\in{\cal V}_s$. Thus, $(\prod v_i) g$ is contained in $s^{-1}{\cal I}s\cap{\cal I}\backslash {\cal I}$ which shows that ${\cal V}_s$ is a system of representatives for $s^{-1}{\cal I}s\cap{\cal I}\backslash {\cal I}$. 
%bea.: nach schritt I.) konnen wir aus eq. (3) noch  nicht schliessen dass die transformierte matrix in $s^{-1}{\cal I}s\cap{\cal I}$ liegt, denn wir haben die eintrage der gestrichenen untermatrizen noch nicht verifiziert; deshalb schritt II wo dies dann am ende mittels eq. (1) in sec (6.1) geschieht.
Taking into account that the entries of a matrix $g=(g_{ij})\in s^{-1}{\cal I}s\cap{\cal I}$ satisfy equation (3) it is not difficult to verify that the elements in ${\cal V}_s$ are different modulo $s^{-1}{\cal I}s\cap{\cal I}$. Thus, the proof of the first part is complete.

2.) Since ${\rm Ad}(s_2)X_\alpha=\alpha(s_2)X_\alpha$ we find that
$$
{\rm Ad}(s_2^{-1})({\cal V}_{s_1})=\{ 1+t_\alpha X_\alpha,\,\alpha\in \Phi^-,\,t_\alpha\in p^{v_p(\alpha(s_2^{-1}))}{\Bbb Z}_p/p^{-v_p(\alpha(s_1))+v_p(\alpha(s_2^{-1}))}{\Bbb Z}_p\}.
$$
Since $v_p(\alpha(s_2^{-1}))=-v_p(\alpha(s_2))$ and $-v_p(\alpha(s_1))+v_p(\alpha(s_2^{-1}))=-v_p(\alpha(s_1s_2))$ an easy calculation yields 
the claim. This finishes the proof of the Lemma.

\bigskip

{\bf (5.3) The local Hecke algebra at $p$. } We set
$$
D_p={\cal I}{\bf T}({\Bbb Q}_p)^+{\cal I}\,\le{\bf G}({\Bbb Q}_p).
$$
In the Proposition below we will see that $D_p$ is a semigroup, hence, we can define the local Hecke algebra 
$$
{\cal H}({\cal I}\backslash D_p/ {\cal I})
$$ 
attached to the pair $(D_p,{\cal I})$. For any $s\in {\bf T}({\Bbb Q}_p)^+$ we define the element $T_s={\cal I}s{\cal I}$, i.e. 
${\cal H}({\cal I}\backslash D_p/ {\cal I})$ is the ${\Bbb Z}$-linear span of the elements $T_s$, $s\in {\bf T}({\Bbb Q}_p)^+$.

\bigskip

{\bf Proposition. }{\it 1.) $D_p\le {\bf G}({\Bbb Q}_p)$ is a semi group.

\medskip

2.) For all $s_1,s_2\in{\bf T}({\Bbb Q}_p)^+$ we have $T_{s_1}T_{s_2}=T_{s_1s_2}$. In particular, the Hecke algebra ${\cal I}\backslash D_p/{\cal I}$ is commutative. 

}

\medskip

{\it Proof. } 1.) We compute
$$
{\cal I}s_1{\cal I}s_2{\cal I}=\bigcup_{v\in{\cal V}_{s_1}} {\cal I}s_1 v s_2{\cal I}=\bigcup_{v\in{\cal V}_{s_1}} {\cal I}s_1 s_2 {\rm Ad}(s_2^{-1})(v) {\cal I}
$$
Since $s_2\in  {\bf T}({\Bbb Q}_p)^+$ we know that $v_p(\beta(s_2^{-1}))\le 0$ for all positive roots $\beta$. Hence, we obtain for any $v=1+t_\alpha X_\alpha\in{\cal V}_{s_1}$:
$$
{\rm Ad}(s_2^{-1})(v)=1+t_\alpha \alpha(s_2^{-1}) X_\alpha\in {\bf B}({\Bbb Z}_p)\subset{\cal I}
$$
(note that by definition of ${\cal V}_{s_1}$ the root $\alpha$ is negative). Thus, ${\rm Ad}(s_2^{-1})(v)\in{\cal I}$ and we obtain ${\cal I}s_1{\cal I}s_2{\cal I}={\cal I}s_1s_2{\cal I}$. Thus $D_p$ is closed under multiplication.

2.) To prove the second claim we compute using part 2.) of the previous Lemma
$$
T_{s_1}T_{s_2}=\bigcup_{v\in {\cal V}_{s_1},w\in{\cal V}_{_2}}{\cal I}s_1vs_2w=
\bigcup_{v\in {\cal V}_{s_1},w\in{\cal V}_{s_2}} {\cal I}s_1s_2{\rm Ad}(s_2^{-1})(v)w=\bigcup_{z\in {\cal V}_{s_1s_2}}{\cal I}s_1s_2 z=T_{s_1s_2}.
$$
Thus, the proof of the Proposition is complete.

\bigskip

{\bf (5.4) The adelic Hecke Algebra. } We fix an integer $N$ which is not divisible by $p$. We select a compact open subgroup 
$U=\prod_{\ell\not=\infty} U_\ell$ of ${\bf G}(\hat{\Bbb Z})$ and a sub semigroup $D=\prod_{\ell\not=\infty} D_\ell$ of ${\bf G}({\Bbb A}_f)$ as 
follows. 
For all primes $\ell$ not dividing $pN$ we set $U_\ell={\bf G}({\Bbb Z}_\ell)$ and  $D_\ell={\bf G}({\Bbb Q}_\ell)$; at the prime $p$ we define 
$U_p={\cal I}$ and $D_p={\cal I}{\bf T}({\Bbb Q}_p)^+{\cal I}$ as in the previous section; 
for primes $\ell|N$ we only assume that $U_\ell\le{\bf G}({\Bbb Z}_\ell)$ is compact open and ${\rm det}(U_\ell)={\Bbb Z}_\ell^*$ and we set 
$D_\ell=U_\ell$. We denote by 
$$
{\cal H}(U_\ell\backslash D_\ell/U_\ell)\quad\mbox{resp.}\quad {\cal H}(U\backslash D/U)
$$ 
the local resp. global (adelic) Hecke algebra attached to the pair $(U_\ell,D_\ell)$ resp. $(U,D)$. We exhibit a set of generators for 
${\cal H}(U\backslash D/U)$ as follows. We define
$$
\Sigma_\ell^+=\{{\rm diag}(\ell^{e_1},\ldots,\ell^{e_n},\ell^{-e_1},\ldots,\ell^{-e_n}),\,e_i\in{\Bbb Z},\,0\le e_1\le e_2\le\cdots\le e_n\}
$$
if $\ell\not|N$ and $\Sigma^+_\ell=\{1\}$ if $\ell|N$. Thus, $v_\ell(\alpha(h))\ge 0$ for all $h\in\Sigma_\ell^+$ and all $\alpha\in\Delta$. Any element 
$T=U_\ell s U_\ell$ in ${\cal H}(U_\ell\backslash D_\ell/U_\ell)$ has a representative $s\in\Sigma_\ell^+$ (for primes $\ell$ not dividing $Np$ this 
is well known and for primes $\ell$ dividing $Np$ this is immediate from the definition of $U_\ell$ and $D_\ell$.) In particular, the ${\Bbb Z}$-algebra 
${\cal H}(U\backslash D/U)$ is generated by the following elements
$$
{\cal H}(U\backslash D/U)=\langle UsU,\,s\in\bigcup_\ell \Sigma_\ell^+  \rangle. \leqno(4)
$$
Moreover, ${\cal H}(U\backslash D/U)$ is commutative, because this holds locally at all primes: for primes $\ell\not|Np$ again this is well known, 
for primes $\ell|N$ this is trivial and for $p$ this follows from (5.3) Proposition.

\bigskip

{\bf (5.5) The global (non-adelic) Hecke algebra. } We set 
$$
\Gamma=U\cap {\bf G}({\Bbb Q})\quad\mbox{and}\quad \Delta=D\cap {\bf G}({\Bbb Q}). 
$$
Thus, $\Gamma$ satisfies the following local condition at the prime $p$:
$$
\Gamma\le{\cal I}\,(\le {\bf G}({\Bbb Q}_p)).\leqno(5)
$$
We note that $U\subset D$, hence, $\Gamma\subset \Delta$. In the next Lemma we compare the Hecke algebras attached to the pairs $(D,U)$ and 
$(\Delta,\Gamma)$
$$
\begin{array}{ccc}
\Gamma&\subset&\Delta\\
\cap&&\cap\\
U&\subset&D.\\
\end{array}
$$

\bigskip

{\bf (5.5.1) Lemma. }{\it The canonical map $\Gamma\alpha\Gamma\mapsto U\alpha U$ induces an isomorphism of rings
$$
{\cal H}(\Gamma\backslash \Delta/\Gamma)\rightarrow {\cal H}(U\backslash D/U).
$$
}

\medskip

{\it Proof. }  According to [M], Theorem 2.7.6, p. 72, we have to show that

i) $D=\Delta U$

ii) $U\alpha U=U\alpha\Gamma$ for all $\alpha\in \Delta$

iii) $U\alpha\cap \Delta=\Gamma\alpha$ for all $\alpha\in \Delta$.

i) is an immediate consequence of strong approximation which holds since ${\rm det}(U)=\hat{\Bbb Z}^*$. 
%"$\supset$": klar, da $\Delta,K\subset D$. "$\subset$": let $d\in D$. Strong approx (beachte dass $K$ offen ist; das genugt fur strong approx) $\Rightarrow$ $d=\gamma k$, $\gamma\in{\bf G}({\Bbb Q})$, $k\in K$. Da \gamma=dk^{-1}\in DK=D$ 
%folgt $\gamma\in D\cap{\bf G}({\Bbb Q})=\Gamma$.
%Beachte: starke Approx gilt wenn die klassenzahl $=1$ ist, d.h. genau dann wenn $|{\Bbb Q}^*\backslash {\Bbb A}^*_f/K_f|=1$; das ist im fall $F={\Bbb Q}$ genau dann der fall wenn ${\rm det}\,K_f=\hat{\Bbb Z}*$
We prove ii). The inclusion "$\supseteq$" is obvious. To prove the reverse inclusion we note that $U\alpha U=\bigcup_v U\alpha v$, where $v$ runs 
over a system of representatives of $\alpha^{-1}U\alpha\cap U\backslash U$. Thus, we have to show that $\Gamma$ contains a system of 
representatives of $\alpha^{-1}U\alpha\cap U\backslash U$. Let $u\in U$ be arbitrary. 
Since $\alpha^{-1}U\alpha\cap U\le {\bf G}({\Bbb A}_f)$ is a compact open subgroup, strong approximation yields $u=\gamma v$ with 
$\gamma\in{\bf G}({\Bbb Q})$ and $v\in \alpha^{-1}U\alpha\cap U$. Hence, $\gamma=uv^{-1}$ is contained in 
$U\cap {\bf G}({\Bbb Q})=\Gamma$. 
Thus, $\gamma$ is a representative of the coset of $u$ in $\alpha^{-1}U\alpha\cap U\backslash U$. It remains to verify that strong approximation 
holds with respect to $\alpha^{-1}U\alpha\cap U$, i.e. ${\rm det}\,(\alpha^{-1}U\alpha\cap U)=\hat{\Bbb Z}^*$. It is sufficient to prove this locally 
for all primes $\ell$. 
If $\ell|N$ we know by definition of $U_\ell,D_\ell$ that $\alpha\in D_\ell=U_\ell$, hence, ${\rm det}(U_\ell\cap \alpha^{-1}U_\ell\alpha)={\Bbb Z}_\ell^*$ 
because the determinant is surjective on $U_\ell$. If $\ell\not|N$, then the Cartan decomposition in case $\ell\not|Np$ and the definition of $D_p$ show 
that $\alpha\in D_\ell$ can be written $\alpha=u_1tu_2$ where $u_1,u_2\in U_\ell$ and $t\in D_\ell$ is a diagonal matrix. 
We obtain ${\rm det}(U_\ell\cap \alpha^{-1}U_\ell\alpha)={\rm det}(U_\ell\cap t^{-1}U_\ell t)$. Let $\lambda\in{\Bbb Z}_\ell^*$ arbitrary; since 
$s={\rm diag}(\lambda,1,\ldots,1)\in U_\ell$ commutes with $t$ we see that it is contained in $U_\ell\cap t^{-1}U_\ell t$. 
Since ${\rm det}(s)=\lambda$ we obtain ${\rm det}(U_\ell\cap \alpha^{-1} U_\ell\alpha)={\Bbb Z}_\ell^*$. Thus, strong approximation holds 
and ii) is proven. Finally, iii) is immediate since $\Gamma$ is contained in $U$ and in $\Delta$ and since $U\cap {\bf G}({\Bbb Q})=\Gamma$. 
%"$\subset". let $g\in U\alpha\cap \Delta$ $\Rightarrow$ $g=u\alpha=d$ with $u\in U$, $d\in D$. Since $\alpha\in \Delta$ we know that $\alpha^{-1}\in{\bf G}({\Bbb Q})^+$, hence, 
%u=d\alpha^{-1}\in{\bf G}({\Bbb Q})^+$. Thus, $u\in U\cap {\bf G}({\Bbb Q})^+=\Gamma$.
Thus, the proof of the lemma is complete.
% was wir brauchen damit das argument im beweis funktioniert ist, dass 1.) jedes element in $D_\ell$, also insbesondere $\alpha$, die form $\alpha=u_1tu_2$ mit $u_1,u_2\in U_\ell$ und $t$ diagonal hat un dass 2.) ${\rm diag}(\lambda,\ldots,\lambda,1,\ldots,1)\in U_\ell$ ist. Dann folgt ${\rm det}(U\ell\cap \alpha^{-1}U_\ell\alpha)={\Bbb Z}_\ell^*$. 

\bigskip

{\bf (5.5.2)} Since the (adelic) Hecke algebra attached to $(D,U)$ is commutative (cf. section (5.4)), we obtain from the above Lemma that 
${\cal H}(\Gamma\backslash \Delta/\Gamma)$ is a commutative algebra. For $s\in \Delta$  we set
$$
T_s=\Gamma s\Gamma;
$$
equation (4) in (5.4) and (5.5.1) Lemma imply that ${\cal H}(\Gamma\backslash \Delta/\Gamma)$ is generated by the following elements
$$
{\cal H}(\Gamma\backslash \Delta/\Gamma)=\langle  \Gamma s\Gamma,\;s\in\bigcup_\ell \Sigma_\ell^+    \rangle
$$
%moreover, we set
%$$
%\Sigma_\ell^<=\left\{\begin{array}{cc}
%\Delta\cap {\bf T}({\Bbb Q}_\ell)^<&l\not|N\\
%\{1\}&\ell|N.
%\end{array}
%\right.
%$$
%Here, ${\bf T}({\Bbb Q}_\ell)^<$ is viewed as embedded in ${\bf T}({\Bbb A}_f)$. We note that an easy calculation shows that for all $\ell\not|N$
%$$
%\Sigma_\ell^<={\bf Z}({\Bbb Q})\{{\rm diag}(\ell^{e_1},\ldots,\ell^{e_n},\ell^{-e_1},\ldots,\ell^{-e_n}),\,e_1\le e_2\le\cdots\le e_n\le 0\}.
%$$
%Since ${\cal H}(U\backslash D/U)$ is the product of the local Hecke algebras ${\cal H}(U_\ell\backslash D_\ell/U_\ell)$ and ${\cal H}(U_\ell\backslash D_\ell/U_\ell)$ is generated by the Hecke operators %$U_\ell s U_\ell$ with $s\in\bigcup_{\ell\not|Np}\Sigma_\ell$, the above Lemma implies that ${\cal H}(\Gamma\backslash \Delta/\Gamma)$ is generated as ${\Bbb Z}$-algebra by the elements  $T_s$, $s\in %\bigcup_\ell \Sigma_\ell^<$. 
(note that $\Sigma_\ell^+\subseteq \Delta$ for all primes $\ell$). We define the element
$$
h_p={\rm diag}(p^1,p^2,\ldots,p^n,p^0,p^{-1},\ldots,p^{-n+1})\in\Sigma_p^+,
$$
hence, $\alpha(h_p)=p$ for all $\alpha\in \Delta$ 
%note that $h_p$ has determinant $p^n$, hence, the corresponding element in ${\bf Sp}_n$ is $h_p'={\rm diag}(p^{1/2},p^{3/2},\ldots,p^{n-1/2},p^{-1/2},p^{-3/2},\ldots,p^{-n+1/2})$ und jetzt siehnt man, dass $2\epsilon_1(h_p')=1$.
and we denote the corresponding Hecke operator by
$$
T_p=T_{h_p}=\Gamma  h_p\Gamma.
$$
Since $T_p$ maps to ${\cal I}h_p{\cal I}\in {\cal H}({\cal I}\backslash D_p/{\cal I})\le {\cal H}(U\backslash D/U)$ under the isomorphism in (5.5.1) Lemma, 
(5.3) Proposition  part 2.) implies that
$$
T_p^e=T_{h_p}^e=T_{h_p^e}.\leqno(6)
$$

\bigskip

{\bf (5.6) The Hecke algebra ${\cal H}$. } We define the ${\Bbb Z}$-algebra ${\cal H}_{\Bbb Z}$ as the ${\Bbb Z}$-subalgebra of ${\cal H}(\Gamma\backslash \Delta/\Gamma)$ 
which is generated by the Hecke operators $T_s$ with $s\in \bigcup_{\ell\not|Np} \Sigma_\ell^+\cup\{h_p\}$; hence,
$$
{\cal H}_{\Bbb Z}\cong\bigotimes_{\ell\not|Np} {\cal H}(U_\ell\backslash D_\ell/U_\ell)\otimes{\Bbb Z}[T_p].
$$ 
%
%Beachte dass ${\cal H}(U_p\backslash D_p/U_p)=\bigoplus_{e\ge 0} \bar{\Bbb Q}_p T_{h_p^e}=\bigoplus_{e\ge 0} \bar{\Bbb Q}_p T_p^e$ 
%(cf. (5.5) equation (6)), woraus folgt, dass ${\cal H}_{\Bbb Z}\cong\bigotimes_{\ell\not|Np} {\cal H}(U_\ell\backslash D_\ell/U_\ell)\otimes{\Bbb Z}[T_p]$
%
We set 
$$
\Sigma^+:=\prod_{\ell\not|Np}\Sigma_\ell^+\cdot\{h_p^m,\,m\in{\Bbb N}_0\}.
$$
As a ${\Bbb Z}$-module, ${\cal H}_{\Bbb Z}$ then is generated by the operators $T_s$ with $s\in\Sigma^+$. 
Finally, we put ${\cal H}={\cal H}\otimes_{\Bbb Z} \bar{\Bbb Q}_p$ and ${\cal H}_{\cal O}={\cal H}_{\Bbb Z}\otimes{\cal O}$.

\bigskip

%${\cal H}(\Gamma\backslash {\bf G}({\Bbb Q}_p))\Gamma)={\cal H}(\Gamma\backslash {\bf G}({\Bbb Q}_p)\cap M_{2n}({\Bbb Z}_p))\Gamma)[T]$
%where $T=\Gamma p1_{2n}\Gamma$ is the hecke operator attached to center. (Any $g\in{\bf G}(){\Bbb Q})$ can mulitplied into $M_{2n}({\Bbb Z}_p)$ by an element in the center $Z=\{\lambda 1_{2n}\}$ of ${\bf G}$.

\section{The cohomology groups ${\bf H}_\lambda$ (attached to the Siegel upper half plane)}

\bigskip

We recall the normalization of the Hecke operators which leads to an action of the Hecke algebra on the $p$-adically integral 
cohomology. We also state the simple topological trace formula of Bewersdorff. 

\bigskip

{\bf (6.1) Normalization of Hecke algebra representations. } We keep the notations from section 5. In particular, 
$\Gamma=U\cap{\bf G}({\Bbb Q})\le {\bf G}({\Bbb Z})$ with $U\le{\bf G}(\hat{\Bbb Z})$ defined as in (5.4) and 
${\cal H}_{\Bbb Z}\cong\bigotimes_{\ell\not|Np} {\cal H}(U_\ell\backslash D_\ell/U_\ell)\otimes{\Bbb Z}[T_p]$ is the Hecke 
algebra defined in (5.6). To ensure that the Hecke algebra later will act on $p$-adically integral cohomology, we have to normalize the 
action of the Hecke algebra. This depends on a choice of a dominant character $\lambda\in X({\bf T})$. We first define a ${\Bbb Z}$-
algebra morphism
$$
\begin{array}{cccc}
\varphi_\lambda=\otimes_{\ell\not|N}\varphi_\ell:&{\cal H}_{\Bbb Z}&\rightarrow&{\cal H}_{\Bbb Z}\\
&T&\mapsto&\{T\}_\lambda\\
\end{array}
$$ 
as follows. For all primes $\ell\not|Np$ we denote by $\varphi_{\lambda,\ell}:\,{\cal H}(U_\ell\backslash D_\ell/U_\ell)\rightarrow {\cal H}(U_\ell\backslash D_\ell/U_\ell)$ the identity map. At the prime $p$ we note that ${\cal H}(U_p\backslash D_p/U_p)={\Bbb Z}[T_p]$ is a polynomial algebra 
generated by $T_p$. We define 
$\varphi_{\lambda,p}:\,{\cal H}(U_p\backslash D_p/U_p)\rightarrow {\cal H}(U_p\backslash D_p/U_p)$ by sending $T_p=T_{h_p}$ to 
$\{T_p\}_\lambda:=\lambda(h_p)T_p$. Equation (6) in (5.5) implies that $\varphi_\lambda(T_{h_p^e})=\varphi_\lambda(T_p^e)=\lambda(h_p^e)T_p^e$. Moreover, 
since $\lambda$ is dominant and $h_p\in\Sigma_p^+$ we obtain $\lambda(h_p)\in{\Bbb Z}$, hence, 
$\varphi_\lambda(T_p)=\lambda(h_p) T_p\in{\cal H}_{\Bbb Z}$ and $\varphi_\lambda$ is defined over ${\Bbb Z}$. 
%
%denn: $\lambda$ dominant $\Rightarrow \lambda=\sum_{\alpha\in\Delta} n_\alpha \omega_\alpha$ mit $n_\alpha\in{\Bbb N}_0$ (s. [Humphreys, Lie Alg], sec. 13.1, p. 67 unten).  Weiter  ist $\omega_\alpha=\sum_{\beta\in\Delta} m_{\alpha,\beta} \beta$ mit $m_{\alpha,\beta}\in{\Bbb Q}_{\ge 0}$ (s. [Humphreys, Lie Alg], sec. 13.2, Table 1). 
%wir erhalten $\lambda=\sum_{\alpha,\beta\in\Delta} n_\alpha m_{\alpha,\beta} \beta$ wobei $c_\beta=\sum_\alpha n_\alpha m_{\alpha,\beta}\ge 0$. Damit folgt
%$v_\ell(\lambda(h))=\sum_\beta c_\beta v_\ell(\beta(h))$. Da $v_\ell(\beta(h))\ge 0$ fur all $\beta\in\Delta$, $h\in\Sigma^+$ und Primzahlen $\ell$, folgt
%$v_\ell(\lambda(h))\ge 0$ fur alle $\ell$ $\Rightarrow \lambda(h)\in{\Bbb Z}$
%
%(Bemerkung:  1.) $\lambda(h)\in{\Bbb Q}$ folgt sofort, da $\Sigma^*\subseteq {\bf G}({\Bbb Q})$) $\Rightarrow \lambda(h)\in{\Bbb Z}$ 
%2.) im Beweis von (7.2) Lemma haben wir nur gebraucht, dass $\alpha(h)\in{\Bbb Z}$ fur wurzeln $\alpha$ d.h. fur Elemente im Wurzelgitter,nicht im Gewichtsgitter)
%
Tensoring with $\bar{\Bbb Q}_p$ we obtain a $\bar{\Bbb Q}_p$-algebra morphism
$$
\varphi_\lambda:\,{\cal H}\rightarrow{\cal H}
$$ 
which is defined over ${\cal O}$.

\medskip

Let ${\bf H}$ be a ${\cal H}$-module. We define the $\lambda$-normalization "$\cdot_\lambda$" of the 
action of ${\cal H}$ on ${\bf H}$ by composing the ${\cal H}$-module structure on ${\bf H}$ with the ${\bar{\Bbb Q}_p}$-algebra morphism 
$\varphi_\lambda$, i.e.
$$
T\cdot_\lambda v=\{T\}_\lambda v\qquad(T\in{\cal H},\; v\in H).
$$
Thus, 
$$
T_h\cdot_\lambda v=T_h v\quad\mbox{if $h\in\Sigma_\ell^+$, $\ell\not|Np$}\quad\mbox{and}\quad T_p^e\cdot_\lambda v=\lambda(h_p^e)T_p^e v.
$$

\bigskip

{\bf (6.2) The Cohomology groups. }  We select a maximal compact subgroup $K_\infty$ of the connected component of the identity of ${\bf G}({\Bbb R})$ and we denote by $X={\bf G}({\Bbb R})/K_\infty {\bf Z}({\Bbb R})$ the symmetric space. 
We denote by ${\cal L}_\lambda$ the sheaf on $\Gamma\backslash \bar{X}$ attached to the irreducible ${\bf G}(\bar{\Bbb Q}_p)$-module $L_\lambda$ of highest weight $\lambda\in X({\bf T})$ (cf. (4.3)). The cohomology groups
$$
{\bf H}_\lambda=H^d(\Gamma\backslash X,{\cal L}_\lambda)
$$
then are modules under the Hecke algebra ${\cal H}$. From now on, by $d=d_n$ we will always understand the {\it middle degree} of the 
locally symmetric space $\Gamma\backslash X$. We denote by 
$$
{\bf H}_{\lambda,{\cal O}}=H^d(\Gamma\backslash X,{\cal L}_\lambda)_{\rm int}
$$ 
the image of $H^d(\Gamma\backslash X,{\cal L}_\lambda({\cal O}))$ in $H^d(\Gamma\backslash X,{\cal L}_\lambda)$. 
$H^d(\Gamma\backslash X,{\cal L}_\lambda)_{\rm int}$ is a ${\cal O}$-lattice in 
$H^d(\Gamma\backslash X,{\cal L}_\lambda)$. On the other hand, the normalized Hecke operators 
$\lambda(h)T_h$, $h\in\Sigma^+$, act on cohomology with integral coefficients $H^d(\Gamma\backslash X,{\cal L}_\lambda({\Bbb Z}))$ 
(cf. e.g. the proof of (7.2) Lemma below; cf. also [Ma 1], (5.4) Lemma for more details). If $h\in\Sigma_\ell^+$ with $\ell\not|Np$ then 
$\lambda(h)\in{\cal O}^*$ 
%even $\lambda(h)\in{\Bbb Z}_{(\ell)}$
is a $p$-adic unit, hence, the not normalized Hecke operator $T(h)$ already acts on cohomology with $p$-adically integral coefficients 
$H^d(\Gamma\backslash X,{\cal L}_\lambda({\cal O}))$. Thus, we only have to normalize the Hecke operator at the prime $p$ 
i.e. we obtain that w.r.t. the $\lambda$-normalized action the Hecke algebra ${\cal H}_{\cal O}$ acts on cohomology with $p$-adically integral 
coefficients and hence, acts on $H^d(\Gamma\backslash X,{\cal L}_\lambda)_{\rm int}$. Thus, ${\bf H}_\lambda$ 
is a finite dimensional ${\cal H}$-module which is defined over ${\cal O}$ with respect to the lattice ${\bf H}_{\lambda,{\cal O}}$. 
% kohomologiegruppen von arithmetischen Gruppen sind immer endlich dimensional (folgt aus Borel Serre kompaktifizierung)

\bigskip

{\bf (6.3) Slope subspace of cohomology. } We define the slope $\le\beta$ subspace
$$
{\bf H}_\lambda^{\le\beta}=H^d(\Gamma\backslash X,{\cal L}_\lambda)^{\le\beta}
$$
with respect to the action of the {\it normalized} Hecke operator $\{T_p\}_\lambda=\lambda(h_p)T_p$. In [Ma 1], 6.5 Theorem it is proven using only elementary means 
from representation theory that there are natural numbers $M(\beta)$, $\beta\in{\Bbb Q}_{\ge 0}$, s. t.
$$
{\rm dim}\,{\bf H}_\lambda^{\le\beta}\le \frac{1}{2}M(\beta)\quad\mbox{for all dominant $\lambda\in X({\bf T})$ and all $\beta\in{\Bbb Q}_{\ge 0}$}.\leqno(1)
$$
Hence, the assumptions in (3.6) are satisfied.

\bigskip

\bigskip

%{\it Note. } In order to avoid confusion with the root "$\alpha$" we denote the slope by "$\beta$". 

{\bf (6.4) A simple topological Trace Formula. } We denote by $\sim_\Gamma$ the equivalence relation on ${\bf G}({\Bbb Q})$ defined by conjugation, i.e. 
$x\sim_\Gamma y$, $x,y\in{\bf G}({\Bbb Q})$, precisely if $x,y$ are conjugate by an element $\gamma\in\Gamma$ and we denote by $[\Xi]_\Gamma=[\Xi]$ the 
conjugacy class of $\Xi\in{\bf G}({\Bbb Q})$.

\bigskip

{\bf Theorem }(cf. [B]). {\it Let $h\in\Sigma^+$. There are integers $c_{[\Xi]}\in{\Bbb Z}$, $[\Xi]\in\Gamma h\Gamma/\sim_\Gamma$, such that the following holds. For all irreducible representations $L_\lambda$ we have
$$
{\rm tr}\,(T_h|H^d(\Gamma\backslash X,{\cal L}_\lambda))=\sum_{[\Xi]\in\Gamma h \Gamma/\sim_\Gamma} c_{[\Xi]} \,{\rm tr}\,(\Xi^{-1}|L_\lambda).\leqno(2)
$$

}

{\it Proof. } This is a direct consequence of 2.6 Satz in [B] taking into account that for regular weight $\lambda$ the cohomology 
of $\Gamma\backslash X$ vanishes in all degrees except for the middle degree $d$

%(the statement in [B] is even more general, but the above special case is sufficient for our purpose).

%${\bf G}$, allg $T=T_\alpha$ mit$\alpha\in{\bf G}({\Bbb Q})$ Twistung mit $\tau$
%${\rm tr}\,(\Xi^{-1}|L_\lambda)$, da $T_h$ via $h^{-1}$ auf dem Koeffizientensystem operiert (aber durch rechtstranslation  mit $h$ (nicht $h^{-1}$) im Argument $g\in{\bf G}$); s. z.B. die formeln in [Bewersdorff], S. 12 oben

\bigskip

{\bf (6.5) Remark. } 1.) We would like to emphasize that the proof of the simple trace formula of Bewersdorff is elementary. The only deeper 
ingredient is the existence of a good compactification of $\Gamma\backslash X$, which is the {\it Borel-Serre compactification}. Apart from 
that the proof only uses very general and basic principles of algebraic topology. (In [B], the formula in equation (2) only serves as a 
starting point for further investigations). 
%(cf. [B], sec. 2). 

2.) The terms appearing on the geometric side of Bewersdorff's trace formula are the archimedean components of orbital integrals on the symplectic group.

\section{Verification of Identity $(\dag)$ in section (3.7)}
 
\bigskip

We show that the family of cohomology groups $(H^d(\Gamma\backslash X,{\cal L}_\lambda))$ satisfies equation $(\dag)$ in (3.7). This is the main technical work. 
Use of Bewersdorff's trace formula reduces the verification of equation ($\dag$) to congruencs between values of irreducible characters and use of the {
\it Weyl character formula} further reduces to a problem about conjugacy of certain symplectic matrices with which we shall begin.

\medskip

We keep the notations from section 6. In particular, $\Gamma={\bf G}({\Bbb Q})\cap U\le{\bf G}({\Bbb Z})$ is an arithmetic subgroup as 
in (5.4), hence, $\Gamma\le {\cal I}\,(\le {\bf G}({\Bbb Q}_p))$ and ${\cal H}_{\cal O}={\cal H}(\Gamma\backslash\Delta /\Gamma)\otimes{\cal O}$ is the Hecke algebra defined in(5.6); in particular $T_p=T_{h_p}$, where $h_p$ is the diagonal matrix defined 
in (5.5.2).

\bigskip

{\bf (7.1) Lemma. }{\it Let $\gamma\in\Gamma$ and $h\in{\bf T}({\Bbb Z}_p)$. The matrix $h^{-1}h_p^{-e}\gamma^{-1}\in{\bf G}({\Bbb Q}_p)\le{\bf GL}_{2n}({\Bbb Q}_p)$, $e\in{\Bbb N}$, is ${\bf G}(\bar{\Bbb Q}_p)$-conjugate to a diagonal matrix
$$
\xi={\rm diag}(\xi_1,\ldots,\xi_{2n})\in{\bf T}({\Bbb Q}_p)
$$
satisfying $v_p(\alpha(\xi))=-e$ for all $\alpha\in \Delta$.

%has $2n$ pairwise different eigenvalues $\lambda_1,\ldots,\lambda_{2n}$ which are all contained in ${\Bbb Q}_p$ and which have $p$-adic valuations
%$$
%-ne,-(n-1)e,\ldots,(n-2)e,(n-1)e
%$$
%
%
%
%$$
%\lambda_i\equiv\ell^{-e_i}\gamma_{ii}^{-1}p^{-im}\pmod{p}\quad{\rm if}\quad i\le n
%$$
%and 
%$$
%\lambda_i\equiv\ell^{-e_i}\gamma_{ii}^{-1}p^{-(n+1-i)m}\pmod{p}\quad{\rm if}\;i>n.
%$$
%
%
%
%$$
%v_p(\lambda_i)=\left\{\begin{array}{cc}
%-ie&i\le n\\
%(i-n-1)e&i>n.\\
%\end{array}
%\right.
%$$

}

\medskip

{\it Proof. } We proceed in several steps.

\medskip

a.) We begin by writing $\gamma=(\gamma_{ij})$ and $h={\rm diag}(h_1,\ldots,h_{2n})$; note that the entries of $h$ as well as the 
diagonal entries $\gamma_{ii}$ of $\gamma$ are $p$-adic units. We set
$$
\underline{h}_p={\rm diag}(p^{(n-1)e},p^{(n-2)e},\ldots,p^0,p^{ne},p^{(n+1)e},\ldots,p^{(2n-1)e}),
$$
hence, $\underline{h}_p$ differs from $h_p^{-e}$ by a scalar multiple $p^{ne}$ and has integer entries. In particular, the matrices 
$h^{-1}h_p^{-e}\gamma^{-1}$ and $A=h^{-1}\underline{h}_p\gamma^{-1}$ differ by a scalar factor $p^{ne}$, hence, we may replace 
$h^{-1}h_p^{-e}\gamma^{-1}$ by $A$. We denote by 
$$
\chi(T)=T^{2n}+c_1T^{2n-1}+\cdots+c_{2n}\in{\Bbb Q}_p[T]
$$ 
the characteristic polynomial of $A$. We write $A=(a_{ij})$. Since $\gamma^{-1}\in\Gamma\le{\cal I}$, we obtain
$$
v_p(a_{ij})\ge \left\{
\begin{array}{cc}
(n-i)e&i\le n\\
(i-1)e&i>n.
\end{array}
\right.\leqno(1)
$$
More precisely, we find
$$
v_p(a_{ii})= \left\{
\begin{array}{cc}
(n-i)e&i\le n\\
(i-1)e&i>n
\end{array}
\right.\leqno(2)
$$
%$$
%a_{ii}= \left\{
%\begin{array}{cc}
%p^{(n-i)e}\gamma_{ii}^{-1}h_i^{-1}&i\le n\\
%p^{(i-1)e}\gamma_{ii}^{-1}h_i^{-1}&i>n
%\end{array}
%\right.\leqno(2)
%$$
and
$$
v_p(a_{ij})\ge \left\{
\begin{array}{cc}
(n-i)e+1&i\le n\\
(i-1)e+1&i>n
\end{array}
\right.\leqno(3)
$$
if $a_{ij}$ is one of the entries "$\mu_{ij}$, $\rho_{ij}$ or $\mu'_{ij}$" below, i.e. if ($j<i\le n$) or ($i\le n$ and $j>n$) or ($i,j>n$ and $j>i$)
$$
A=\left(
\begin{array}{ccccccc}
\tau_1&&\beta_{ij}&|&&&\\
&\ddots&&|&&\rho_{ij}&\\
\mu_{ij}&&\tau_n&|&&&\\
-&-&-&-&-&-&-\\
&&&|&\tau_{n+1}&&\mu'_{ij}\\
&*&&|&&\ddots&\\
&&&|&\beta'_{ij}&&\tau_{2n}\\
\end{array}
\right).\leqno(4)
$$
We recall that 
$$
\chi(T)={\rm det}\,(A-T{\bf I})=\sum_{\pi\in S_{2n}} {\rm sgn}(\pi)\,\prod_i (a_{i,\pi(i)}-\delta_{i,\pi(i)}T).
$$
A permutation $\pi\in S_{2n}$ contributes to $c_i$ (i.e. the summand corresponding to $\pi$ contributes in degree $T^{2n-i}$) only if $\pi$ has at 
least $2n-i$ fixed points. We denote by ${\rm Fix}(\pi)$ the set of fixed points of $\pi$ and obtain
$$
c_i=\sum_{\pi\in S_{2n}\atop|{\rm Fix}(\pi)|\ge 2n-i} \sum_{I\subseteq {\rm Fix}(\pi) \atop |I|=2n-i} {\rm sgn}(\pi)\,\prod_{i=1\atop i\not\in I}^{2n} a_{i\pi(i)}.
$$
%beachte: eine permutation $\pi$ kann mit mehreren summanden zu $c_i$ beitragen, namlich genau dann wenn $J\subset {\rm Fix}(\pi)$ eine echte teilmenge ist !!!

%If a permutation $\pi$ contributes to $\chi$ in degree $T^{2n-i}$ then there are $2n-i$ integers $k_1,\ldots,k_{2n-i}\in\{1,\ldots,2n\}$ such that $\pi(k_i)=k_i$ (there may be even more $k$ satisfying $\pi(k)=k$). We denote by $J$ the complement of the set of these $k_i$ in $\{1,\ldots,2n\}$. Hence, $J$ is $\pi$-invariant and has cardinality $i$ and the contribution of $\pi$ to the Leibniz formula is
%$$
%{\rm sgn}(\pi)\,\prod_{i\not\in J} (T-a_{ii})\; \prod_{i\in J} \left\{
%\begin{array}{cc}
%a_{i\pi(i)}&{\rm if}\;\pi(i)\not=i\\
%T-a_{ii}&{\rm if}\;\pi(i)=i.\\
%\end{array}
%\right.
%$$
In particular, the coefficient $c_i$ of $\chi$ is a sum of terms of the form
$$
\pm \prod_{i\in J} a_{i\sigma(i)},\leqno(5)
$$
where $J\subseteq\{1,\ldots,2n\}$ is a subset of cardinality $i$ and $\sigma\in {\rm Sym}(J)$ is a permutation of $J$ ($J$ is the complement of $I$; 
note that $\pi(I)=I$, hence, the complement $J$ of $I$ too is $\pi$-invariant and $\sigma=\pi|_J$ is the restriction of $\pi$ ). 
%%Es gilt $\chi(T)={\rm det}\,(A-T{\rm id})$. Um den koeffz $a_i$ zu bestimmen mittels leibniz muss man $2n-i$ mal die unbestimmte $T$ erwischen, d.h  die permutation $\pi$ in %der Leibniz formel man muss $2n-i$-mal auf die diagonale zugreifen. 
%Seien$ k_1,\ldots, k_{2n-i}$ die indices fur die $\pi(k)=k$ ist und sei $J$ das komplement davon in $\{1,%\ldots,2n\}$. dann ist $|J|=i$ und $\pi$ schrankt sich zu einer permutation $\sigma$ von $J$ ein und kann dort eine bliebige permutation sein. Die summe der beitrage uber alle %$J$ der machtigkeit $i$ und alle $\sigma\in{\rm Sym}(J)$ liefert dann den koeffz von $a_i$. 

\medskip

b.) We select $i\in\{1,\ldots,2n\}$ and look closer at the coefficient $c_i$. We first assume $i\le n$ and we define the subset 
$J_{\rm min}=\{n-i+1,n-i+2,\ldots,n\}$; note that $|J_{\rm  min}|=i$ and let $\sigma\in {\rm Sym}(J_{\rm min})$. 
If $\sigma={\rm id}_{J_{\rm min}}$ then (5) yields a term which has $p$-adic value $e\,\sum_{k=0}^{i-1} k$ by equation (2). 
If $\sigma\not={\rm id}_{J_{\rm min}}$ then $\sigma$ picks at least one entry "$\mu_{ij}$" below the diagonal (cf. equation (4)) which is divisible 
by one more $p$ (cf. equation (3)) and, hence, the $p$-adic value of the term in equation (5) is larger than $e\,\sum_{k=0}^{i-1} k$. 
Thus, 
$$
\sum_{\sigma\in{\rm Sym}(J_{\rm min})} (\pm 1) \prod_{i\in J_{\rm min}} a_{i\sigma(i)}
$$
has $p$-adic value $e\,\sum_{k=0}^{i-1} k$. If $J\not=J_{\rm min}$ with $|J|=i$, then equation (1) implies that $\prod_{i\in J} a_{i\sigma(i)}$ has $p$-adic value bigger than $e\,\sum_{k=0}^{i-1} k$ for any $\sigma\in{\rm Sym}(J)$. Thus, we obtain
$$
v_p(c_i)=e\,\sum_{k=0}^{i-1} k,\quad i=1,\ldots,n.
$$
Next we assume $i>n$ and we define the subset $J_{\rm min}=\{1,\ldots,n,\ldots,i\}$. We claim: if $\sigma\in{\rm Sym}(J_{\rm min})$ is not the identity then there is $i_0\in J_{\rm min}$ such that 
$$
v_p(a_{i_0,\pi(i_0)})\ge 
 \left\{
\begin{array}{ccc}
(n-i_0)e+1&\mbox{if}&i_0\le n\\
(i_0-1)e+1&\mbox{if}&i_0>n.\\
\end{array}
\right.\leqno(6)
$$
To prove the claim we assume that $\sigma\in{\rm Sym}(J_{\rm min})$ is a permutation such that equation (6) does not hold. Equation (3) 
then implies that $\sigma(i)\le n$ for all $i\le n$ ("$\rho_{ij}$" has $p$-adic value greater than or equal to $(n-i)e+1$). 
Thus, $\sigma$ maps $\{1,\ldots,n\}$ to itself and also maps $\{n+1,\ldots,i\}$ to itself, i.e. $\sigma$ defines permutations $\sigma|_{\{1,\ldots,n\}}$ resp. $\sigma|_{\{n+1,\ldots,i\}}$ of $\{1,\ldots,n\}$ resp. of $\{n+1,\ldots,i\}$. Since equation (6) does not hold, equation (3) further implies that 
$\sigma|_{\{1,\ldots,n\}}$ and $\sigma|_{\{n+1,\ldots,i\}}$ are the identity. Hence, $\sigma$ is the identity which proves the claim. As above we then obtain for all $\sigma\in{\rm Sym}(J_{\rm min})$, $\sigma\not={\rm id}$, that
$$
v_p(\prod_{j\in J_{\rm min}} a_{i,\sigma(i)}) >v_p(\prod_{i\in J_{\rm min}} a_{i,i})=e\,\sum_{k=0}^{i-1} k.
$$
Thus,
$$
v_p(\sum_{\sigma\in{\rm Sym}(J_{\rm min})}(\pm 1) \prod_{i\in J_{\rm min}} a_{i\sigma(i)})=e\,\sum_{k=0}^{i-1} k,
$$
As in the case $i\le n$ equation (1) implies that for any $J\not=J_{\rm min}$, $|J|=i$, and any $\sigma\in{\rm Sym}(J)$ 
$$
v_p(\prod_{i\in J} a_{i\sigma(i)})>e\,\sum_{k=0}^{i-1} k.
$$
Thus,
$$
v_p(c_i)=e \sum_{k=0}^{i-1} k
$$
for all $i=n+1,\ldots,2n$. Hence, this holds for all $i=1,\ldots,2n$.

\medskip

%We view $\chi$ as polynomial with coefficients in ${\Bbb Q}_p$ (using the embedding $i$). 
c.) In particular, the Newton polygon of $\chi$ consists of $2n$ segments which have slopes $0,e,2e,\ldots,(2n-1)e$. 
%Beachte: wir haben die koeffz des (char) polynomes genau umgekehrt zu [Neukirch: Alg Zahlentheorie], II,6 p. 151 durchnumeriert, d.h. unser Newton polgon lauft in die andere richtung als das von Neukirch, bzw. die $p$-adischen bewertungen der NS sind =den slopes nicht =-slopes.
Thus, there are $2n$ roots $\lambda_1',\ldots,\lambda_{2n}'\in\bar{\Bbb Q}_p$ of $\chi$ which have $p$-adic valuations 
$0,e,2e,\ldots,(2n-1)e$. Since $h^{-1}h_p^{-e}\gamma^{-1}$ and $A$ differ by a scalar factor $p^{ne}$, we deduce that $h^{-1}h_p^{-e}\gamma^{-1}$ has $2n$ pairwise different eigenvalues $\lambda_1,\ldots,\lambda_{2n}$ with $p$-adic values $-ne,-(n-1)e,\ldots,(n-1)e$. 

\medskip

d.) We claim that $\chi$ splits over ${\Bbb Q}_p$ as a product of linear factors. In fact, if $\chi$ does not split completely then there is an 
irreducible factor of $\chi$ of degree $\ge 2$, hence, there are roots $\lambda_i,\lambda_j$ which are conjugate by an automorphism 
$\sigma\in{\rm Aut}(\bar{\Bbb Q}_p/{\Bbb Q}_p)$. This would imply that $\lambda_i,\lambda_j$ have the same $p$-adic value which is 
a contradiction. Thus, $\chi$ is split over ${\Bbb Q}_p$, hence, all roots $\lambda_i$ are contained in ${\Bbb Q}_p$.

\medskip

%Finally, we consider the general case, i.e. we drop the assmption $p\not| h_i$. We write $h=h_{(p)} h^{(p)}$ where the entries of $h^{(p)}$ 
%are prime to $p$ and $h_{(p)}={\rm diag}(p^{e_1},\ldots,p^{e_n},p^{v-e_1},p^{v-e_{2n}})$ with $e_1\le e_2\le\cdots\le e_n$ ($v=v_p(\nu)$). 
%Repeating the above proof with $h$ replaced by $h^{(p)}$ and $h_p^{-m}$ replaced by $h_{(p)}^{-1}h_p^{-m}$ we obtain the result in the general case. Thus, the proof of the lemma is complete

\medskip

\medskip

e.) In c.) we have seen that the image of $h^{-1}h_p^{-e}\gamma^{-1}$ under ${\bf G}({\Bbb Q}_p)\subseteq{\bf GL}_{2n}({\Bbb Q}_p)$ 
has $2n$ different eigenvalues. Thus, $h^{-1}h_p^{-e}\gamma^{-1}$ is a semi simple element in ${\bf GL}_{2n}({\Bbb Q}_p)$ and, hence, 
in ${\bf G}({\Bbb Q}_p)$. 
%cf [S], 2.4.8 Theorem ii), p. 34; bea.: unser $\Phi$ ist die einbettung ${\bf G}({\Bbb Q}_p)\subseteq{\bf GL}_{2n}({\Bbb Q}_p)$ die injektiv ist; deshalb folgt fur $g$ mit $\Phi(g)$ h.e.: $\Phi(g)=\Phi(g)_s=\Phi(g_s)\Rightarrow $g=g_s$, d.h. $g$ ist h.e.
In particular, $h^{-1}h_p^{-e}\gamma^{-1}$ is (${\bf G}(\bar{\Bbb Q}_p)$-)conjugate to an element
$$
\xi={\rm diag}(\xi_1,\ldots,\xi_{2n})
$$
in ${\bf T}(\bar{\Bbb Q}_p)$. Since conjugate matrices in ${\bf GL}_{2n}({\Bbb Q}_p)$ have the same eigenvalues, c.) implies that 
$\xi\in{\Bbb Q}_p$ and the $p$-adic values of the $\xi_i\in{\Bbb Q}_p$ are contained in the sequence $-ne,-(n-1)e,\ldots,(n-1)e$. Conjugating 
the regular element $\xi\in{\bf T}({\Bbb Q}_p)$ with a suitable element in the Weyl group ${\cal W}$ of ${\bf G}$ we may assume 
that $v_p(\alpha(\xi))<0$ for all $\alpha\in\Delta$. 
%
%The Weyl group of ${\bf Sp}_n$ is a semidirect product $W=S_n\times ({\Bbb Z}/2{\Bbb Z})^n$, where $\tau\in S_n$ acts on $\xi$ by permuting the entries $\xi_1,\ldots,\xi_n$ and $\tau_i\in({\Bbb Z}/2{\Bbb Z})^n$ permutes $\xi_i$ and $\xi_{n+i}$. Since 
%$$
%v_p(\xi_i)=-e-v_p(\xi_{n+i})
%$$ 
%(notice that ${\rm det}\,h_p=n$, hence, ${\rm det}\,h^{-1}h_p^{-e}\gamma^{-1}=p^{-en}$ and the multiplier $\nu$ of $h^{-1}h_p^{-e}\gamma^{-1}$ therefore has $p$-adic value $v_p(\nu)=-e$) we know that exctly one of $\xi,\xi_{n+i}$ has $p$-adic value $<0$ (note that $v_p(\xi)$, $v_p(\xi_{n+i})$ are divisible by $e$). By applying a suitable $\tau_i$ we can achieve that the entries $\x_1,\ldots,\xi_n$ all have $p$-adic value $<0$. By applying a permutation $\sigma\in S_n$ we then can achieve that $0>v_p(\xi_1)> \cdots >v_p(\xi_n)$.
%
This then implies that
$$
(v_p(\xi_1),\ldots,v_p(\xi_{2n}))=(-e,-2e,\ldots,-ne,0,e,\ldots,(n-1)e)
$$
which shows that $v_p(\alpha(\xi))=-e$ for all simple roots $\alpha$. 
%1.) bea.: $v_p(-2\epsion_1(\xi))>0$ erzwingt $v_p(x_1)<0$ und da $v_p(ßalpha(\xi))>$ folgt dann dass $v_p(\xi_i)<0$ fuir all $i=1,\ldots,n$. das erzwingt dann dass $v_p(x_i)$ wir oben angegeben sind
%
%2.) bea.: $\xi$ hat determinante $p^{-en}$, d.h. um die wurzel $-2\epsilon_1$ von $\xi$ zu bestimmen mussen wir $\xi$ durch der skalarmatrix $p^{e/2}$ mulitplizieren; das ergibt dann $-2\epsipon_1(\xi)=e$. 
This completes the proof of the Lemma.

\bigskip

{\it Remark. } The entries $\xi_i$ of $\xi$ are even algebraic integers (contained in ${\Bbb Q}_p$).
%but not in ${\Bbb Q}$

%VVVVVVVVVVVVVVV brauch ich das ? VVVVVVVVVVVVVVV
%Since the semi simple element $h^{-1}h_p^{-e}\gamma^{-1}\in{\bf G}({\Bbb Q}_p)$ is contained in a maximal ${\Bbb Q}_p$-torus (cf. [S], Corollary 13.3.8, p. 231) and since any two maximal ${\Bbb Q}_p$-tori in ${\bf G}$ are ${\bf G}({\Bbb Q}_p)$-conjugate (cf. [S], 14.4.3 Theorem, p. 250) we finally obtain that $h^{-1}h_p^{-e}\gamma^{-1}$ is ${\bf G}({\Bbb Q}_p)$-conjugate to an element in ${\bf T}({\Bbb Q})$ which must be $\xi$. This complete the proof of the Corollary.
%AAAAAAAAAAAAAAAAAAAAAAAA

\bigskip

\bigskip

{\bf (7.2) Lemma. } {\it Let $\lambda\in X({\bf T})$ be algebraic and dominant. 
%, i.e. the restriction $\lambda^\circ=\lambda|_{{\bf T}^0}$ of $\lambda$ to ${\bf T}^0$ is a dominant character of the semi simple group ${\bf G}^0$ and $\lambda|_{\bf Z}$ is a (algebraic) character 
Then, for any $\gamma\in \Gamma$ and $h\in\Sigma^+\,(\subseteq{\bf T}({\Bbb Q}))$ we have
$$
\lambda(hh_p^e)\,{\rm tr}(\pi_\lambda(h^{-1}h_p^{-e}\gamma^{-1})|L_\lambda)\in{\Bbb Z}.
$$
}

\medskip

{\it Proof. } Since $hh_p^e\gamma\in{\bf G}({\Bbb Q})$ we know that $\pi_\lambda(hh_p^e\gamma)$ leaves $L_\lambda({\Bbb Q})$ invariant. 
The subspace $L_\lambda({\Bbb Z})$
%definiert als ${\Bbb Z}$-struktur auf $L_\lambda$; vgl. section (5.3)
decomposes
$$
L_\lambda({\Bbb Z})=\bigoplus_\mu L_\lambda(\mu,{\Bbb Z}),
$$
where $\mu\in X({\bf T})$ runs over all weights of the form 
$$
\mu=\lambda-\sum_{\alpha\in\Delta} c_\alpha\alpha\leqno(7)
$$ 
with $c_\alpha\in{\Bbb N}_0$ for all $\alpha\in \Delta$ and where ${\bf T}({\Bbb Q})$ acts on $L_\lambda({\mu},{\Bbb Z})\otimes{\Bbb Q}$ via $\mu$:
$$
t v_\mu=\mu(t)v_\mu
$$
for all $t\in {\bf T}({\Bbb Q})$ and all $v_\mu\in L_\lambda(\mu,{\Bbb Z})\otimes{\Bbb Q}$. 
%$L_\lambda(\mu,{\Bbb Z})$ is the intersection of the weight $\mu$ subspace $L_\lambda(\mu,{\Bbb Q})\le L_\lambda({\Bbb Q})$ with $L_\lambda({\Bbb Z})$, 
%cf. [Humphries, Lia alg], 27.1 Theorem, p. 158; represent ${\bf G}$ as Chevalley group, hence, ${\bf G}({\Bbb Z})\le {\mathfrak U}_{\Bbb Z}^\times$.
In particular, if $\mu$ is as in equation (7) then we obtain
$$
\lambda(hh_p^e)\pi_\lambda(h^{-1}h_p^{-e})v_\mu=\prod_{\alpha\in\Delta} \alpha(h h_p^e)^{c_\alpha}.\leqno(8)
$$

%Since $\alpha(hh_p^e)\in{\Bbb Z}_p$ and since any $\gamma^{-1}\in\Gamma$ leaves $L_\lambda({\Bbb Z})$ invariant we deduce that
%$$
%\lambda(hh_p^e)\pi_\lambda(h^{-1}h_p^{-e}\gamma^{-1})L_\lambda({\Bbb Z}_p)\subseteq L_\lambda({\Bbb Z}_p).\leqno(invariant=7)
%$$
%for all $\gamma\in\Gamma$ and all $h\in \Sigma_\ell$ with integer coefficients. If $h\in {\bf T}({\Bbb Q}_p)^{++}\cap M_{2n}({\Bbb Z})$ then even obtain
%$$
%\lambda(hh_p^e)\pi_\lambda(h^{-1}h_p^{-e}\gamma^{-1})L_\lambda({\Bbb Z})\subseteq L_\lambda({\Bbb Z}).
%$$

Since $h\in\Sigma^+$ we obtain $\alpha(h)\in {\Bbb Z}$ and equation (8) implies that $\lambda(hh_p^e)\,\pi_\lambda(h^{-1}h_p^{-e}\gamma^{-1})$ 
leaves ${L_\lambda}({\Bbb Z})$ invariant which yields the claim (note that $\gamma\in{\bf G}({\Bbb Z})$). Thus, the proof of the lemma is complete.

%Etwas allgemeiner konnen wir folgendes zeigen:
%Thus, if $h\in {\bf T}({\Bbb Q}_p)^+$ resp. $\in{\bf T}({\Bbb Q})^+$ resp. $\in \Sigma^+$ then $\alpha(h)\in {\Bbb Z}_p$ resp. $\in {\Bbb %Z}_{(p)}$ resp. $\in{\Bbb Z}$ for all $\alpha\in \Delta$. Hence, equation (8) implies that 
%$\lambda(hh_p^e)\,\pi_\lambda(h^{-1}h_p^{-e}\gamma^{-1})$ leaves ${L_\lambda}({\Bbb Q}_p)$ resp. ${L_\lambda}({\Bbb Z}_{(p)})$ resp. %${L_\lambda}({\Bbb Z})$ invariant which yields 
%$$
%\lambda(hh_p^e)\,{\rm tr}(\pi_\lambda(h^{-1}h_p^{-e}\gamma^{-1})|L_\lambda)\in
%\left\{\begin{array}{ccc}
%{\Bbb Z}_p&\mbox{if}&h\in {\bf T}({\Bbb Q}_p)^+\\
%{\Bbb Z}_{(p)}&\mbox{if}&h\in {\bf T}({\Bbb Q})^+\\
%{\Bbb Z}&\mbox{if}&h\in \Sigma^+.\\
%\end{array}\right.
%$$

\bigskip

{\bf (7.3)  } We denote by $h_\alpha$ the coroot corresponding to the root $\alpha$ and for any root $\alpha$ and any
 $\lambda\in X({\bf T})$ we set $\langle \lambda,\alpha\rangle=\lambda(h_\alpha)$. 
%$\langle \lambda,\alpha\rangle=2\frac{(\lambda,\alpha)}{(\alpha,\alpha)}=2\frac{(t_\lambda,t_\alpha)}{(t_\alpha,t_\alpha)}=2\frac{\lambda(t_\alpha)}{(t_\alpha,t_\alpha)}=\lambda(h_\alpha)$
%s. auch [Springer; alg. groups], p. 124
We further denote by $\omega_\alpha$ the fundamental dominant weights corresponding to the basis $\Delta$ and ${\cal W}$ is the Weyl group 
of ${\bf G}$. 
%
%
%
%and by $s_\alpha$ the reflection corresponding to the root $\alpha$. Thus,
%$$
%s_\alpha(\lambda)=\lambda-\langle \lambda,\alpha  \rangle\alpha.
%$$
%and $W$ is generated by the reflections $s_\alpha$ corresponding to simple roots $\alpha\in \Delta$. 
%%cf. [Humphreys, Lie Alg] 10.3 Theorem d), p. 51
We also write $\rho$ for the half sum of the positive roots and we put $w\cdot\lambda=w(\lambda+\rho)-\rho$.  
%????????? Finally we denote by $[\Xi]=[\Xi]_\Gamma$ the $\Gamma$-conjugacy class of $\Xi$. ???????????????

\bigskip

{\bf Lemma. } {\it Let $\lambda\in X({\bf T})$ be a dominant character. For any $w\in {\cal W}$, $w\not=1$, we obtain
$w\cdot\lambda=\lambda-\sum_{\alpha\in \Delta} c_\alpha \alpha$, where $c_\alpha=c_{\alpha,w}\in{\Bbb N}_0$ and 
$$
c_{\alpha_0}\ge  \frac{\langle\lambda,\alpha_0\rangle}{2}
$$ 
for at least one root $\alpha_0\in\Delta$. 
}

\medskip

{\it Proof. } Since $w\lambda$ is a weight of the irreducible ${\bf G}$-module of highest weight $\lambda$ we know that
$$
w\lambda=\lambda-\sum_{\alpha\in\Delta} b_\alpha\alpha
$$
for certain $b_\alpha\in{\Bbb N}_0$. Since $\lambda$ is dominant we may write $\lambda=\sum_{\alpha\in\Delta} d_\alpha\omega_\alpha$ where 
$d_\alpha=\langle \lambda,\alpha\rangle\in{\Bbb N}_0$. On the other hand, $w\not=1$ implies that $w\lambda$ is not contained in the Weyl chamber 
corresponding to the basis $\Delta$, hence, $\langle w\lambda,\alpha_0\rangle\le 0$ for some root $\alpha_0\in\Delta$.
We obtain
$$
0\ge \langle w\lambda,\alpha_0 \rangle=\langle \sum_{\alpha\in \Delta} d_\alpha \omega_\alpha -\sum_{\alpha\in\Delta} b_\alpha \alpha,\alpha_0 \rangle=d_{\alpha_0}-\sum_{\alpha\in\Delta} b_\alpha \langle \alpha,{\alpha_0}\rangle.
$$
Since $\langle\alpha,\alpha_0\rangle=\alpha(h_{\alpha_0})=2$ if $\alpha=\alpha_0$ and $\langle\alpha,\alpha_0\rangle\le 0$ if $\alpha\not=\alpha_0$ 
%s. Table 1 in [Humphries], sec. 11.4, p. 59 (der fall $C_n$)
this yields $0\ge  d_{\alpha_0}-2b_{\alpha_0}$. Thus, 
$$
b_{\alpha_0}\ge \frac{1}{2} d_{\alpha_0}=\frac{1}{2}\langle \lambda,\alpha_0\rangle.
$$
Since $w\cdot\lambda=w\lambda+w\rho-\rho$ and
$$
w\rho-\rho=-\sum_{\alpha\in\Phi^+\atop \alpha\in w\Phi^-}\alpha
$$
%details dazu:
%By definition of $\rho$ we obtain 
%$$
%w\rho-\rho=\frac{1}{2}\left(\sum_{\alpha\in\Phi^+\atop w\alpha\in \Phi^-}w \alpha-\sum_{\beta\in\Phi^+\atop \beta\not\in w\Phi^+}\beta\right)
%=\frac{1}{2}\left(\sum_{\alpha\in\Phi^-\atop \alpha\in w\Phi^+}\alpha-\sum_{\beta\in\Phi^+\atop \beta\in w\Phi^-}\beta\right).
%$$
%The assignment $\alpha\mapsto -\alpha$ is a bijection form the index set of the first summand to the index set of the second summand, hence, we obtain
%$$
%w\rho-\rho=\frac{1}{2}\left(-\sum_{\alpha\in\Phi^+\atop \alpha\in w\Phi^-}\alpha-\sum_{\beta\in\Phi^+\atop \beta\in w\Phi^-}\beta\right)
%=-\sum_{\alpha\in\Phi^+\atop \alpha\in w\Phi^-}\alpha.
%$$
this yields the claim and the Lemma is proven.

\bigskip

{\bf (7.4) Congruences between characters values for different weights. } Using the {\it Weyl character formula} we obtain the following congruences 
between character values of irreducible algebraic representations of ${\bf G}$ for varying highest weights.

\bigskip

{\bf Proposition. }{\it Let $\lambda,\lambda'\in X({\bf T})$ be dominant characters and let $C\in{\Bbb Q}_{>0}$ such that 
$\langle\lambda,\alpha\rangle>2C$ and $\langle\lambda',\alpha\rangle>2C$ for all $\alpha\in\Delta$. If $\lambda\equiv\lambda'\pmod{(p-1)p^m X({\bf T})}$, 
then for any $\Gamma$-conjugacy class $[\Xi]_\Gamma\subseteq\Gamma hh_p^e\Gamma$, $e\in{\Bbb N}$, $h\in\Sigma^+$ the following congruence holds
$$
\lambda(hh_p^{e}) {\rm tr}\,(\Xi^{-1}|{L_\lambda})\equiv \lambda'(hh_p^e) {\rm tr}\,(\Xi^{-1}|{L_{\lambda'}})\pmod{p^{{\rm min}(m+1,Ce)}}.
$$
}

\medskip

{\it Proof. } We note that by (7.2) Lemma 
$$
\lambda(hh_p^{e}) {\rm tr}\,(\Xi^{-1}|{L_\lambda}),\, \lambda'(hh_p^e) {\rm tr}\,(\Xi^{-1}|{L_{\lambda'}})\in{\Bbb Z}.
$$
We may assume that the representative $\Xi$ of the $\Gamma$-conjugacy class $[\Xi]_\Gamma\subseteq \Gamma hh_p^e\Gamma$ is 
of the form $\Xi=\gamma h_p^eh$ for some $\gamma\in\Gamma$. By definition of $\Sigma^+$ we may write $h=h_{(p)}h_p^c$ with 
$h_{(p)}\in\prod_{\ell\not|Np}\Sigma_\ell^+\subseteq {\bf T}({\Bbb Z}_p)$ and $c\in{\Bbb N}_0$. Hence, 
$$
\Xi=\gamma h_{(p)}h_p^{e'}
$$
with $e'=e+c\ge e(>0)$. Thus, (7.1) Proposition implies that $\Xi^{-1}$ is ${\bf G}(\bar{\Bbb Q}_p)$-conjugate to an element $\xi\in {\bf T}({\Bbb Q}_p)$ 
satisfying 
$$
v_p(\alpha(\xi))=-e'\leqno(9)
$$
for all $\alpha\in\Delta$. 
%wir sind immer uber ${\Bbb Q}$ bzw. ${\Bbb Z}$; erst wenn $\xi$ ins Spiel kommt haben wir zum ersten Mal und unvermeidlich Koordinaten in ${\Bbb Q}_p$ (namlich die Eintrage von $\xi$, wharend die eintrage von $h,h_p,\gamma$ in  ${\Bbb Q}$ bzw. ${\Bbb Z}$ liegen.
Using the {\it Weyl character formula} we therefore obtain
$$
\Delta(\xi)\cdot{\rm tr}\, (\Xi^{-1}|{L_\lambda})={\sum_{w\in W} (-1)^{\ell(w)}\,(w\cdot\lambda)(\xi)},
$$
where
$$
{\Delta}(\xi)=\prod_{\alpha\in\Phi^+} (1-\alpha^{-1}(\xi)).
$$
%cf. [Springer, Encycl. Math Sciences 55: Algebraic geo IV], I. 4.6.4, p. 45 
%s. auch [Humphries, Lie Alg], sec. 22.5, p. 124
% hier sind die details:
%
%\begin{itemize}
%
%\item Es gilt ${\rm tr}(h|L_\lambda)=\sum_\mu m(\mu) \mu(h)$ wobei $m(\mu)={\rm dim}\, L_\lambda(\mu)$ und $\mu$ uber alle gewichte in $L_\lambda$ lauft
%
%\item Sei $X=X({\bf T})$ die Charaktergruppe von ${\bf T}$ und ${\Bbb Z}[X]$ der Gruppenring zu $X$; $e^\mu$ sei das basis element in ${\Bbb Z}[X]$ zu $\mu$; nach definition der Muliplikation in %${\Bbb Z}[X]$ gilt also $e^\mu e^\tau=e^{\mu+\tau}$ ("+" ist die verknupfung in $X$ ($X$ ist additiv geschrieben)). Wir definieren fur jedes $h\in {\bf T}(\bar{\Bbb Q})$ eine Abbildung
%$$
%\begin{array}{cccc}
%\Psi_h:&{\Bbb Z}[X]&\rightarrow& \bar{\Bbb Q}\\
%&\sum_\mu c_\mu e^\mu&\mapsto&\sum_\mu c_\mu \mu(h)
%\end{array}
%$$
%(Auswertung bei $h$)
%
%\item Es gilt in ${\Bbb Z}[X]$: $\Delta\cdot\sum_\mu m(\mu) e^\mu=\sum_w{\rm sgn}(w) e^{w\cdot\lambda}$, $\Delta=\prod_{\alpha>0}(1-\alpha^{-1}(h))$ (benutze die Formel fur $a(\rho)$ in [Springer], p. 45 mitte)
%
%\item Anwenden von $\Psi_h$ liefert
%$$
%\Delta(h) \cdot\sum_\mu m(\mu) \mu(h)=\sum_w{\rm sgn}(w) w\cdot\lambda(h)
%$$
%
%\item Damit folgt ${\rm tr}(h|L_\lambda)=\sum_w{\rm sgn}(w) (w\cdot)\lambda(h)/\prod_\alpha(1-\alpha^{-1}(h))$
%
%\end{itemize}
%
%s. auch [Rossman], sec. 6.7 Theorem 9, p. 242; definition of $e^\alpha$ on p. 108). 
%
%
Note that equation (9) implies that $\Delta(\xi)$ is a $p$-adic unit, hence, $\Delta(\xi)\not=0$. Using (7.3) Lemma we can write 
$$
w\cdot \lambda=\lambda-\sum_{\alpha\in\Delta} c_{\alpha,w} \alpha
$$ 
with $c_{\alpha,w}\in{\Bbb N}_0$ and $c_{\alpha_w,w}\ge \langle \lambda,\alpha_w \rangle/2$ for some root $\alpha_w\in\Delta$. We obtain
$$
\lambda(hh_p^e)\,{\rm tr}\, (\Xi^{-1}|{L_\lambda})=\frac{\lambda(hh_p^e\xi)}{\Delta(\xi)} \,\left(1+\sum_{w\not=1} {\rm sgn}(w)\prod_{\alpha\in\Delta} \alpha(\xi)^{-c_{\alpha,w}}\right).
$$
Since $v_p(\alpha^{-1}(\xi))=e'\ge 1$ for all $\alpha\in \Delta$ we find that $\Delta(\xi)\in{\Bbb Z}_p$ is a $p$-adic unit. Moreover, for any $\alpha\in\Delta$ we have $v_p(\alpha(hh_p^e))=e'$, hence,
$$
v_p(\alpha(hh_p^e))=-v_p(\alpha(\xi))
$$
for all $\alpha\in\Delta$. In particular, this equality holds for all $\beta$ contained in the root lattice of ${\bf G}$ and since a (integral) multiple of any integral weight is contained in the root lattice we obtain
$$
v_p(\chi(hh_p^e))=-v_p(\chi(\xi))
$$
for all $\chi\in X({\bf T})$. Thus, $\chi(hh_p^e\xi)\in{\Bbb Z}_p^*$ is a $p$-adic unit. 
%bea.: $h,h_p,\xi\in{\bf T}({\Bbb Q}_p)$ $\Rightarrow \chi(hh_p^e\xi)\in{\Bbb Q}_p$
In particular, $\lambda(hh_p^e\xi)$ is a $p$-adic unit. Taking into account that $c_{\alpha,w}\ge 0$ for all $\alpha\in\Delta$, $w\in {\cal W}$ and 
that $c_{\alpha_w,w}\ge \langle \lambda,\alpha_w \rangle/2\ge C$ we thus obtain using equation (9)
$$
\lambda(hh_p^e)\,{\rm tr}\, (\Xi^{-1}|{L_\lambda})\equiv \frac{\lambda(hh_p^e\xi)}{\Delta(\xi)}\pmod{p^{Ce'}{\Bbb Z}_p}.\leqno(10)
$$
%
%
%
%Hence,
%$$
%\lambda(hh_p^m)\,{\rm tr}\, A|_{L_\lambda}-\lambda'(hh_p^m)\,{\rm tr}\, A|_{L_{\lambda'}}\equiv \lambda(hh_p^m\xi)-\lambda'(hh_p^m\xi)\pmod{p^{Cm}}.
%$$
%
%
%
Since $\lambda\equiv\lambda'\pmod{(p-1)p^m X(\bf T)}$ there is a $\chi\in X(\bf T)$ such that $\lambda-\lambda'=(p-1)p^m\chi$. Taking into account that 
$\chi(hh_p^e\xi)$ is a $p$-adic unit this yields
$$
\frac{\lambda(hh_p^e\xi)}{\lambda'(hh_p^e\xi)}
%=(\lambda-\lambda')(hh_p^e\xi)
=\chi(hh_p^e\xi)^{(p-1)p^m}\in 1+p^{m+1}{\Bbb Z}_p.
$$
%by little fermat
Hence,
$$
\lambda(hh_p^e\xi)\equiv\lambda'(hh_p^e\xi)\pmod{p^{m+1}{\Bbb Z}_p}.
$$
%denn $$
%\lambda(hh_p^e\xi)-\lambda'(hh_p^e\xi)\in p^{m+1}{\Bbb Z}_p.
%$$
%
%
%
%
%
%equation (8) in addition implies that $(\lambda-\lambda')(hh_p^e\xi)\in p^{m+1}{\Bbb Z}_p$.
%denn:
%$$
%(\lambda-\lambda')(hh_p^e\xi)=\prod_{i=1}^n\delta_i^{z_i}
%$$ 
%with $z_i\in(p-1)p^m{\Bbb Z}$. Hence, $(\lambda-\lambda')(hh_p^e\xi)\in p^{m+1}{\Bbb Z}_p$. 
%
%
%
Together with equation (10) which also holds with $\lambda$ replaced by $\lambda'$ we obtain
$$
\lambda(hh_p^e)\,{\rm tr}\, (\Xi^{-1}|{L_\lambda})\equiv\lambda'(hh_p^e)\,{\rm tr}\, (\Xi^{-1}|{L_{\lambda'}})\pmod{p^{{\rm min}(m+1,Ce')}{\Bbb Z}_p}.
$$
Since $e'\ge e$ this completes the proof.

\bigskip

{\bf (7.5) Congruences for Hecke operators in different weights. } Let $\beta\in{\Bbb Q}_{\ge 0}$. For any pair of dominant characters 
$\lambda,\lambda'\in X({\bf T})$ we denote by ${\bf e}_{\lambda,\lambda'}={\bf e}_{{\bf H}_\lambda,{\bf H}_{\lambda'}}^{\le\beta}$ the approximate idempotent projecting to
the slope $\le\beta$ subspaces of ${\bf H}_\lambda$ and ${\bf H}_{\lambda'}$ as defined in (3.2); $\{{\bf e}_{\lambda,\lambda'}\}_\lambda$ 
resp. $\{{\bf e}_{\lambda,\lambda'}\}_{\lambda'}$ then is the approximate 
idempotent projecting to the slope subspaces ${\bf H}_\lambda^{\le\beta}$ and ${\bf H}_{\lambda'}^{\le\beta}$ which are now 
defined with respect to the normalized action of $T_p\in{\cal H}$ (cf. (6.3)).

\bigskip

{\bf Theorem. }{\it Let $C\in{\Bbb Q}_{>0}$. Assume that the dominant characters $\lambda,\lambda'\in X({\bf T})$ satisfy 
\begin{itemize}

\item $\langle \lambda,\alpha\rangle>2C$ and $\langle \lambda',\alpha\rangle>2C$ for all $\alpha\in \Delta$. 

\item $\lambda\equiv\lambda'\pmod{(p-1)p^m X({\bf T})}$.

\end{itemize} 

Then for all Hecke operators $T\in {\cal H}_{\cal O}$ and all slopes $\beta\in{\Bbb Q}_{\ge 0}$ the following congruence holds:
$$
{\rm tr}\,(\{{\bf e}_{\lambda,\lambda'}^{\lceil\frac{m+1}{C}\rceil}T\}_\lambda|{\bf H}_\lambda)\equiv {\rm tr}\,(\{{\bf e}_{\lambda,\lambda'}^{\lceil\frac{m+1}{C}\rceil}T\}_{\lambda'}|{\bf H}_{\lambda'})\pmod{p^{\Box}},
$$
where
$$
\Box=(1-\frac{\beta M(\beta)}{C})(m+1)-\beta M(\beta).
$$
}
 
\medskip

{\it Proof. } Since ${\cal H}_{\cal O}$ is generated as ${\cal O}$-module by the Hecke operators $T_h$, $h\in\Sigma^+$ (cf. (5.6)), 
we may assume that $T=T_{h}$ for some ${h}\in\Sigma^+$. 
We write $h=h_p^fh_{(p)}$ with $f\in{\Bbb N}_0$ and $h_{(p)}\in\prod_{\ell\not|Np}\Sigma_\ell^+$ and we set $L=\lceil\frac{m+1}{C}\rceil$. We 
recall that ${\bf e}_{\lambda,\lambda'}=p(T_p)$, where the polynomial $p=\sum_{e=1}^{t}c_e X^e$ satisfies the following properties: 
%its coefficients $c_e$ lie in $E_{H_\lambda,H_{\lambda'}}^{\le\beta}$ (the common splitting field of $T_p$ acting on 
%$H_\lambda^{\le\beta}$ and $H_{\lambda'}^{\le\beta}$), 
its degree $t$ is bounded by $M(\beta)$, its 
constant term $p(0)=0$ and ${\bf S}(p)\ge -\beta$ (cf. (3.3) Lemma). Since $\{{\bf e_{\lambda,\lambda'}}\}_\lambda=p(\{T_p\}_\lambda)$ we obtain
$$
\{{\bf e}_{\lambda,\lambda'}^LT_h\}_\lambda=\{{\bf e}_{\lambda,\lambda'}\}_\lambda^L \{T_h\}_\lambda=p^L(\{T_p\}_\lambda)\{T_h\}_\lambda=\sum_{e=L}^{tL} b_e \{T_p\}_\lambda^e\{T_h\}_\lambda,
$$
where $v_p(b_e)\ge -e\beta$. Since $\{T_p\}_\lambda^e=\lambda(h_p^e)T_p^e$, $\{T_h\}=\lambda(h_p^f)T_h$ and $T_hT_p^e=T_hT_{h_p^e}=T_{hh_p^e}$ 
%(cf. equation (6) in section (6.4); note that $h\in\Sigma^+$) 
this yields
$$
\{{\bf e}_{\lambda,\lambda'}^LT_h\}_\lambda=\sum_{e=L}^{tL} b_e \lambda(h_p^{e+f}) T_{hh_p^e}.
$$
Applying the Topological trace formula of Bewersdorff (cf. (6.4) Theorem) we obtain
\begin{eqnarray*}
(3)\qquad{\rm tr}\,(\{T_h {\bf e}_{\lambda,\lambda'}^L\}_\lambda|{\bf H}_{\lambda})&=&\sum_{e=L}^{tL} b_e\,\lambda(h_p^{e+f})\,{\rm tr}\,(T_{hh_p^e}|{H^d(\Gamma\backslash X,{\cal L}_\lambda)})\\
&=&\sum_{e=L}^{tL} b_e \,\lambda(h_p^{e+f})  \sum_{[\Xi]\in\Gamma hh_p^e\Gamma/\sim_\Gamma} c_{[\Xi]}\, {\rm tr}\,(\Xi^{-1}|{L_\lambda}).\\
\end{eqnarray*}
Let $[\Xi]\in\Gamma hh_p^e\Gamma/\sim_\Gamma$. Since $\lambda-\lambda'\in(p-1)p^m X({\bf T})$ we know that $\lambda'(h_{(p)})\lambda(h_{(p)})^{-1}\in 1+p^{m+1}{\Bbb Z}_p$
%
%$\lambda-\lambda'=(p-1)p^m\chi$, $\chi\in X({\bf T})$. Da $h_{(p)}\in{\bf G}({\Bbb Q})$ folgt $\chi(h_{(p))\in{\Bbb Q}$. Da $h_{(p)}$ keine 
%eintrage mit nicht trivialer $v_p$-bewertung hat folgt $v_p(\chi(h_{(p)))=0$, d.h. $\chi(h_{(p))\in{\Bbb Z}_p^*$. dann folgt die Behauptung aus 
%dem kleinen fermat.
%
and (7.4) Proposition implies (note that $e\ge L\ge 1$)
$$
\lambda(h_p^{e+f})\,b_e\,{\rm tr}\,(\Xi^{-1}|{L_\lambda}) \equiv \lambda'(h_p^{e+f})\,b_e\,{\rm tr}\,(\Xi^{-1}|{L_{\lambda'}})\pmod{p^{\S}},\leqno(4)
$$
where 
$$
\S={\rm min}(m+1,Ce)-e\beta
$$ 
(note that $v_p(b_e)\ge-e\beta$). Since $e\ge L\ge (m+1)/C$ we obtain $Ce\ge m+1$. Hence, $\S=m+1-e\beta$. Recalling that 
$e\le tL$, $t\le M(\beta)$ and $L\le \frac{m+1}{C}+1$ we further obtain 
$$
\S
%=m+1-\beta e\ge 
%(m+1)-\beta tL
%\ge m+1-\beta \lceil\frac{m+1}{C}\rceil 2M(\beta+1)
%\ge m+1-\beta 2M(\beta+1)(\frac{m+1}{C}+1) 
%=
\ge(1-\frac{\beta M(\beta)}{C})(m+1)-\beta M(\beta).
$$
Thus, equation (3) (which also holds with $\lambda$ replaced by $\lambda'$) and equation (4) yield the claim. This completes the proof of the 
Proposition.

\bigskip

{\bf (7.6)  $p$-Adic families of Siegel eigen classes. } We select a slope $\beta\in{\Bbb Q}_{\ge 0}$.
The above Theorem immediately implies that the family of ${\cal H}$-modules 
$({\bf H}_\lambda)=(H^d(\Gamma\backslash\bar{X},{\cal L}_\lambda))$ satisfies equation ($\dag$) in (3.7), where $a'=\frac{1}{C}$,  
%(note that $L=\lceil\frac{m+1}{C}\rceil\ge \frac{m+1}{C}$), 
$a=1-\frac{\beta M(\beta)}{C}$, $b=-\beta M(\beta)$. Thus, (3.7) Proposition yields the existence of $p$-adic continuous families of 
finite slope with (to simplify, we may omit a factor "$2$")
$$
{\sf a}=\frac{1}{M(\beta)}{\rm min}\,(a,\frac{a'}{M(\beta+1)}) \quad\mbox{and}\quad {\sf b}
%=\frac{b}{M(\beta)}-(M(\beta)+2)\log_p (M(\beta))
=-\beta-(M(\beta)+2)\log_p (M(\beta)).\leqno(5)
$$
We want to choose $C>0$ such that ${\sf a}$ becomes large. To this end we set
$$
C=C(\beta)=\beta M(\beta)+\frac{1}{M(\beta+1)}
$$
and obtain
$$
{\sf a}=\frac{1}{(1+\beta M(\beta+1)M(\beta))M(\beta)}.\leqno(6)
$$
Since with the above choice of $C(\beta)$ the numbers $a'/M(\beta+1)$, $a$ and $b$ are decreasing in $\beta$ the assumptions of (3.7) 
Proposition are satisfied. We denote by ${\cal E}(\lambda)^\beta={\cal E}({\bf H}_\lambda^\alpha)$ the set of all characters 
$\Theta:\,{\cal H}\rightarrow\bar{\Bbb Q}_p$ which are defined over ${\cal O}$ and such that the corresponding eigenspace 
$H^d(\Gamma\backslash X,{\cal L}_\lambda)^\beta(\Theta)$ w.r.t the normalized action of ${\cal H}$ does not vanish. We then obtain from (3.7) 
the following

\bigskip

%{\it Remark. } Obviously, ${\sf a}'>0$ and with the above choice of $C$ we also obtain ${\sf a}>0$. To see that also $C\ge 1$ we note that for any $\alpha\in{\Bbb Q}_{\ge 0}$ we have
%$\alpha\, {\rm dim}\,{\bf H}_\lambda^\alpha\in{\Bbb N}_0$ (cf. [G-M], p. 797). In particular, if $\alpha>0$ 
%is a non trivial slope for ${\bf H}_\lambda$ then $\alpha\,{\rm dim}\, {\bf H}_\lambda^\alpha\ge 1$. 
%Thus, if there is a slope $\alpha\in{\Bbb Q}_{\ge 0}$ such that $0<\alpha\le \beta$ and ${\bf H}_\lambda^\alpha\not=0$ for at least one $\lambda$ then 
%$\beta\,{\rm dim}\,{\bf H}_\lambda^{\le\beta+1}\ge 1$ and, hence, $\beta M(\beta+1)\ge 1$. This implies $C(\beta)\ge 1$ in case $\beta>0$. 
%If there is no such slope $\alpha$ then ${\bf H}_\lambda^{\le\beta}={\bf H}_\lambda^{\rm ord}$ for all weights $\lambda$ and we are in the ordinary case, 
%i.e. we can choose $C=1$.

{\bf (7.7) Corollary. }{\it Let $\beta\in{\Bbb Q}_{\ge 0}$.  

\medskip

1.) ${\rm dim}\, H^d(\Gamma\backslash\bar{X},{\cal L}_\lambda)^\beta$ is locally constant, i.e. there is $D=D(\beta)$ only depending on $\beta$ such 
that $\lambda\equiv\lambda'\pmod{(p-1)p^D X({\bf T})}$ implies 
$$
{\rm dim}\, H^d(\Gamma\backslash\bar{X},{\cal L}_\lambda)^\beta={\rm dim}\, H^d(\Gamma\backslash\bar{X},{\cal L}_{\lambda'})^\beta
$$

\medskip

2.) Any $\Theta\in{\cal E}(\lambda_0)^{\beta}$ fits in a p-adic continuous family of eigencharacters of slope $\beta$, i.e. there are 
$\Theta_\lambda\in{\cal E}(\lambda)^{\beta}$ such that $\Theta_{\lambda_0}=\Theta$ and $\lambda\equiv\lambda'\pmod{(p-1)p^m X({\bf T})}$ implies 
$$
\Theta_\lambda\equiv\Theta_{\lambda'}\pmod{p^{{\sf a}(m+1)+{\sf b}}};
$$ 
here, ${\sf a}$ and ${\sf b}$ as in equation (5) and (6) and $\lambda$ runs 
over all dominant characters satisfying $\langle\lambda,\alpha\rangle>2\beta M(\beta)+2$ for all simple roots $\alpha$.

}

\bigskip

{\it Remark. } The congruence between two eigencharacters $\Theta_\lambda$ and $\Theta_{\lambda'}$ is non trivial only if ${\sf a}(m+1)+{\sf b}>0$, i.e.
only if
$$
m+1>E(\beta)=(-b+M(\beta)(M(\beta)+2)\log_p M(\beta))(1+\beta M(\beta+1)M(\beta)).
$$
%beachte: hier brauchen wir, dass ${\rm min}\,(\frac{2 a'}{M(\alpha+1)},a)> 0$ ist (insbesondere also $\not=0$), da sonst beim muliplizieren mit 
%diesem minimum die vorzeichen umgekehrt wurden 
Thus, only the existence of the family $(\Theta_\lambda)$ with $\lambda\equiv \lambda_0\pmod{p^{E(\beta)}}$ is a non trivial statement. 
%
%Outside the disc $\lambda\equiv \lambda_0\pmod{p^{E(\alpha)}}$ i.e. if $\lambda_1\not\equiv \lambda_0\pmod{E(\alpha)}$ there are no congruences 
%between elements in the families $(\Theta_\lambda)$, $\lambda\equiv \lambda_0\pmod{p^{E(\alpha)}}$ and $(\Theta_\lambda)$, 
%$\lambda\equiv \lambda_1\pmod{p^{E(\alpha)}}$ i.e. the families are totally independent of each other.
%
Note that $M(\beta)\rightarrow\infty$ if $\beta\rightarrow\infty$ and that we may assume $M(\beta)\le M(\beta+1)$; these two statements imply that 
$E(\beta)\sim M(\beta+1)^4\log_p M(\beta+1)$ for $\beta\rightarrow \infty$.

\bigskip

{\bf \Large References}

\bigskip

[A-S] Ash, A., Stevens, G., $p$-adic Deformations of arithmetic cohomology, preprint 2008

[B] Bewersdorff, J., Eine Lefschetzsche Fixpunktformel f{\"u}r Hecke Operatoren, Bonner Mathematische Schriften {\bf 164} 1985 

%[G-M] Gouvea, F., Mazur, B., Families of modular eigenforms. Math Comp. {\bf 58} (1992), 793 - 805 

%[J-L] Jacquet, H., Langlands, R. P., Automorphic Forms on ${\rm GL}(2)$, LNM 114, Springer, 1970

[H 1] Hida, H., Iwasawa modules attached to congruences of cusp forms, Ann. Sci. Ecole Norm. Sup. {\bf 19} (1986) 231 - 273

[H 2] -, On $p$-adic Hecke algebras for ${\rm GL}_2$ over totally real fields, Ann. of. Math. {\bf 128} (1988) 295 - 384

%Elementary theory of Eisenstein series and $L$-functions. London Mathematical Society (1993)

[K] Koike, M., On some $p$-adic properties of the Eichler Selberg trace formula I, II, Nagoya Math. J. {\bf 56} (1974), 45 - 52, {\bf 64} (1976), 87 - 96

[M] Miyake, T., Modular forms, Springer 1989

[Ma 1] Mahnkopf, J., On Truncation of irreducible representations of Chevalley groups, J. of Number Theory, 2013

[Ma 2] -, $P$-adic families of modular forms, in: Oberwolfach Reports, 2011

[Ma 3] -, On $p$-adic families of modular forms, in: Proceedings of RIMS Conference on ``Automorphic forms, automorphic representations and related topics", 
93 - 108, Kyoto University, 2010 

%[R] Rossmann, W., Lie Groups, Oxford University Press 2002

\end{document}